\documentclass[preprint,review,12pt]{elsarticle}

\usepackage{amssymb}
\usepackage{amsmath}
\usepackage{mathrsfs}
\usepackage{mathtools}

\usepackage{lineno}

\usepackage{multirow}
\usepackage[margin=1in]{geometry}
\usepackage{hyperref}

\usepackage{subfig}

\numberwithin{equation}{section}
\newcommand{\mbf}[1]{\mathbf{#1}}
\newcommand{\mbfit}[1]{\boldsymbol{#1}}
\newcommand{\vfunt}[2]{\mathbf{#1}(\mathbf{#2},t)}
\newcommand{\vfun}[2]{\mathbf{#1}(\mathbf{#2})}
\newcommand*\diff{\mathop{}\!\mathrm{d}}
\newcommand*\vol{\mathop{}\!\mathrm{Vol}}
\newcommand{\defeq}{\vcentcolon=}
\DeclareMathOperator{\nbor}{nbor}

\newcommand{\ddiv}[1]{\mbf{D}^{#1} \cdot }
\newcommand{\dcurl}[1]{\mbf{D}^{#1} \times}
\newcommand{\dgrad}[1]{\mbf{G}^{#1}}
\newcommand{\grad}{\nabla}
\let\div\undefined
\newcommand{\div}{\nabla \cdot}
\newcommand{\lapl}{\nabla^2}
\newcommand{\dlapl}[1]{L^{#1}}
\newcommand{\vlapl}[1]{\mbf{L}^{#1}}
\newcommand{\adv}{\mbfit{u}\cdot \nabla}
\newcommand{\curl}{\nabla \times}

\newtheorem{thm}{Theorem}
\newtheorem{lem}[thm]{Lemma}
\newproof{pf}{Proof}

\newcommand{\cellgrid}{\mathbb{C}}
\newcommand{\facegrid}{\mathbb{F}}
\newcommand{\edgegrid}{\mathbb{E}}
\newcommand{\nodegrid}{\mathbb{N}}

\newcommand{\deleted}[1]{}

\newcommand{\ie}{\textit{i.e.}}

\usepackage{color}
\PassOptionsToPackage{normalem}{ulem}
\usepackage{ulem}

\definecolor{mygreen}{rgb}{0,.6,0}

\allowdisplaybreaks

\newcommand{\newsix}{\phi^{\text{new}}_{6h}}
\newcommand{\newfive}{\phi^{\text{new}}_{5h}}
\newcommand{\bspline}{\phi^{\text{B}}_{4h}}
\newcommand{\bsplinesix}{\phi^{\text{B}}_{6h}}
\newcommand{\stndfour}{\phi_{4h}}
\newcommand{\cosfour}{\phi^{\cos}_{4h}}

\usepackage[nameinlink,noabbrev]{cleveref}
\Crefname{equation}{eq.}{eqs.} 
\Crefname{equation}{Eq.}{Eqs.}

\journal{J. Comput. Phys.}

\begin{document}

\begin{frontmatter}



\title{An Immersed Boundary Method with Divergence-Free Velocity Interpolation {and Force Spreading}}


\author{Yuanxun Bao\corref{cor1}}
\ead{billbao@cims.nyu.edu}

\author{Aleksandar Donev}
\ead{donev@courant.nyu.edu}

\author{David M.~McQueen}
\ead{mcqueen@cims.nyu.edu}

\author{Charles S.~Peskin}
\ead{peskin@cims.nyu.edu}

\address{Courant Institute of Mathematical Sciences, New York University, 251 Mercer Street, New York, NY, USA}

\author{Boyce E.~Griffith}
\ead{boyceg@unc.edu}
\address{Departments of Mathematics and Biomedical Engineering, Carolina Center for Interdisciplinary Applied Mathematics, and McAllister Heart Institute, University of North Carolina, Chapel Hill, NC, USA}

\cortext[cor1]{Corresponding author}


\begin{abstract}
The Immersed Boundary (IB) method is a mathematical framework for constructing robust numerical methods to study fluid-structure interaction in problems involving an elastic structure immersed in a viscous fluid.
The IB formulation uses an Eulerian representation of the fluid and a Lagrangian representation of the structure. The Lagrangian and Eulerian frames are coupled by integral transforms with delta function kernels. The discretized IB equations use approximations to these transforms with regularized delta function kernels to interpolate the fluid velocity to the structure, and to spread structural forces to the fluid. It is well-known that the conventional IB method can suffer from poor volume conservation since the interpolated Lagrangian velocity field  is not generally divergence-free, and so this can cause spurious volume changes. In practice, the lack of volume conservation is especially pronounced for cases where there are large pressure differences across thin structural boundaries. The aim of this paper is to greatly reduce the volume error of the IB method by introducing velocity-interpolation and force-spreading schemes with the properties that the interpolated velocity field in which the structure moves is at least $\mathscr{C}^1$ and satisfies a continuous divergence-free condition, and that the force-spreading operator is the adjoint of the velocity-interpolation operator. We confirm through numerical experiments in two and three spatial dimensions that this new IB method   is able to achieve substantial improvement in volume conservation compared to other existing IB methods, at the expense of a modest increase in the computational cost. {Further, the new method provides smoother Lagrangian forces (tractions) than traditional IB methods. The method presented here is restricted to periodic computational domains. Its generalization to non-periodic domains is important future work. }
\end{abstract}


\end{frontmatter}



\section{Introduction}
The \textit{Immersed Boundary} (IB) method \cite{Peskin2003_IBreview} is a general mathematical framework for the numerical solution of \textit{fluid-structure interaction} problems arising in biological and engineering applications. The IB method was introduced to simulate flow patterns around the heart valves \cite{Peskin1972,Peskin1977_IBfirst}, and since its success in modeling cardiac fluid dynamics \cite{McQueen1997, IBAMR,IBMDelta_Boundary,IBAMR_HeartValve}, it has been extended and applied to various other applications, including but not limited to motion of biological swimmers \cite{Bhalla2013, Lushi2013}, dynamics of red-blood cells \cite{Fai2013} and dry foam \cite{Kim2010_2Dfoam, Kim2014_foam3D}, and rigid body motion \cite{RigidIBM, RigidMultiblobs}.

The essence of the IB method as a numerical scheme lies in its simple way of coupling an Eulerian representation of the fluid and a Lagrangian representation of the structure. {The \textit{force spreading} linear operator {$\mbfit{S}$} that spreads forces (stresses)} from the structure to the fluid and {the \textit{velocity interpolation} linear operator {$\mbfit{S}^{\star}$} that interpolates velocities from the fluid to the structure} are carried out via a regularized delta function $\delta_h$.
One effective way to construct $\delta_h$ is to require the regularized delta function to satisfy a set of moment conditions to achieve approximate grid translation-invariance and desired interpolation accuracy \cite{C3IBkernels,Bao2016_new6pt}, thereby avoiding special grid treatment near the fluid-structure interface. In spite of its wide applicability and ease of implementation, the conventional IB method with a collocated-grid discretization (referred to herein as IBCollocated) has two well-known shortcomings in accuracy: it achieves only first-order convergence for problems that possess sharp-interface solutions \cite{Lai2000_IBaccuracy,Griffith2005_IBaccuracy}, and it can be relatively poor of volume conservation \cite{Peskin1993_IBmodified}. Much research effort has been put into improving the convergence rate of the IB method to second order or even higher order for problems with singular forcing at the sharp interface. Notable examples include, the \textit{Immersed Interface Method} (IIM) \cite{Li2001, Lee2003_IIM}, and more recently, a new method known as \textit{Immersed Boundary Smooth Extension} \cite{IBSE,IBSEflows}. Our focus here, however, is on improving the volume conservation properties of the IB method.  


As an immediate consequence of fluid incompressibility, which is one of the basic assumptions of the IB formulation, the volume enclosed by the immersed structure is exactly conserved as it deforms and moves with the fluid {in the continuum setting}. Thus, a desirable feature of an IB method is to conserve volume as nearly as possible.
In practice, however, it is observed that, even in the simplest case of a quasi-static pressurized membrane \cite{Griffith2012_IBMACvolume}, the conventional IB method (regardless of collocated- or staggered-grid discretization) produces volume error that persistently grows in time, as if fluid ``leaks'' through the boundary. 
{An intuitive explanation for this ``leak'' is that fluid is ``squeezing'' between the marker points used to discretize the boundary in a conventional IB method; however, this is {\it not} the full story, because refining the Lagrangian discretization does not improve the volume conservation of the method for a fixed Eulerian discretization.}

{In the conventional IB method, we can extend the notion of velocity interpolation to any point in the domain (not restricted to the immersed structure), denoted here with an italic $\mbfit{X}$. The {\it continuous} interpolated velocity field can be written as $\mbfit{U}(\mbfit{X}) = (\mbfit{\mathcal{J}} \mbf{u})(\mbfit{X})$, where $\mbfit{\mathcal{J}}$ denotes the {\it continuous} interpolation operator that interpolates the velocity at $\mbfit{X}$ from the discrete fluid velocity $\mbf{u}$. If a closed surface moves with velocity that is {\it continuously} divergence-free with respect to the {\it continuum} divergence operator, \ie, ($\div \mbfit{U})(\mbfit{X}) = 0$, then the volume enclosed by the (deformed) surface will be exactly conserved. However, in the discrete setting, even if the interpolated velocity field is {\it continuously} divergence-free, exact volume conservation is generally not achieved because of the time-stepping error from the temporal integrator. Another source of error comes from discretizing the surface itself. In the IB method, only a discrete collection of points on the surface, \ie, the Lagrangian markers, move according to the interpolated velocity field. A closed discretized surface can be constructed by simply connecting the Lagrangian markers defining a facet, and the resulting faceted surface by this construction does not enclose a constant volume. In the absence of temporal integration errors, this kind of volume-conservation error will approach zero as the discretization of the surface is refined. Peskin and Printz realized that the major cause of poor volume conservation of IBCollocated is that the {\it continuous} interpolated velocity field given by the conventional IB interpolation operator (denoted by $\mbfit{\mathcal{J}}_{\text{IB}}$) is not {\it continuously} divergence-free \cite{Peskin1993_IBmodified}, despite that the discrete fluid velocity is enforced to be {\it discretely} divergence-free with respect to the {\it discrete} divergence operator by the fluid solver.}

To improve the volume conservation of the conventional IB method, Peskin and Printz proposed a modified finite-difference approximation to the discrete divergence operator to ensure that the \textit{average} of the continuous divergence of the interpolated velocity is equal to zero in a small control volume with size of a grid cell \cite{Peskin1993_IBmodified}. Their IB method with modified finite-difference operators (herein referred to as IBModified) was applied to a two-dimensional model of the heart, and it achieved improvement in volume conservation by one-to-two orders of magnitude compared to IBCollocated. Nevertheless, a major drawback of IBModified that limits its use in applications is its complex, non-standard finite-difference operators that uses coefficients derived from the regularized delta function (but see \cite{Kim2010_2Dfoam,Kim2014_foam3D} for applications). To address the issue of  spurious currents across immersed structure supporting extremely large pressure differences, Guy and Strychalski \cite{Strychalski2016} developed a different extension of the IB method that uses non-uniform Fast Fourier Transform \cite{Dutt1993,NUFFT} (NUFFT) to generate ``spectral'' approximations to the delta function, which also has superior volume conservation.

Over the past two decades, the staggered-grid (MAC) discretization has been widely adopted by the IB community \cite{IBAMR, IBAMR_HeartValve, Bhalla2013, RigidIBM, RigidMultiblobs, Devendran2012_IBviscoelastic}. In addition to its most celebrated feature of avoiding the odd-even decoupling in the Poisson solver that can otherwise occur with collocated-grid discretization, which leads to ``checkerboard'' instability in the solutions, Griffith \cite{Griffith2012_IBMACvolume} concluded from his numerical studies that the improvement in volume conservation of the IB method with staggered-grid discretization (IBMAC) is essentially the same as that of IBModified. In practice, IBMAC is more practical than IBModified in that the improvement in volume conservation directly comes as a byproduct of grid discretization without any modification to the finite-difference operators, and it is relatively straightforward to extend IBMAC to include adaptive mesh refinement \cite{IBAMR,Roma1999} and physical boundary conditions \cite{NonProjection_Griffith}. However, we emphasize that the nature of Lagrangian velocity interpolation of IBMAC remains the same as that of IBCollocated, and, hence, there is much room for further improvement in volume conservation by ensuring that the interpolated velocity is constructed to be nearly or exactly divergence-free. We note that the methods designed to improve the convergence rate of IB methods, such as IIM \cite{Li2001,Lee2003_IIM} and the Blob-Projection method \cite{Cortez2000_blobprojection}, also improve volume conservation, because the solution near the interface is computed more accurately. These methods, however, are somewhat more complex and less generalizable than the conventional IB method.

This paper is concerned with further improving volume conservation of IBMAC by constructing a {{\it continuous}} velocity-interpolation operator {$\mbfit{\mathcal{J}}$} that is divergence-free in the \textit{continuous} sense. {The discrete IB interpolation operator $\mbfit{S}^{\star}$ is simply the restriction of $\mbfit{\mathcal{J}}$ to the Lagrangian markers.}
The key idea introduced in this paper is first to construct a \textit{discrete vector potential} that lives on an \textit{edge-centered} staggered grid from the discretely divergence-free fluid velocity, and then to apply the conventional IB interpolation scheme to obtain a \textit{continuum} vector potential, from which the interpolated velocity field is obtained by applying the continuum curl operator.  Note that the existence of the discrete vector potential relies on the fact that the discrete velocity field is discretely divergence-free. The interpolated velocity field obtained in this manner is guaranteed to be continuously divergence-free, since the divergence of the curl of any vector field is zero.
{We also propose a new force-spreading operator {$\mbfit{S}$} that is defined to be the new adjoint of the interpolation operator $\mbfit{S}^{\star}$, so that Lagrangian-Eulerian interaction conserves energy.} The Eulerian force density that is the result of applying this force-spreading operator to a Lagrangian force field turns out to be discretely divergence-free, so we refer to this new force-spreading operation as divergence-free force spreading. We name the IB method equipped with the new interpolation and spreading operators as the \textit{Divergence-Free Immersed Boundary} (DFIB) method. {As presented here, the DFIB method is limited to periodic domains.}

In contrast to the local nature of interpolation and spreading in the conventional IB method, the spreading and interpolation operators of the DFIB method turn out to be non-local in that their construction requires the solution of discrete Poisson equations, although these operators can be evaluated efficiently using the Fast Fourier Transform (FFT) 
or multigrid methods.
Another new feature of our method is that transferring information between the Eulerian grid and the Lagrangian mesh involves derivatives of the regularized delta function $\grad \delta_h$ instead of only $\delta_h$. We confirm through various numerical tests in both two and three spatial dimensions that the DFIB method is able to reduce volume error by several orders of magnitude compared to IBMAC and IBModified at the expense of only a modest increase in the computational cost. Moreover, we confirm that the volume error for DFIB decreases as the Lagrangian mesh is refined with the Eulerian grid size held fixed, which is not the case in the conventional IB method \cite{Peskin1993_IBmodified}.
In addition to the substantial improvement in volume conservation, the DFIB method is quite straightforward to realize from an existing modular IB code with staggered-grid discretization, that is, by simply switching to the new velocity-interpolation and force-spreading schemes while leaving the fluid solver and time-stepping scheme unchanged. 

The rest of the paper is organized as follows. In \autoref{sec_eqnmotion}, we begin by giving a brief description of the continuum equations of motion in the IB framework. Then we define the staggered grid on which the fluid variables live and introduce the spatial discretization of the equations of motion.  \autoref{sec_divfreeinterp} introduces the two main contributions of this paper: divergence-free velocity interpolation and force spreading. In \autoref{sec_timestepping}, we present a formally second-order time-stepping scheme that is used to evolve the spatially-discretized equations, followed by a cost comparison of DFIB and IBMAC. Numerical examples of applying DFIB to problems in two and three spatial dimensions are presented in \autoref{sec_numericalresults}, where the volume-conserving characteristics of the new scheme are assessed.



\section{Equations of motion and spatial discretization} \label{sec_eqnmotion}
\subsection{Equations of motion}
This section provides a brief description of the continuum equations of motion in the IB framework \cite{Peskin2003_IBreview}. We assume a neutrally-buoyant elastic structure $\Gamma$ that is described by the Lagrangian variables $\mbfit{s}$, immersed in a viscous incompressible fluid occupying the whole fluid domain $\Omega \subset \mathbb{R}^3$ that is described by the Eulerian variables $\mbfit{x}$. \Cref{momentum_eq,divfree_eq} are the incompressible Navier-Stokes equations describing mass and momentum conservation of the fluid, in which $\mbfit{u}(\mbfit{x},t)$ denotes the fluid velocity, $p(\mbfit{x},t)$ is the pressure, and $\mbfit{f}(\mbfit{x},t)$ is the Eulerian force density (force per unit volume) exerted by the structure on the fluid. In this formulation, we assume that the density $\rho$ and the viscosity $\mu$ of the fluid are constant. The fluid-structure coupled equations are:
\begin{align}
& \rho\left( \frac{\partial \mbfit{u} }{\partial t} + \mbfit{u}\cdot \nabla \mbfit{u}  \right) + \grad p = \mu  \lapl \mbfit{u} + \mbfit{f}, \label{momentum_eq} \\
& \div \mbfit{u} =0, \label{divfree_eq} \\
& \vfunt{\mbfit{f}}{\mbfit{x}}  = \int_{\Gamma} \vfunt{\mbfit{F}}{\mbfit{s}}\, \delta( {\mbfit{x} -\vfunt{ {\mbfit{\mathcal{X}}}}{\mbfit{s}}} ) \diff \mbfit{s}, \label{forcespreading_eq} \\
& \frac{\partial { \mbfit{\mathcal{X}}} }{\partial t}(\mbfit{s},t) = \mbfit{u}({\mbfit{\mathcal{X}}}(\mbfit{s},t),t) = \int_{\Omega}  \mbfit{u}(\mbfit{x},t) \, \delta(\mbfit{x} - {\mbfit{\mathcal{X}}}(\mbfit{s},t)) \diff \mbfit{x}, \label{velocityinterp_eq}\\
& \mbfit{F}(\mbfit{s},t) = \boldsymbol{\mathcal{F}}[\vfunt{{\mbfit{\mathcal{X}}}}{\cdot} \, ; \mbfit{s}] = -\frac{\delta E}{\delta {\mbfit{\mathcal{X}}}}(\mbfit{s},t).
\label{Lagrangian_force_eq}
\end{align}
\Cref{forcespreading_eq,velocityinterp_eq} are the fluid-structure interaction equations that couple the Eulerian and the Lagrangian variables. \Cref{forcespreading_eq} relates the Lagrangian force density $\mbfit{F}(\mbfit{s},t)$ to the Eulerian force density $\mbfit{f}(\mbfit{x},t)$ using the Dirac delta function, where ${\mbfit{\mathcal{X}}}(\mbfit{s},t)$ is the physical position of the Lagrangian point $\mbfit{s}$.  \Cref{velocityinterp_eq} is simply the no-slip boundary condition of the Lagrangian structure, \ie, the Lagrangian point ${\mbfit{\mathcal{X}}}(\mbfit{s},t)$ moves at the same velocity as the fluid at that point. In \Cref{Lagrangian_force_eq}, the system is closed by expressing the Lagrangian force density $\mbfit{F}(\mbfit{s},t)$ in the form of a force density functional $\mbfit{\mathcal{F}}[ {\mbfit{\mathcal{X}}}(\cdot,t) \, ; \mbfit{s}]$, which in many cases can be derived from an elastic energy functional $E[{\mbfit{\mathcal{X}}}(\cdot,t) \,; \mbfit{s}]$  by taking  the variational derivative, denoted here by $\delta/\delta {\mbfit{\mathcal{X}}}$, of the elastic energy. 

\subsection{Spatial discretization}

\begin{figure}
\centering
\subfloat[][]{\includegraphics[width = .25\linewidth]{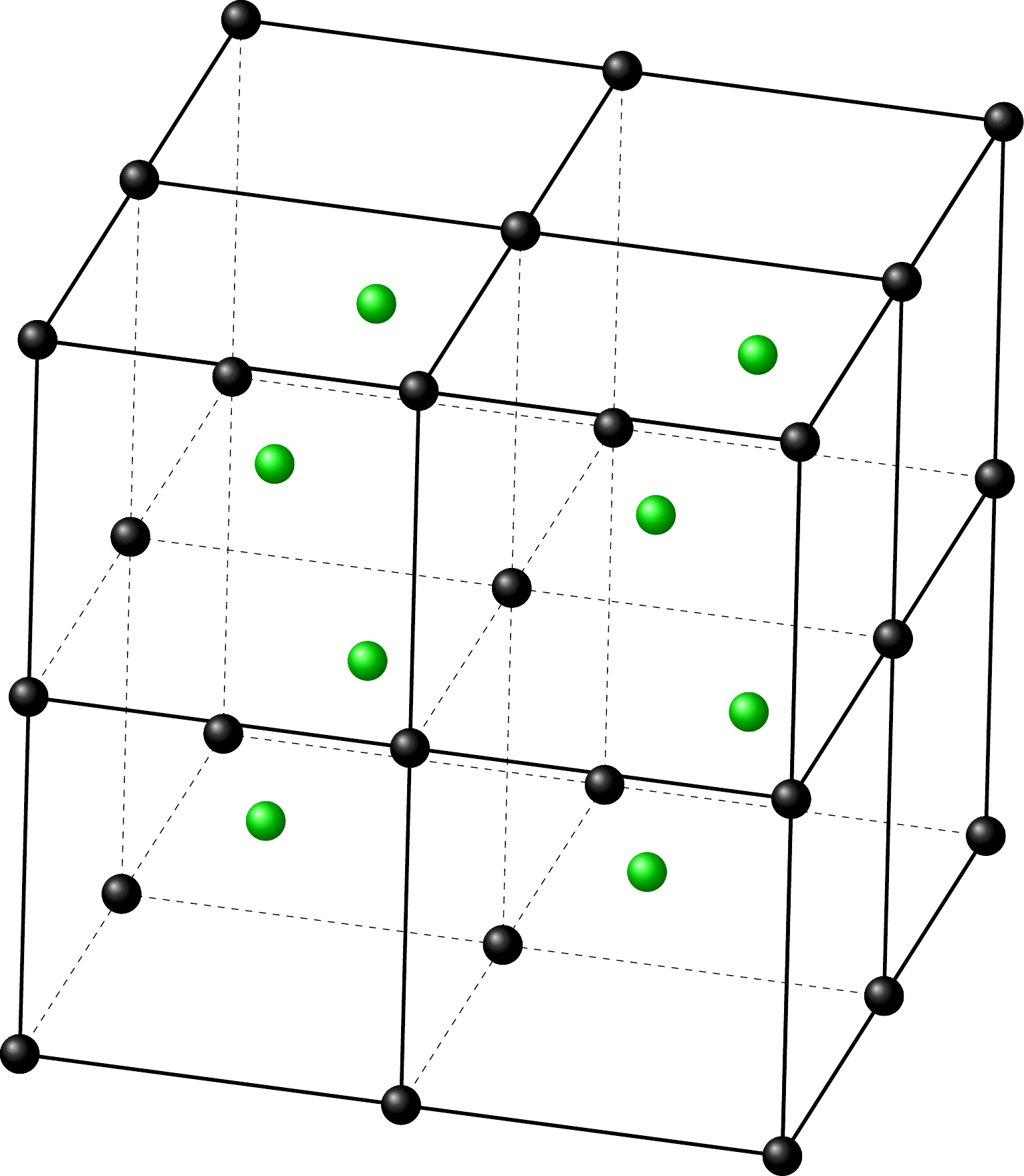} \label{fig:cell_node_grid}} \hspace{2em}
\subfloat[][]{\includegraphics[width = .25\linewidth]{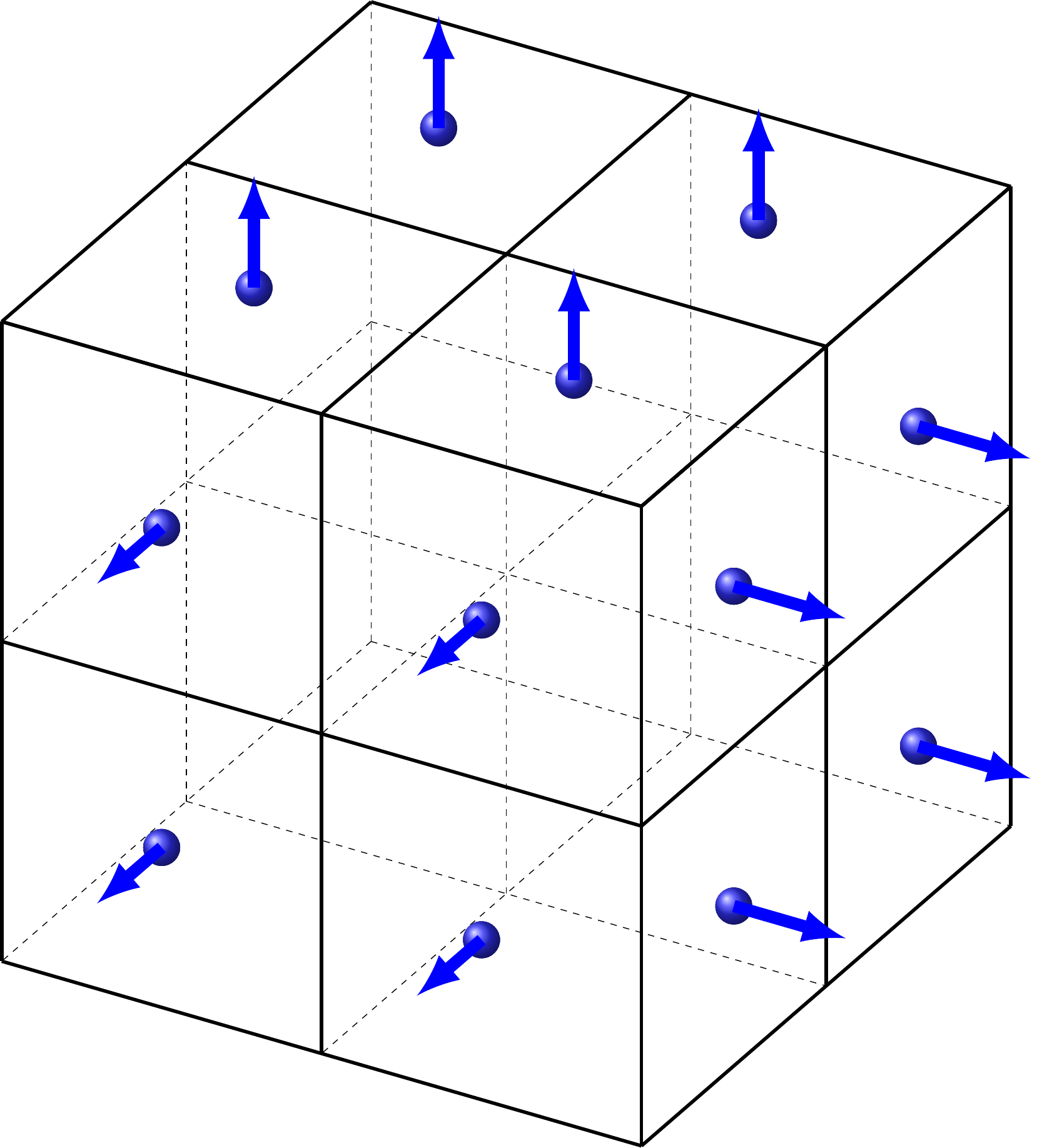} \label{fig:face_grid}} \hspace{2em}
\subfloat[][]{\includegraphics[width = .25\linewidth]{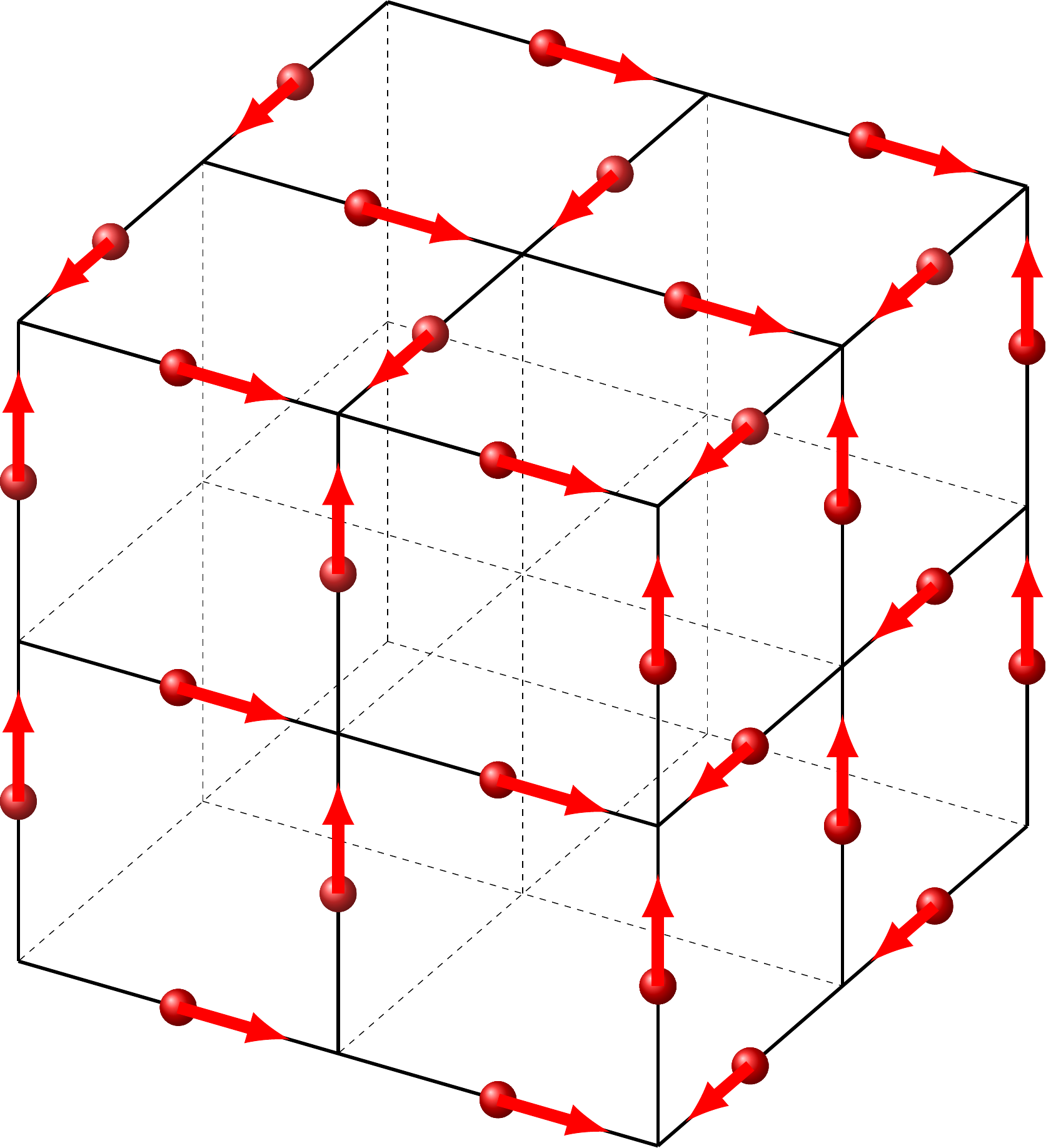} \label{fig:edge_grid}
}
\caption{Staggered grids on which discrete grid functions are defined. (a) Cell-centered (green) and node-centered (black) grids for scalar functions. (b) Face-centered grid for vector grid functions. (c) Edge-centered grid for vector grid functions.}
\end{figure}

Throughout the paper, we assume the fluid occupies a {\it periodic} domain $\Omega = [0,L]^3$ that is discretized by a uniform $N \times N \times N$ Cartesian grid with meshwidth $h = \frac{L}{N}$. Each grid cell is indexed by $(i,j,k)$ for $i, j , k = 0,\dots, N-1$. 
For the Eulerian fluid equations, we use the staggered-grid discretization, in which the pressure $p$ is defined on the cell-centered grid (\autoref{fig:cell_node_grid}), denoted by $\cellgrid$, \ie, at positions $\mbf{x}_{i,j,k} = ( (i+\frac{1}{2})h, (j+\frac{1}{2})h,  (k+\frac{1}{2})h)$. The discrete fluid velocity $\mbf{u}$ is defined on the face-centered grid (\autoref{fig:face_grid}), denoted by $\facegrid$, with each component perpendicular to the corresponding cell faces, \ie, at positions $\mbf{x}_{i-\frac{1}{2},\, j,\, k}$\,, $\mbf{x}_{i,\, j-\frac{1}{2},\, k}$ and $\mbf{x}_{i,\, j, \, k-\frac{1}{2}}$ for each velocity component respectively. 
We also introduce two additional shifted grids: the node-centered grid (\autoref{fig:cell_node_grid}) for scalar grid functions, denoted by $\nodegrid$ , \ie, at positions $\mbf{x}_{i-\frac{1}{2},\, j-\frac{1}{2}, \, k-\frac{1}{2}}$, and the edge-centered grid (\autoref{fig:edge_grid}) for vector grid functions, denoted by $\edgegrid$, with each component defined to be parallel to the corresponding cell edges \ie, at positions $\mbf{x}_{i,\, j-\frac{1}{2}, \, k-\frac{1}{2}}$,  $\mbf{x}_{i-\frac{1}{2}, \, j, \, k-\frac{1}{2}}$ and $\mbf{x}_{i-\frac{1}{2}, \, j-\frac{1}{2}, \, k}$ for each component respectively. In \autoref{sec_divfreeinterp}, we will use these half-shifted staggered grids to construct divergence-free velocity interpolation and force spreading.

To discretize the differential operators in \Cref{momentum_eq,divfree_eq}, we introduce the central difference operators corresponding to the partial derivatives $\partial / \partial x_{\alpha}$,
\begin{equation}
 D_{\alpha}^h \varphi \defeq \frac{\varphi(\mbf{x} + \frac{h}{2}\mbf{e}_{\alpha}) - \varphi(\mbf{x} - \frac{h}{2}\mbf{e}_{\alpha})}{h}, \quad \alpha = 1, 2, 3, 
\end{equation}
where $\varphi$ is a scalar grid function and $\{ \mbf{e}_1, \mbf{e}_2, \mbf{e}_3 \}$ is the standard basis of $\mathbb{R}^3$. We can use $D^h_{\alpha}$ to define the discrete gradient, divergence and curl operators:
\begin{align}
 &\dgrad{h} \varphi \defeq (D^h_1\varphi,\, D^h_2\varphi,\, D^h_3\varphi), \label{dgrad}\\
 &\ddiv{h} \mbf{v} \defeq  D^h_{\alpha}v_{\alpha}, \label{ddiv} \\
 &\dcurl{h} \mbf{v}  \defeq \epsilon_{ijk} D^h_j v_k \label{dcurl}, 
\end{align}
where $\mbf{v}$ is a vector grid function$, \epsilon_{ijk}$ is the totally antisymmetric tensor, and the Einstein summation convention is used here. The discrete differential operators may be defined on different pairs of domain and range (half-shifted staggered grids), and therefore, in a slight abuse of notation, we will use the same notation to denote the different operators, 
\begin{align}
 \dgrad{h} & : \varphi(\cellgrid) \longrightarrow \mbf{v}({\facegrid})  \text{ or } \varphi(\nodegrid) \longrightarrow \mbf{v}(\edgegrid), \\
 \ddiv{h}  & : \mbf{v}(\edgegrid) \longrightarrow \varphi(\nodegrid)  \text{ or } \mbf{v}(\facegrid) \longrightarrow \varphi(\cellgrid), \\
 \dcurl{h} & : \mbf{v}(\edgegrid) \longrightarrow \mbf{v}(\facegrid)  \text{ or } \mbf{v}(\facegrid) \longrightarrow \mbf{v}(\edgegrid).
\end{align}
Although the curl operator does not appear in the equations of motion explicitly, we define it here for use in \autoref{sec_divfreeinterp}. 
The discrete scalar Laplacian operator can be defined by $L^h = \ddiv{h} \dgrad{h}$, which yields the familiar compact second-order approximation to $\lapl$: 
\begin{equation}
  L^h \varphi \defeq \sum_{\alpha = 1}^3 \frac{\varphi(\mbf{x} + h\mbf{e}_{\alpha}) - 2\varphi(\mbf{x})+\varphi(\mbf{x} - h\mbf{e}_{\alpha})}{h^2}. 
\end{equation}
Note that the range and domain of $\dlapl{h}$ are a set of grid functions defined on the same grid, and that grid can be $\cellgrid$ or $\nodegrid$ or any of the three subgrids of $\edgegrid$ or $\facegrid$ on which the different components of vector-valued functions are defined. We will use the notation $\vlapl{h}$ to denote the discrete vector Laplacian operator that applies (the appropriately shifted) $\dlapl{h}$ to each component of a vector grid function.

\subsubsection{Advection}

We follow the same treatment of discretization of the advection term as in earlier presentations of the IB method  \cite{Devendran2012_IBviscoelastic}. From the incompressibility of the fluid flow $\div \mbfit{u} = 0$, we can write the advection term  in the skew-symmetric form 
\begin{equation}
 [(\adv) \mbfit{u}]_{\alpha} = \frac{1}{2}  \mbfit{u}\cdot (\grad u_{\alpha}) + \frac{1}{2} \div (\mbfit{u} u_{\alpha}), \quad \alpha = 1,2,3. \label{adv_skewsymm} 
\end{equation} 
Let {$\mbfit{N}(\mbf{u})$} denote the discretization of \Cref{adv_skewsymm}, and we define 
\begin{equation}
 [ \mbfit{N}(\mbf{u}) ]_{\alpha} = \frac{1}{2} \tilde{\mbf{u}} \cdot \dgrad{2h} u_{\alpha} + \frac{1}{2} \ddiv{2h} (\tilde{\mbf{u}} u_{\alpha}), \quad \alpha = 1,2,3, \label{skewsymm_discretization}
\end{equation}
where $\tilde{\mbf{u}}$ denotes an averaged collocated advective velocity whose components all live on the same grid as $u_{\alpha}$. The advective velocity $\tilde{\mbf{u}}$ in \cite{Devendran2012_IBviscoelastic} is obtained by using the same interpolation scheme as the one used for moving the immersed structure. In our work, we simply take the average of $\mbf{u}$ on the grid. For example, the three components of $\tilde{\mbf{u}}$ in the $x$-component equation are 
\begin{align*}
\tilde{u}_1 &= u_1(\mbf{x}_{i-\frac{1}{2},j,k}) \,, \\
\tilde{u}_2 &= \frac{u_2(\mbf{x}_{i,j-\frac{1}{2},k}) + u_2(\mbf{x}_{i,j+\frac{1}{2},k}) + u_2(\mbf{x}_{i-1,j-\frac{1}{2},k}) + u_2(\mbf{x}_{i-1,j+\frac{1}{2},k})}{4}, \\
\tilde{u}_3 &= \frac{u_3(\mbf{x}_{i,j,k-\frac{1}{2}}) + u_3(\mbf{x}_{i,j,k+\frac{1}{2}}) + u_3(\mbf{x}_{i-1,j,k-\frac{1}{2}}) + u_3(\mbf{x}_{i-1,j,k+\frac{1}{2}})}{4} .
\end{align*} 
Note that in the $y$- and $z$-component equations, we need different averages of $\mbf{u}$ to construct $\tilde{\mbf{u}}$.
We choose to use the wide-stencil operators in \Cref{skewsymm_discretization} so that the resulting grid functions are all defined on the same grid as $u_{\alpha}$.
A more compact discretization of the advection term has been previously described in \cite{Griffith2012_IBMACvolume,NonProjection_Griffith,LLNS_Staggered}.

\subsubsection{Fluid-Structure Interaction}

{The immersed structure $\Gamma$ is discretized by a Lagrangian mesh of $M$ points or markers, {denoted here by a non-italic $\mbf{X}=\left\{ \mbf{X}_m  \right\}_{m=1}^M$, and the discrete Lagrangian force densities defined on the Lagrangian markers are $\mbf{F} = \left\{ \mbf{F}_m \right\}_{m=1}^M$}. 
As discussed in the introduction, we can extend the notion of velocity interpolation to any point $\mbfit{X}$ in the domain, not just restricted to the Lagrangian markers $\mbf{X}$, and define a {\it continuous} interpolated velocity field $\mbfit{U}(\mbfit{X}) = (\mbfit{\mathcal{J}} \mbf{u})(\mbfit{X})$. 
In the conventional IB method, the {\it continuous} velocity-interpolation operator $\mbfit{\mathcal{J}}_{\text{IB}}$ can be defined as 
\begin{equation}
 \mbfit{U}(\mbfit{X}) = (\mbfit{\mathcal{J}_{\text{IB}}} \mbf{u})(\mbfit{X})  \defeq \sum_{\mbf{x} \in \facegrid} \mbf{u}(\mbf{x}) \delta_h( \mbf{x} - \mbfit{X}) h^3.
 \label{conventionalIB_interp_continuous}
\end{equation}
}
{We note that  the interpolated velocity field given by \Cref{conventionalIB_interp_continuous} is {\it not} generally divergence-free with respect to the continuum divergence operator\footnote{The interpolated velocity  given by \Cref{conventionalIB_interp_continuous} has the same regularity as the regularized delta function $\delta_h$ which are generally at least $\mathscr{C}^1$ in the IB method, and therefore, the divergence of $\mbfit{U}$ is well-defined.},
\ie, generally
\begin{equation}
(\div \mbfit{U}) (\mbfit{X}) = -\sum_{\mbf{x} \in \facegrid} \mbf{u}(\mbf{x}) \cdot (\grad\delta_h) ( \mbf{x} - \mbfit{X}) h^3 \neq 0, \label{traditionalIBinterp_div}
\end{equation}
}
even if $\mbf{u}$ is {discretely} divergence-free with respect to the discrete divergence operator.
{The restriction of $\mbfit{\mathcal{J}}$ to the collection of Lagrangian markers $\mbf{X}$ defines the discrete IB interpolation operator
\begin{equation}
(\mbfit{S}^{\star}[\mbf{X}]\mbf{u})(\mbf{X}) = (\mbfit{\mathcal{J}} \mbf{u})(\mbf{X}). 
\end{equation}
We will also develop a new force-spreading operator $\mbfit{S}[\mbf{X}]$ that is the adjoint of the new velocity-interpolation operator $\mbfit{S}^{\star}[\mbf{X}]$. Here we use the notation $[\mbf{X}]$ to emphasize that these linear operators are parametrized by the position of the markers, as will be important when discussing temporal integration.}


The {discretization of the} interaction equations (\Cref{forcespreading_eq,velocityinterp_eq}) can be concisely written in the form
\begin{align}
\mbf{f}(\mbf{x}) &= \left(\mbfit{S}[\mbf{X}]\mbf{F} \right) (\mbf{x}), \label{force_spreading} \\
\mbf{U}(\mbf{X}) &= \left( \mbfit{S}^{\star}[\mbf{X}]\mbf{u} \right)(\mbf{X}), \label{velocity_interpolation}
\end{align}
{where $\mbf{f}$ is the discrete Eulerian force density defined on the appropriate subgrid of $\facegrid$ for each component, and $\mbf{U} = \left\{ \mbf{U}_m\right\}_{m=1}^M$ denotes the interpolated velocities at the Lagrangian markers $\mbf{X}$. }
In the conventional IB method, {the force-spreading operator $\mbfit{S}_{\text{IB}}$ and the velocity-interpolation operator $\mbfit{S}^{\star}_{\text{IB}}$} are simply discrete approximations of the surface and volume integrals in \Cref{forcespreading_eq,velocityinterp_eq}, \ie,
\begin{align}
{(\mbfit{S}_{\text{IB}}[\mbf{X}]\mbf{F} )(\mbf{x})} &\defeq \sum_{m=1}^M  \mbf{F}_m \, \delta_h( \mbf{x} - \mbf{X}_m) \Delta \mbf{s}, \label{traditionalIBforcespreading} \\
{(\mbfit{S}^{\star}_{\text{IB}}[\mbf{X}]\mbf{u})(\mbf{X})} & \defeq \sum_{\mbf{x} \in \facegrid} \mbf{u}(\mbf{x}) \delta_h( \mbf{x} - \mbf{X}) h^3, \label{traditionalIBinterp}
\end{align}
{and they are adjoint operators with respect to the power identity (inner product) defined later in \Cref{power_identity}.}
Note that \Cref{traditionalIBinterp} is a vector equation. For each of the three components of the equation, the sum $\mbf{x} \in \facegrid$ is to be understood here and in similar expressions as the sum over the appropriate subgrid of $\facegrid$.
In \Cref{traditionalIBinterp,traditionalIBforcespreading}, the Dirac delta function is replaced by a regularized delta function $\delta_h$ to facilitate the coupling between the Eulerian and Lagrangian grids,  which is taken to be of the  tensor-product form
\begin{equation}
 \delta_h(\mbf{x}) = \frac{1}{h^3} \phi\left(\frac{x_1}{h}\right) \phi\left(\frac{x_2}{h}\right) \phi\left(\frac{x_3}{h}\right),
\end{equation}
where $\phi(r)$ denotes the one-dimensional immersed-boundary kernel that is constructed from a set of moment conditions to achieve approximate grid translation invariance \cite{Peskin2003_IBreview, Bao2016_new6pt}. 

{In the following section, we will develop a new velocity-interpolation operator $\mbfit{\mathcal{J}}$ that produces a {\it continuously} divergence-free interpolated velocity field constructed from a discretely divergence-free discrete fluid velocity.}

In summary, the spatially-discretized equations of motion are 
\begin{align}
& \rho\left( \frac{\diff \mbf{u} }{\diff t} +  \mbfit{N}( \mbf{u}) \right) + \dgrad{h} p = \mu  \mbf{L}^h \mbf{u} + {\mbfit{S}[\mbf{X}]\mbf{F}}, \label{momentum_spatial} \\
& \ddiv{h} \mbf{u} =0, \label{discrete_divzero} \\
& \frac{\diff {{\mbf{X}}}}{\diff t} = \mbf{U}({\mbf{X}},t) = {\mbfit{S}^{\star}[\mbf{X}]\mbf{u}}. \label{discrete_interp} 
\end{align}

\section{Divergence-free velocity interpolation and force spreading} \label{sec_divfreeinterp}
This section presents the two main contributions of this paper: divergence-free velocity interpolation and force spreading. Familiarity with discrete differential operators on staggered grids and with some discrete vector identities, reviewed and summarized in \ref{appendix_vector_identity}, will facilitate the reading of this section.


\subsection{Divergence-free velocity interpolation} \label{subsec_divfreeinterp}

Here we introduce a new recipe for {constructing} an interpolated velocity field  {$\mbfit{U}(\mbfit{X}) = (\mbfit{\mathcal{J}}\mbf{u})(\mbfit{X})$} that is {{\it continuously}} divergence-free {with respect to the continuum divergence operator, \ie, $(\nabla \cdot \mbfit{U}) (\mbfit{X}) = 0$ for all $\mbfit{X}$}. For now we drop the dependence on time and {emphasize again that $\mbfit{X}$} is an arbitrary position in the domain $\Omega \subset \mathbb{R}^3$, not just on the Lagrangian structure $\Gamma$. The main idea is first to construct a discrete vector potential $\mbf{a}(\mbf{x})$ that is defined on the edge-centered staggered grid $\edgegrid$, and then to apply the conventional IB interpolation to $\mbf{a}(\mbf{x})$ to obtain a continuum vector potential {$\vfun{\mbfit{A}}{\mbfit{X}}$}, so that the Lagrangian velocity defined by {$\vfun{\mbfit{U}}{\mbfit{X}} = (\curl{\mbf{\mbfit{A}}})(\mbfit{X})$} is automatically divergence-free. 

Suppose the discrete velocity field $\vfun{u}{x}$ is defined on $\facegrid$ and is discretely divergence-free, \ie, $\ddiv{h} \mbf{u} = 0$. Let $\mbf{u}_0$ be the mean of $\vfun{u}{x}$,
\begin{equation}
\mbf{u}_0 = \frac{1}{V} \sum_{\mbf{x} \in \facegrid} \vfun{u}{x} h^3,
\end{equation}
where  $ V = \sum_{\mbf{x} \in \facegrid} h^3$ is the volume of the domain.
Using the Helmholtz decomposition,  we construct a discrete velocity potential $\mbf{a}(\mbf{x})$ for $\mbf{x} \in \edgegrid$ that satisfies
\begin{equation}
\left\{
\begin{array}{lcl}
	\dcurl{h} \mbf{a} &=& \mbf{u} - \mbf{u}_0 ,\\
         \ddiv{h} \mbf{a} &=& 0,
\end{array}
\right.
\label{vpot_def}
\end{equation}
where the requirement that $\mbf{a}(\mbf{x})$ is {discretely} divergence-free is an arbitrary gauge condition that makes $\mbf{a}(\mbf{x})$ uniquely defined up to a constant. {If the gauge condition of $\mbf{a}(\mbf{x})$ is omitted in \Cref{vpot_def}, then the discrete velocity potential $\mbf{a}(\mbf{x})$ is only uniquely defined up to $\dgrad{h} \psi$, where $\psi$ is some unknown scalar grid function defined on $\nodegrid$.} Note that $\ddiv{h} \mbf{a}$ is a scalar field defined on $\nodegrid$. In \ref{appendix_existence}, we prove that the discrete vector potential $\mbf{a}(\mbf{x})$ defined by \Cref{vpot_def} exists (see \autoref{appendix_vpot_existence}). To determine $\mbf{a}(\mbf{x})$ explicitly, we take the discrete curl of the first equation in \Cref{vpot_def} and use the identity \Cref{curlcurl} {with the gauge condition of $\mbf{a}({\mbf{x}})$}, which leads to a vector Poisson equation for $\mbf{a}(\mbf{x})$,
\begin{equation}
- \vlapl{h} \, \mbf{a} = \dcurl{h} \mbf{u}, \label{vpot_poisson}
\end{equation} 
that can be efficiently solved. Note that the solution of the Poisson problem \Cref{vpot_poisson} determines $\mbf{a}(\mbf{x})$ up to an arbitrary constant (it is not necessary to uniquely determine $\mbf{a}(\mbf{x})$ because the constant term vanishes {upon subsequent differentiation}). 

The next step is to interpolate the discrete vector potential $\mbf{a}(\mbf{x})$ to obtain the continuum vector potential 
\begin{equation}
{\vfun{\mbfit{A}}{\mbfit{X}} = \sum_{\mbf{x} \in \edgegrid}  \vfun{a}{x} \, \delta_h({\mbf{x} - \mbfit{X}} )h^3}. \label{interpcontVP}
\end{equation}
Lastly, we take the continuum curl of {$\mbfit{A}(\mbfit{X})$} with respect to {$\mbfit{X}$},
\begin{align}
 (\curl{\mbfit{A}})(\mbfit{X}) &= \sum_{\mbf{x} \in \edgegrid} \vfun{a}{x} \times (\grad \delta_h) (\mbf{x} - \mbfit{X}) h^3,
\end{align}
and our new interpolation is completed by adding the mean flow $\mbf{u}_0$, that is,
\begin{equation}
 \vfun{\mbfit{U}}{\mbfit{X}} = (\mbfit{\mathcal{J}}\mbf{u})(\mbfit{X}) = \mbf{u}_0  + \sum_{\mbf{x} \in \edgegrid} \vfun{a}{x} \times (\grad \delta_h) (\mbf{x} - \mbfit{X}) h^3. \label{divfree_interp}
\end{equation} 
We note that the interpolation \Cref{interpcontVP} is not performed in the actual implementation of the scheme. Instead, $\grad \delta_h$ is computed on the edge-centered staggered grid $\edgegrid$ in \Cref{divfree_interp}.
Notice that, by construction, the interpolated velocity in \Cref{divfree_interp} is {continuously} divergence-free. 

There are two important features of our new interpolation scheme that are worth mentioning. First, in comparison to locally interpolating the velocity from the nearby fluid grid in the conventional IB method, our new interpolation scheme is non-local, in that it involves solving the discrete Poisson problem \Cref{vpot_poisson}. Second, {if the regularized delta function $\delta_h$ is $\mathscr{C}^k$,  we note that the interpolated velocity field given by \Cref{divfree_interp} is a globally-defined function that is $\mathscr{C}^{k-1}$. We can think of the regularized delta function concentrated at $\mbfit{X}$ as being defined everywhere with zero  outside a cube of fixed edge length (e.g. $6h$ for the $\mathscr{C}^3$ 6-point kernel \cite{Bao2016_new6pt}).  Moreover, the continuity of derivatives of $\delta_h$ also applies globally, including at the edges for the cube. Since the continuum vector potential defined by \Cref{interpcontVP}  is a finite sum of such $\mathscr{C}^k$ functions, and the interpolated velocity field $\mbfit{U}(\mbfit{X})$ is obtained by differentiating $\mbfit{A}(\mbfit{X})$ once, then the resulting interpolated velocity field must have $k-1$ continuous derivatives.} 
{Note that  if we use an IB kernel that is $\mathscr{C}^1$, then the interpolated velocity $\mbfit{U}$ is $\mathscr{C}^0$, and $\div \mbfit{U}$ is defined in only a piecewise manner. This naturally brings into question whether the volume of a closed surface is strictly conserved as the surface passes over the discontinuity of the velocity derivatives. Indeed, we observe numerically that the DFIB method offers only marginal improvement in volume conservation for $\mathscr{C}^1$ kernel functions, such as the standard 4-point kernel \cite{Peskin2003_IBreview}, unless the Lagrangian mesh is discretized with impractically high resolution (8 markers per fluid meshwidth, see \autoref{fig:circleFine}). By contrast, we will show that with only a moderate Lagrangian mesh size (1 to 2 markers per fluid meshwidth), the DFIB method offers a substantial improvement in volume conservation for kernels of higher smoothness, which gives a continuously differentiable interpolated velocity $\mbfit{U}$. Further, we observe that volume conservation of the DFIB method improves with the smoothness of the interpolated velocity field.}

{In addition to the standard 4-point kernel (denoted by $\stndfour$),} the IB kernels considered in this paper include the $\mathscr{C}^3$  5-point and 6-point kernels \cite{C3IBkernels,Bao2016_new6pt} (denoted by $\newfive$ and $\newsix$ respectively), and the $\mathscr{C}^2$ 4-point B-spline kernel \cite{Unser1999}, 
\begin{equation}
\displaystyle
\phi^{B}_{4h}(r) =
\begin{cases}
\frac{2}{3} - r^2 + \frac{1}{2}r^3 & \quad  0 \leq |r| < 1, \\
\frac{4}{3} - 2r + r^2 - \frac{1}{6}r^3 & \quad  1 \leq |r| < 2, \\
0 &  \quad |r| \geq 2,
\end{cases} \label{bspline4pt}
\end{equation}
and the $\mathscr{C}^4$ 6-point B-spline kernel, 
\begin{equation}
\displaystyle
\phi^{B}_{6h}(r) =
\begin{cases}
\frac{11}{20}-\frac{1}{2}r^2+\frac{1}{4}r^4 - \frac{1}{12}r^5  &  \quad 0 \leq |r| < 1, \\
\frac{17}{40} + \frac{5}{8}r - \frac{7}{4}r^2 + \frac{5}{4}r^3 - \frac{3}{8}r^4 + \frac{1}{24}r^5 & \quad 1 \leq |r| < 2, \\
\frac{81}{40} - \frac{27}{8}r + \frac{9}{4}r^2 - \frac{3}{4}r^3 + \frac{1}{8}r^4 - \frac{1}{120}r^{{5}} &  \quad 2 \leq |r| < 3, \\
0 & \quad |r| \geq 3.
\end{cases} \label{bspline6pt}
\end{equation}
These B-spline kernels are members of a sequence of functions obtained by recursively convolving each successive kernel function against a rectangular pulse (also known as the window function), starting from the window function itself \cite{Unser1999}. The limiting function in this sequence is a Gaussian \cite{Unser1992}, which is exactly translation-invarant and isotropic. The family of IB kernels with nonzero even moment conditions, such as $\stndfour$ and $\newsix$, also have a Gaussian-like shape, but it is not currently known whether this sequence of functions also converges to a Gaussian. 


\subsection{The force-spreading operator}
The force-spreading operator {$\mbfit{S}$} is constructed to be adjoint to the velocity-interpolation opeartor {$\mbfit{S}^{*}$} so that energy is conserved by the Lagrangian-Eulerian interaction,
\begin{equation}
(\mbf{u}, \mbfit{S} \mbf{F})_{\mbf{x}} = (\mbfit{S}^*\mbf{u}, \mbf{F})_{\mbf{X}},
\end{equation}
where $(\cdot,\cdot)_{\mbf{x}}$ and $(\cdot,\cdot)_{\mbf{X}}$ denote the corresponding discrete inner products on the Eulerian and Lagrangian grids.
In other words, the power generated  by the elastic body forces is transferred to the fluid without loss,\footnote{Here and in similar expressions, $ \sum_{\mbf{x} \in \facegrid} \vfun{u}{x} \cdot \vfun{f}{x} h^3$ is a shorthand for $\sum_{i=1}^3 \sum_{\mbf{x} \in \facegrid} u_i(\mbf{x}) f_i(\mbf{x}) h^3$.}
\begin{align}
 \sum_{\mbf{x}\in \facegrid} \vfun{u}{x} \cdot \vfun{f}{x} \, h^3 &= \sum_{m=1}^M \mbf{U}_m \cdot \mbf{F}_m \Delta \mbf{s},
 \label{power_identity}
\end{align}
where $\mbf{U}_m$ is the Lagrangian marker velocity at $\mbf{X}_m$, and $\mbf{F}_m \Delta \mbf{s}$ is the Lagrangian force applied to the fluid by the Lagrangian marker $\mbf{X}_m$.
Our goal is to find an Eulerian force density $\vfun{f}{x}$ that satisfies the power identity \Cref{power_identity}. To see what \Cref{power_identity} implies about $\vfun{f}{x}$, we rewrite both sides in terms of $\vfun{a}{x}$. On the left-hand side of \Cref{power_identity}, we use \Cref{vpot_def} to obtain
\begin{align}
 \sum_{\mbf{x} \in \facegrid} \vfun{u}{x} \cdot \vfun{f}{x} \, h^3  &= \mbf{u}_0 \cdot \sum_{\mbf{x} \in \facegrid} \vfun{f}{x} h^3 + \sum_{\mbf{x} \in \facegrid } ( \dcurl{h} \mbf{a} ) (\mbf{x}) \cdot \vfun{f}{x} \, h^3  \notag \\
 &= \mbf{u}_0 \cdot \mbf{f}_0 V + \sum_{\mbf{x} \in \edgegrid} \vfun{a}{x} \cdot ( \dcurl{h} \mbf{f}) (\mbf{x}) \, h^3, \label{power_lhs}
\end{align}
where the average of $\vfun{f}{x}$ over the domain is
\begin{equation}
\mbf{f}_0 = \frac{1}{V} \sum_{\mbf{x} \in \facegrid} \vfun{f}{x} \, h^3. \label{sum_eulerianf}
\end{equation}
Note that we have used the summation-by-parts identity \Cref{sumbyparts2} to transfer the discrete curl operator $\dcurl{h}$ from  $\mbf{a}(\mbf{x})$ to $\mbf{f(\mbf{x})}$, and thus, the grid on which the summation is performed in \Cref{power_lhs} is $\edgegrid$ not $\facegrid$.
On the the right-hand side of \Cref{power_identity}, we substitute for $\mbf{U}_m$ by using the divergence-free velocity interpolation \Cref{divfree_interp},  
\begin{align}
\sum_{m=1}^M \mbf{U}_m \cdot \mbf{F}_m \Delta \mbf{s} &= \mbf{u}_0 \cdot \sum_{m=1}^M \mbf{F}_m \Delta \mbf{s} + \sum_{m=1}^M \sum_{\mbf{x} \in \edgegrid}  
 \vfun{a}{x} \times (\grad \delta_h) (\mbf{x} - \mbf{X}_m) \cdot (\mbf{F}_m \Delta \mbf{s})\, h^3 \notag \\
 &= \mbf{u}_0 \cdot \sum_{m=1}^M \mbf{F}_m \Delta \mbf{s} + \sum_{\mbf{x} \in \edgegrid} \vfun{a}{x} \cdot \sum_{m=1}^M (\grad \delta_h) (\mbf{x} - \mbf{X}_m) \times (\mbf{F}_m \Delta \mbf{s}) \, h^3. \label{power_rhs}
\end{align}
Since $\mbf{u}_0$ and $\vfun{a}{x}$ are arbitrary (except for $\ddiv{h} \mbf{a} =0$), the power identity \Cref{power_identity} is satisfied if and only if
\begin{equation}
 \mbf{f}_0 = \frac{1}{V} \sum_{m=1}^M \mbf{F}_m \Delta \mbf{s} \label{sum_lagrangianf}
\end{equation}
and 
\begin{equation}
 \dcurl{h} \mbf{f}  = \sum_{k=1}^M (\grad \delta_h) (\mbf{x} - \mbf{X}_m) \times (\mbf{F}_m \Delta \mbf{s}) + \dgrad{h} \varphi, \, \text{ for all } \mbf{x} \in \edgegrid,
 \label{fspreading_eq}
\end{equation}
where $\varphi$ is an arbitrary scalar field that lives on the node-centered grid $\nodegrid$.
Note that we have the freedom to add the term $\dgrad{h} \varphi$ in \Cref{fspreading_eq}, since from the identity \Cref{sumbyparts1} and $\ddiv{h} \mbf{a} = 0$, we have
\begin{equation*}
 \sum_{\mbf{x} \in \edgegrid} \vfun{a}{x} \cdot  \left( \dgrad{h} \varphi \right) h^3 = -\sum_{\mbf{x} \in \nodegrid} \left( \ddiv{h} \mbf{a} \right) (\mbf{x}) \, \varphi(\mbf{x}) \, h^3  = 0.
\end{equation*}
Indeed, we are required to include this term since the left-hand side of \Cref{fspreading_eq} is discretely divergence-free but there is no reason to expect the first term on the right-hand side of \Cref{fspreading_eq} is also divergence-free. Note that it is not required to find $\varphi$ in order to determine $\vfun{f}{x}$, because we can eliminate $\varphi$ by taking the discrete curl on both sides of \Cref{fspreading_eq},
\begin{equation}
 \dcurl{h} ( \dcurl{h} \mbf{f}) = \dcurl{h} \left( \sum_{m=1}^M (\grad \delta_h) (\mbf{x} - \mbf{X}_m) \times (\mbf{F}_m \Delta \mbf{s}) \right), \, \text{ for all } \mbf{x} \in \edgegrid.
\label{fspreading_eq2}
\end{equation}
By imposing the gauge condition
\begin{equation}
 \ddiv{h} \mbf{f} = 0, \label{force_gauge}
\end{equation}
we obtain a vector Poisson equation for $\vfun{f}{x}$,
\begin{equation}
 -(\vlapl{h} \, \mbf{f})(\mbf{x})= \dcurl{h} \left( \sum_{m=1}^M (\grad \delta_h) (\mbf{x} - \mbf{X}_m) \times (\mbf{F}_m \Delta \mbf{s}) \right), \, \text{ for all } \mbf{x} \in \edgegrid.
\label{fvectorpoisson}
\end{equation}
 Note again that $\grad \delta_h$ is computed on $\edgegrid$, so that the cross-product with $\mbf{F}_m$ is face-centered, which agrees with the left-hand side of \Cref{fvectorpoisson}. Note that the solution of \Cref{fvectorpoisson} can be uniquely determined by the choice of $\mbf{f}_0$.
Like our velocity interpolation scheme, the new force-spreading scheme is also non-local because it requires solving discrete Poisson equations. We remark that the new force-spreading scheme is also constructed so that the resulting force density $\vfun{f}{x}$ is discretely divergence-free. This means that $\vfun{f}{x}$ includes the pressure gradient that is generated by the Lagrangian forces. We do not see a straightforward way to separate the pressure gradient from $\vfun{f}{x}$ in case it is needed for output purposes. 


\section{Time-stepping scheme} \label{sec_timestepping}
In this section, we present a second-order time-stepping scheme, similar to the ones developed previously \cite{IBAMR_HeartValve}, that evolves the spatially-discretized system \Cref{momentum_spatial,discrete_divzero,discrete_interp,}. 
Let $\mbf{u}^n, \mbf{X}^n$ denote the approximations of the fluid velocity and Lagrangian marker velocities at time $t_n = n\Delta t$. To advance the solutions to $\mbf{u}^{n+1}$ and $\mbf{X}^{n+1}$, we perform the following steps:

\begin{enumerate}[Step 1.]
 \item First, update the Lagrangian markers to the intermediate time step $n+\frac{1}{2}$ using the interpolated velocity, 
\begin{equation}
\widetilde{\mbf{X}}^{n+\frac{1}{2}} = \mbf{X}^n + \frac{\Delta t}{2} \mbfit{S}^{\star}\left[\mbf{X}^n\right]\mbf{u}^n. \label{timestep_1}
\end{equation}  \label{firststep}

 \item \label{secondstep} Evaluate the intermediate Lagrangian force density at $\widetilde{\mbf{X}}^{n+\frac{1}{2}}$ from the force density functional or the energy functional, and spread it to the Eulerian grid using the force-spreading scheme to get
\begin{equation}
 \mbf{f}^{n+\frac{1}{2}} = \mbfit{S}\left[\widetilde{\mbf{X}}^{n+\frac{1}{2}}\right] \mbf{F}^{n+\frac{1}{2}}. \label{timestep_2}
\end{equation}

 \item \label{thirdstep} Solve the fluid equations  on the periodic grid \cite{Devendran2012_IBviscoelastic},
\begin{equation}
\left\{
\begin{array}{l}
\displaystyle \rho \left( \frac{\mbf{u}^{n+1}-\mbf{u}^{n}}{\Delta t} + \widetilde{\mbf{N}}^{n+\frac{1}{2}} \right) + \dgrad{h}{p}^{n+\frac{1}{2}} = \mu \vlapl{h} \left(\frac{\mbf{u}^{n+1} + \mbf{u}^n}{2}\right) + \mbf{f}^{n+\frac{1}{2}}, \\
\\
\ddiv{h} \mbf{u}^{n+1} = 0,
\end{array}
\right. \label{momentum_spatial_time}
\end{equation}
where the second-order Adams-Bashforth (AB2) method is applied to approximate the nonlinear advection term 
\begin{equation}
 \widetilde{\mbf{N}}^{n+\frac{1}{2}} = \frac{3}{2}\mbf{N}^{n} - \frac{1}{2}\mbf{N}^{n-1},
\end{equation}
and $\mbf{N}^n = {\mbfit{N}(\mbf{u}^n) }$.

\item \label{fourthstep} In the last step, update the Lagrangian markers $\mbf{X}^{n+1}$ by using the mid-point approximation
\begin{equation}
 \mbf{X}^{n+1} = \mbf{X}^n + \Delta t \ \mbfit{S}^{\star}\left[\widetilde{\mbf{X}}^{n+\frac{1}{2}}\right] \left( \frac{\mbf{u}^{n+1} + \mbf{u}^n}{2} \right). \label{timestep_4}
\end{equation}
\end{enumerate}
Note that the time-stepping scheme described above requires two starting values because of the treatment of the nonlinear advection term using the AB2. To get the starting value at $t=\Delta t$, we can use the second-order Runge-Kutta (RK2) scheme described in \cite{Peskin2003_IBreview,Devendran2012_IBviscoelastic}. 

In \autoref{table:DFIBcost} we compare the cost of DFIB and IBMAC for the above IB scheme in terms of the number of the two cost-dominating procedures: the scalar Poisson solver which costs $\mathcal{O}(N^d \log N)$ using FFT on the periodic domain, where $d\in\{2,3\} $ is the spatial dimension, and spreading/interpolation of a scalar field between the Eulerian grid and the Lagrangian mesh which costs $\mathcal{O}(M)$. In summary, DFIB is only more expensive than IBMAC by 4 scalar Poisson solves for two-dimensional (2D) problems, and is more expensive by 9 scalar Poisson solves and 9 scalar interpolation and spreading for three dimensional (3D) problems. Therefore, the DFIB method is about two times slower than IBMAC per time step in 3D. We point out that if the RK2 scheme \cite{Peskin2003_IBreview,Devendran2012_IBviscoelastic} is employed rather than the scheme above, then we can save one interpolation step per time step,  but the fluid equations need to be solved twice.


\renewcommand{\arraystretch}{1.3}
\begin{table}
\centering
\begin{tabular}{| c | c | c | c | c | c | c | c | c |}
\cline{2-9}
\multicolumn{1}{c|}{} & \multicolumn{4}{c|}{ \# of scalar Poisson solves } & \multicolumn{4}{c|}{ \# of scalar interpolation/spreading} \\
\cline{2-9}
\multicolumn{1}{c|}{} & \multicolumn{2}{c|}{ 2D } & \multicolumn{2}{c|}{ 3D } & \multicolumn{2}{c|}{ 2D } & \multicolumn{2}{c|}{ 3D } \\
\cline{2-9}
\multicolumn{1}{c|}{} & DFIB & IBMAC & DFIB & IBMAC & DFIB & IBMAC & DFIB & IBMAC \\
\hline
$\mbfit{S}^{\star}$ in \Cref{timestep_1} & 1 & - & 3 & - & 2 & 2 & 6 & 3 \\
$\mbfit{S}$ in \Cref{timestep_2} & 2 & - & 3 & - & 2 & 2 & 6 & 3 \\
Fluid solver & 3 & 3 & 4 & 4 & - & - & - & - \\
$\mbfit{S}^{\star}$ in \Cref{timestep_4} & 1 & - & 3 & - & 2 & 2 & 6 & 3 \\
\hline
Total & 7 & 3 & 13 & 4 & 6 & 6 & 18 & 9 \\
\hline
\end{tabular}
\caption{Cost of DFIB versus IBMAC in terms of the number of scalar Poisson solves and interpolation/spreading of a scalar from/to the Eulerian grid. }
\label{table:DFIBcost}
\end{table}


\section{Numerical Results} \label{sec_numericalresults}

This section presents numerical results of the DFIB method for various benchmark problems in 2D and 3D. We first consider in 2D a thin elastic membrane subject to surface tension of the membrane only. The continuum solution of this simple 2D problem has the special feature that the tangential component of the elastic force vanishes, and therefore, the normal derivative of the tangential fluid velocity does not suffer any jump across the immersed boundary. This has the effect that second-order convergence in the fluid velocity $\mbfit{u}$ and the Lagrangian deformation map ${\mbfit{\mathcal{X}}}$ can be achieved \cite{Griffith2005_IBaccuracy}. In the second set of tests, we compare volume conservation in 2D, \ie, area conservation of DFIB and IBMAC by applying them to a circular membrane under tension, and we discuss the connection between area conservation and the choice of Lagrangian marker spacing relative to the Eulerian grid size. In the third set of computations, we apply the DFIB method to a problem in which a 2D elastic membrane actively evolves in a parametrically-forced system. In the last set of numerical experiments, we extend the surface tension problem to 3D, and compare volume conservation of DFIB with that of IBMAC.

\subsection{A thin elastic membrane with surface tension in 2D} \label{sec_convergence}

It is well-known that the solutions to problems involving an infinitely thin massless membrane interacting with a viscous incompressible fluid possess jump discontinuities across the interface in the pressure and in the normal derivative of the velocity due to singular forcing at the interface \cite{Lee2003_IIM,Lai2001}. These sharp jump discontinuities cannot be fully resolved by the conventional IB method because of the use of the regularized delta function at the interface. Consequently, the numerical convergence rate for the Lagrangian deformation map ${\mbfit{\mathcal{X}}}$ is generally only first order even if the discretization is carried out with second-order accuracy. To achieve the expected rate of convergence, we consider problems with solutions that possess sufficient smoothness. 

As a simple benchmark problem with a sufficiently smooth continuum solution we consider a thin elastic membrane that deforms in response to surface tension only.
Suppose that the elastic interface $\Gamma$ is discretized by a collection of Lagrangian markers $\mbf{X} = \left\{\mbf{X}_1, \dots, \mbf{X}_{M}\right\}$. The discrete elastic energy functional associated with the surface tension of the membrane is the total (polygonal) arc-length of the interface \cite{Kim2010_2Dfoam}, 
\begin{equation}
E[\mbf{X}_1,\dots,\mbf{X}_{M}] = \gamma \sum_{m=1}^{M} \left| \mbf{X}_m - \mbf{X}_{m-1} \right|,
\end{equation} 
where $\mbf{X}_0 = \mbf{X}_{M}$ and $\gamma$ is the surface tension constant (energy per unit length). The Lagrangian force generated by the energy functional at the marker $\mbf{X}_m$ is 
\begin{equation}
\mbf{F}_m \Delta {s} = -\frac{\partial E}{\partial \mbf{X}_m} = \gamma \left( \frac{\mbf{X}_{m+1} - \mbf{X}_m}{|\mbf{X}_{m+1} - \mbf{X}_m|} -  \frac{\mbf{X}_{m} - \mbf{X}_{m-1}}{|\mbf{X}_{m} - \mbf{X}_{m-1}|}\right).
\end{equation}
In our tests, we set the initial configuration of the membrane to be the ellipse
\begin{equation}
{\mbfit{\mathcal{X}}}(s,0) = L \cdot \left( \frac{1}{2}+\frac{5}{28} \cos(s), \ \frac{1}{2}+\frac{7}{20}\sin(s)  \right),  \quad s \in [0, 2\pi].
\end{equation}
The Eulerian fluid domain $\Omega = [0,L]^2$ is discretized by a uniform $N \times N$ Cartesian grid with meshwidth $h=\frac{L}{N}$ in each direction. The elastic interface $\Gamma$ is discretized by a uniform  Lagrangian mesh of size  $M = \lceil\pi N\rceil$ in the Lagrangian variable $s$, so that the Lagrangian markers $\mbf{X}= \left\{\mbf{X}_1, \dots, \mbf{X}_{M}\right\}$ are physically separated by a distance of approximately $\frac{h}{2}$ in the equilibrium circular configuration. 
In all of our tests, we set $L=5, ~ \rho = 1, ~\gamma = 1, ~ \mu = 0.1$.  The time step size is chosen to be $\Delta t = \frac{h}{2}$ to ensure the stability of all simulations up to $t = 20$ when the elastic interface is empirically observed to be in equilibrium. 

\begin{figure}
\centering
\includegraphics[width = \linewidth]{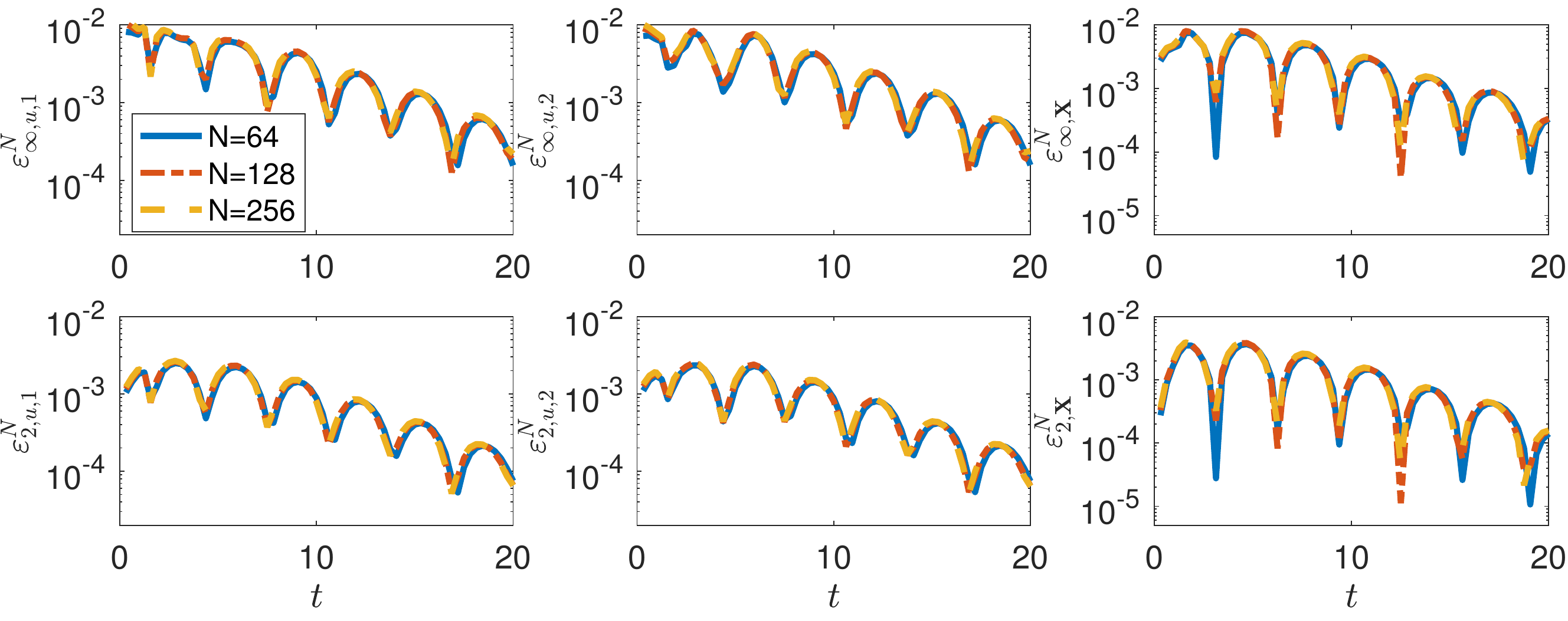}
\caption{Rescaled $l_{\infty}$-norm (top panel) and $l_2$-norm (bottom panel) errors of the $x,y$-component of the fluid velocity defined by \Cref{errorvel}, and errors of the Lagrangian deformation map defined by \Cref{errorX} for the 2D surface tension problem are plotted as a function of time from $t=0$ to $t=20$. 
The left and middle columns show errors of the fluid velocity components and the right column shows errors of the Lagrangian deformation map.
The Eulerian grid sizes are $N=64,128, 256$ and the corresponding Lagrangian mesh sizes are $M= 202, 403, 805$, so that the spacing between two Lagrangian markers is kept at a distance of approximately $\frac{h}{2}$ in the equilibrium configuration. For the finer grid resolution $N = 128, 256$,  the errors in each norm are multiplied by a factor of 4 and $4^2$ respectively. After rescaling, the error curves of the finer grid resolution almost align with the error curves of grid resolution $N=64$, which indeed confirms that second-order convergence in $\mbf{u}$ and $\mbf{X}$ is achieved. 
For this set of computations, we use the $\mathscr{C}^3$ 6-point IB kernel in the discrete delta function, and the time step size is chosen to be $\Delta t = \frac{h}{2}$.  }
\label{fig:surf_vXconv}
\end{figure}

We denote by $\mbf{u}^N(t)$ the computed fluid velocity field and by $\mathcal{I}^{2N\rightarrow N}$  a restriction operator from the finer grid of size $2N\times 2N$ to the coarser grid of size $N\times N$. The discrete $l_p$-norm of the successive  error in the velocity component $u_i$ is defined by
\begin{equation}
 \varepsilon^N_{p,u,i}(t) = \left\| u_i^N(t) - \mathcal{I}^{2N\rightarrow N} u_i^{2N}(t)\right\|_p. \label{errorvel}
\end{equation}
To avoid artifacts in the error-norm computation because of Lagrangian markers getting too clustered during the simulation, we reparametrize the interface (for the purpose of the error computation only) from the computed markers using periodic cubic splines after each time step, 
and compute the $l_p$-norm error of $\mbf{X}$ based on a collection of  $M'$ uniformly sampled markers $\widetilde{\mbf{X}}$  from the reparametrized interface, that is,
\begin{equation}
 \varepsilon^N_{p,\mbf{X}}(t) = \left\| \widetilde{\mbf{X}}^N(t) - \widetilde{\mbf{X}}^{2N}(t)\right\|_p, \label{errorX}
\end{equation}
where $M'$ does not change with $N$. We emphasize that the resampled markers are only used to compute the error norm and are discarded after each time step.
In \autoref{fig:surf_vXconv} the successive $l_{\infty}$-norm and $l_2$-norm errors of the $x,y$-component of the fluid velocity and of the deformation map are plotted  as a function of time from $t=0$ to $t=20$ for grid resolution $N=64,128$ and $256$. The number of resampled markers for computing $\varepsilon^N_{p,\mbf{X}}(t)$ is $M'=128$. To clearly visualize that second-order convergence is achieved by our scheme, we multiply the computed errors for the finer grid resolution $N=128$ and $256$  by a factor of 4 and $4^2$ respectively, and plot them along with errors for the coarser grid resolution $N=64$ in \autoref{fig:surf_vXconv}. The observation that all three error curves almost align with each other (as shown in \autoref{fig:surf_vXconv})  confirms that  second-order convergence in $\mbf{u}$ and $\mbf{X}$ is achieved. 



\subsection{Area conservation and IB marker spacing} \label{sec_spacing}

As an immediate consequence of fluid incompressibility, the volume enclosed by a closed immersed boundary should be exactly conserved as it deforms and moves with the fluid. However, it is observed that even in the simplest 
scenario of a pressurized membrane in its circular equilibrium configuration \cite{Griffith2012_IBMACvolume}, the volume error of an IB method with conventional interpolation and spreading  systematically grows at a rate proportional to the pressure jump across the elastic interface \cite{Peskin1993_IBmodified}. In this set of tests, we demonstrate that, the ``volume'' or area enclosed by a 2D membrane is well-conserved  by the DFIB method when the Lagrangian interface is sufficiently resolved.

We follow the same problem setup as in the test described in \cite{Griffith2012_IBMACvolume}. A thin elastic membrane ${\mbfit{\mathcal{X}}}(s,t)$, initially in a circular equilibrium configuration, 
\begin{equation}
{\mbfit{\mathcal{X}}}(s,0) = \left( \frac{1}{2}+\frac{1}{4} \cos(s), \frac{1}{2}+\frac{1}{4}\sin(s)  \right),  \quad s \in [0, 2\pi], 
\label{membrane_circle}
\end{equation}
is immersed in a periodic unit cell $\Omega = [0,1]^3$ with zero initial background flow. The Lagrangian force density on the interface is described by 
\begin{equation}
 \mbfit{F}(s,t)  =  \kappa \frac{\partial^2 {\mbfit{\mathcal{X}}}}{\partial s^2}, \label{force2nd}
\end{equation}
in which $\kappa$ is the uniform stiffness coefficient. The elastic membrane is discretized by a uniform Lagrangian mesh of $M$ points in the variable $s$. We approximate the Lagrangian force density by 
\begin{equation}
 \mbf{F}_{m}  = \frac{\kappa}{(\Delta s)^2} \, ( \mbf{X}_{m+1} - 2\mbf{X}_m + \mbf{X}_{m-1} ), \label{springforce_2nd}
\end{equation}
which corresponds to a collection of Lagrangian markers connected by linear springs of zero rest length with stiffness $\kappa$. For this problem, since the elastic interface is initialized in the equilibrium configuration with zero background flow, any spurious fluid velocity and area loss incurred in the simulation are regarded as numerical errors.  

In our simulations, we set $\rho = 1, ~\mu = 0.1, ~\kappa = 1$. The size of the Eulerian grid is fixed at $128 \times 128$ with meshwidth $h = \frac{1}{128}$. The size of the Lagrangian mesh $M$ is chosen so that  two adjacent Lagrangian markers are separated by a physical distance of $h_s$  in the equilibrium configuration, that is, $M \approx 2\pi R / h_s$, where $R$ is the radius of the circular membrane.  {In addition to the Lagrangian markers, we also include a dense collection of passive tracers with $N_{\text{tracer}} = 20M$ to address the limiting case of moving the entire interface. These tracers are initially in the same configuration as the circular membrane in \Cref{membrane_circle}, and they move passively with the interpolated velocity according to \Cref{timestep_1,timestep_4}.} The time step size is set to be $\Delta t = \frac{h}{4}$ for stability. In all computations, we use the $\mathscr{C}^3$ 6-point IB kernel $\newsix$ to form the regularized delta function $\delta_h$.  

\begin{figure}
\subfloat[][IBMAC]{\includegraphics[width = \linewidth]{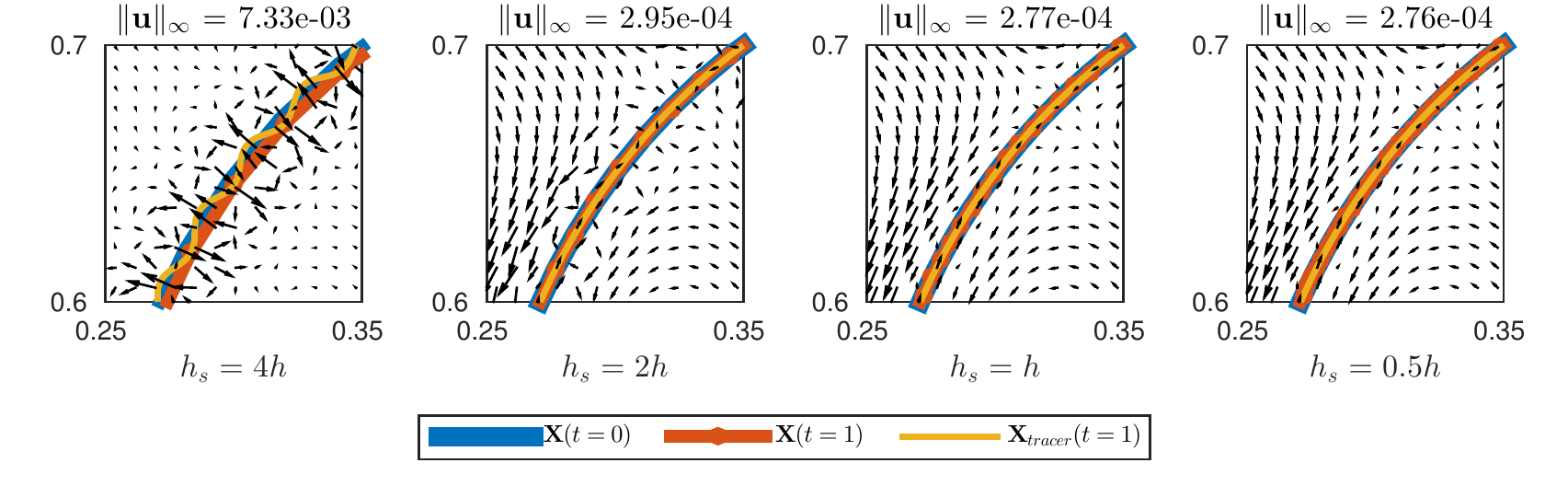}}\\
\subfloat[][DFIB]{\includegraphics[width = \linewidth]{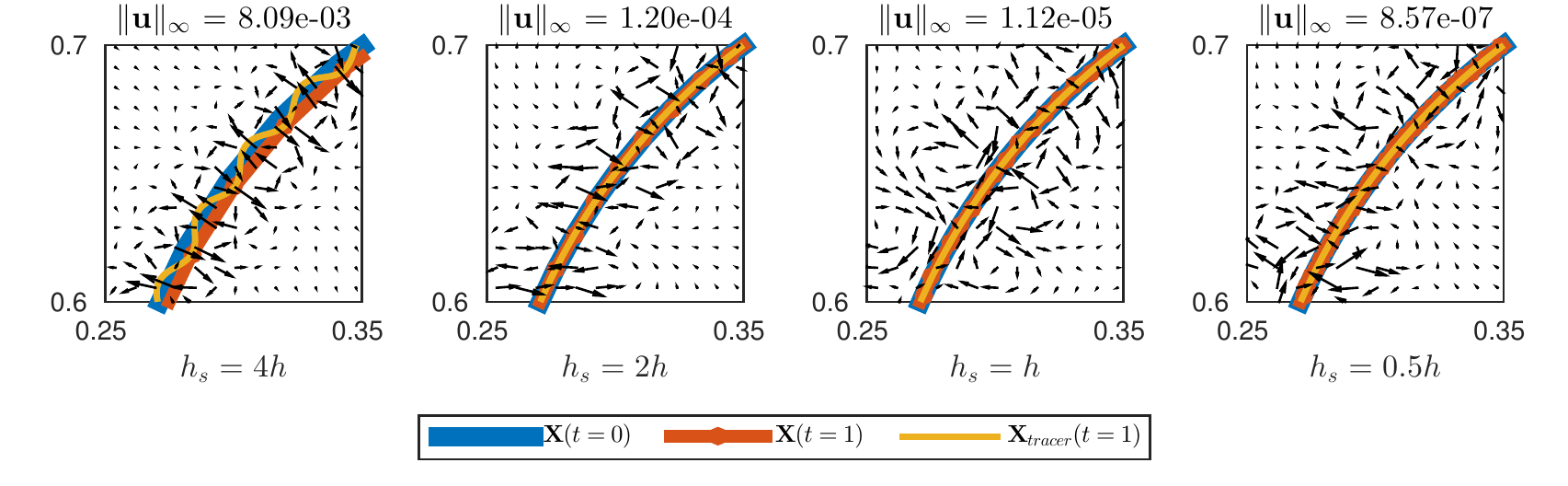}}
\caption{A magnified view of the quasi-static circular membrane and its nearby spurious velocity field for different Lagrangian mesh spacing $h_s = 4h, \, 2h, \, h$ and $\frac{h}{2}$, as indicated below each figure panel, while keeping $h=\frac{1}{128}$ fixed. The top panel (a) shows the computational results from IBMAC, and the bottom panel (b) shows the results from DFIB. The interface represented by the Lagrangian markers $\mbf{X}(t=1)$ is shown in red, the initial configuration $\mbf{X}(t=0)$ is shown in blue, and the interface represented by {$N_{\text{tracer}} = 20M$} passive tracers is shown in yellow, {where $M$ is the number of Lagrangian markers.} The time step size is set to be $\Delta t = \frac{h}{4}$ for stability. In the above computations, the $\mathscr{C}^3$ 6-point IB kernel $\newsix$ is used in IBMAC and  DFIB.}
\label{fig:circIBMAC_IBDF}
\end{figure}

In \autoref{fig:circIBMAC_IBDF} we compare the computational results of DFIB with those of IBMAC for different $h_s = 4h, \, 2h, \, h$ and $\frac{h}{2}$ (from left to right in \autoref{fig:circIBMAC_IBDF}). Each subplot of \autoref{fig:circIBMAC_IBDF} shows a magnified view of the same arc of the circular interface along with its nearby spurious fluid velocity field. The interface represented by the Lagrangian markers $\mbf{X}(t=1)$ is shown in red  and the initial configuration $\mbf{X}(t=0)$  is shown in the blue curve. 
{The interface represented by the passive tracers $\mbf{X}_{\text{tracer}}(t=1)$ is shown in the yellow curve. } In the first column of \autoref{fig:circIBMAC_IBDF} in which $h_s = 4h$, we see that the maximum spurious velocity $\|\mbf{u} \| _{\infty}$ of  IBMAC is of the same magnitude as that of DFIB. At such coarse resolution in the Lagrangian mesh, fluid apparently leaks through the gap between two adjacent markers, as can be observed by the wiggly pattern in the passive tracers. As the the Lagrangian mesh is refined gradually from $h_s = 4h$ to $\frac{h}{2}$ (from left to right in \autoref{fig:circIBMAC_IBDF}), we see that $\|\mbf{u}\|_{\infty}$ decreases from $10^{-3}$ to $10^{-7}$ in the DFIB method, whereas $\|\mbf{u}\|_{\infty}$ stops improving around $10^{-4}$ in IBMAC. Moreover, in the columns where $h_s = 2h,\, h, \, \frac{h}{2}$, we see a clear global pattern in the spurious velocity field in IBMAC, while the spurious velocity field of DFIB  appears to be much smaller in magnitude and random in pattern. 

\begin{figure}
\subfloat[][]{\includegraphics[width =.33\linewidth]{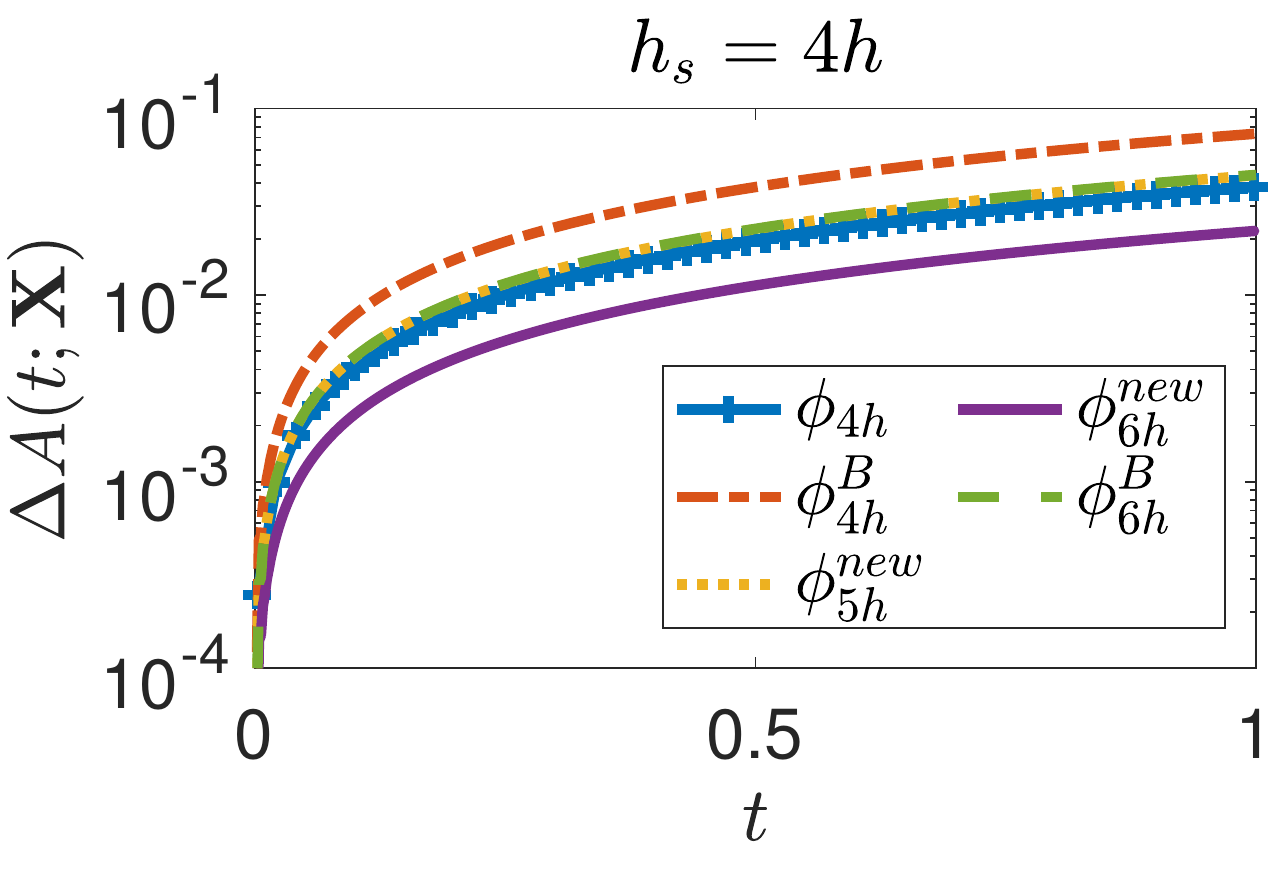} \label{fig:IBMACcircle_a} }
\subfloat[][]{\includegraphics[width =.33\linewidth]{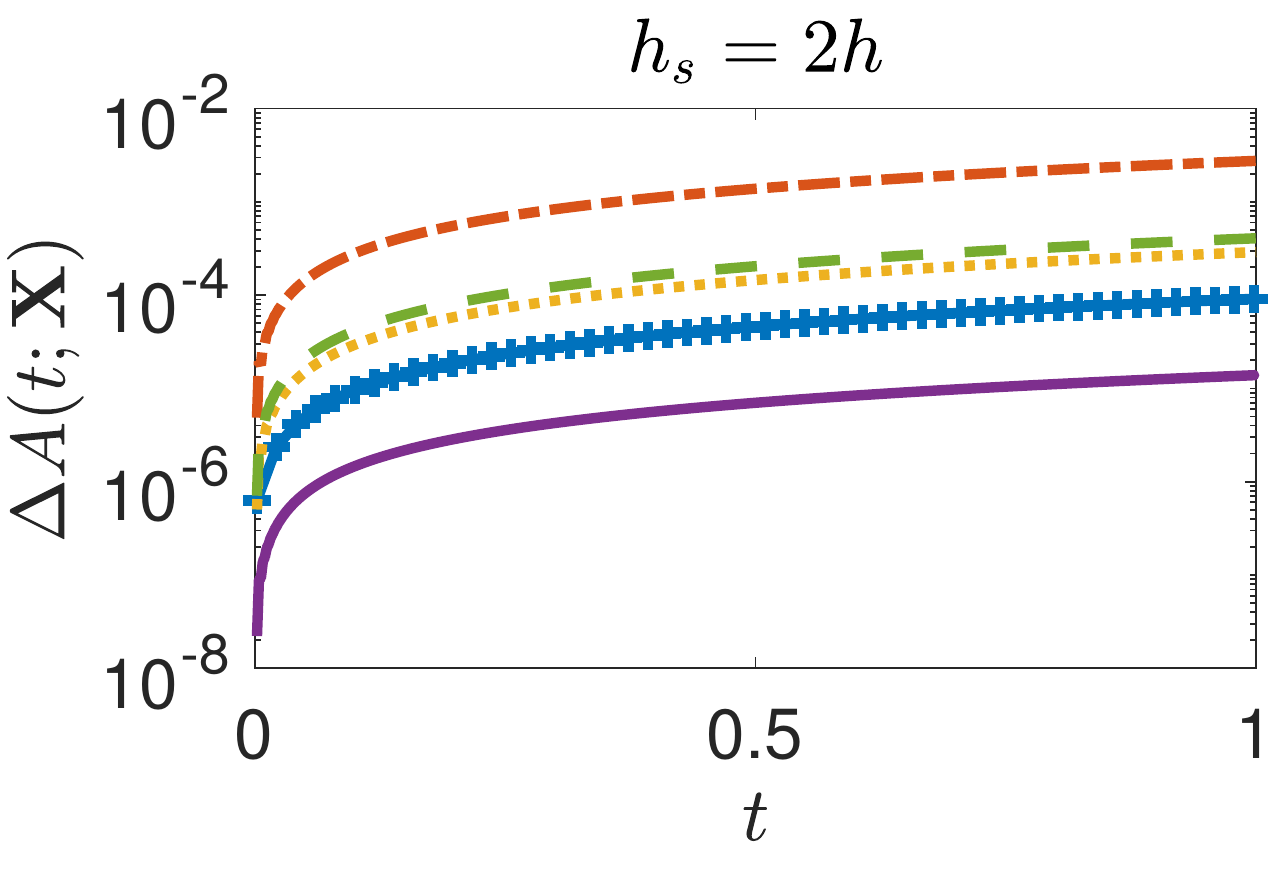} \label{fig:IBMACcircle_b} }
\subfloat[][]{\includegraphics[width =.33\linewidth]{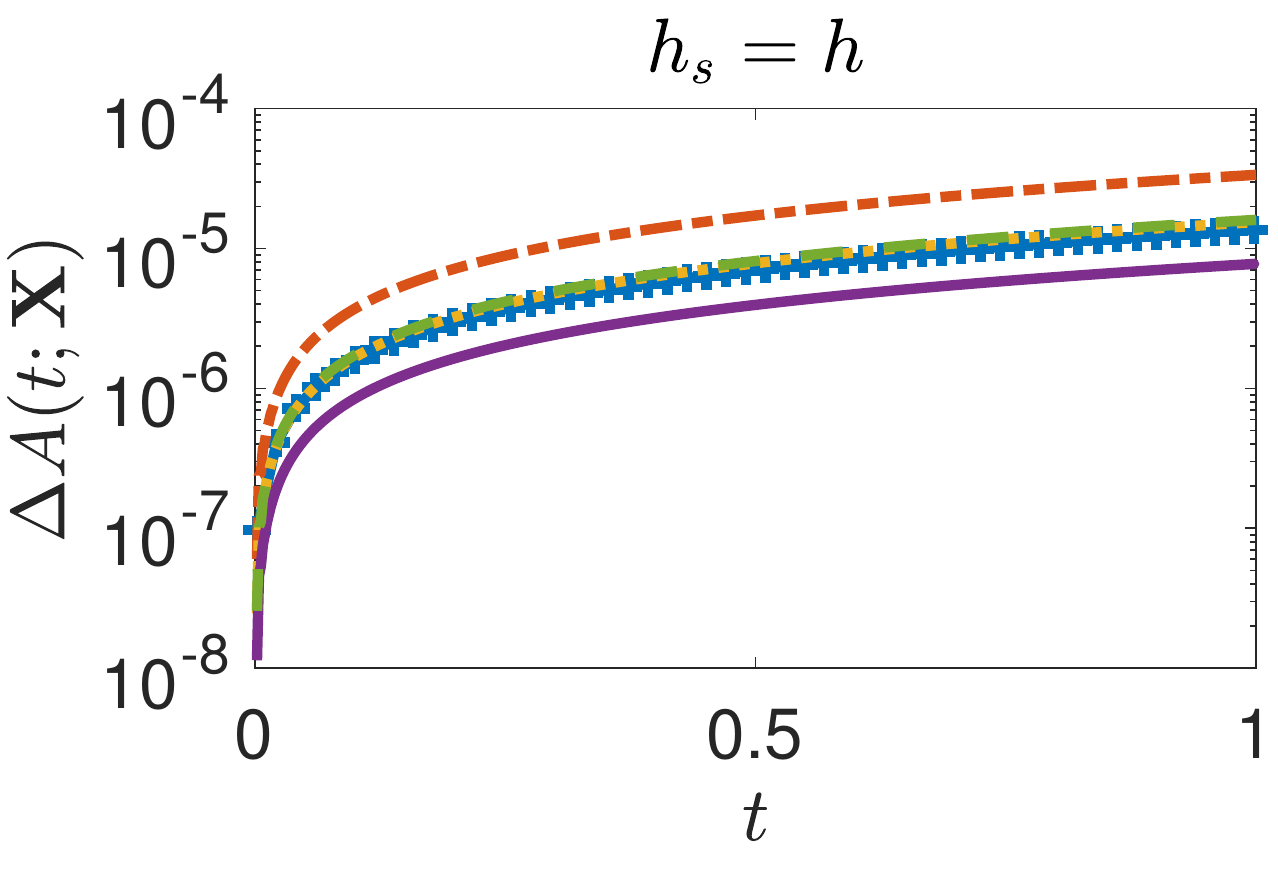} \label{fig:IBMACcircle_c}} \\
\subfloat[][]{\includegraphics[width =.33\linewidth]{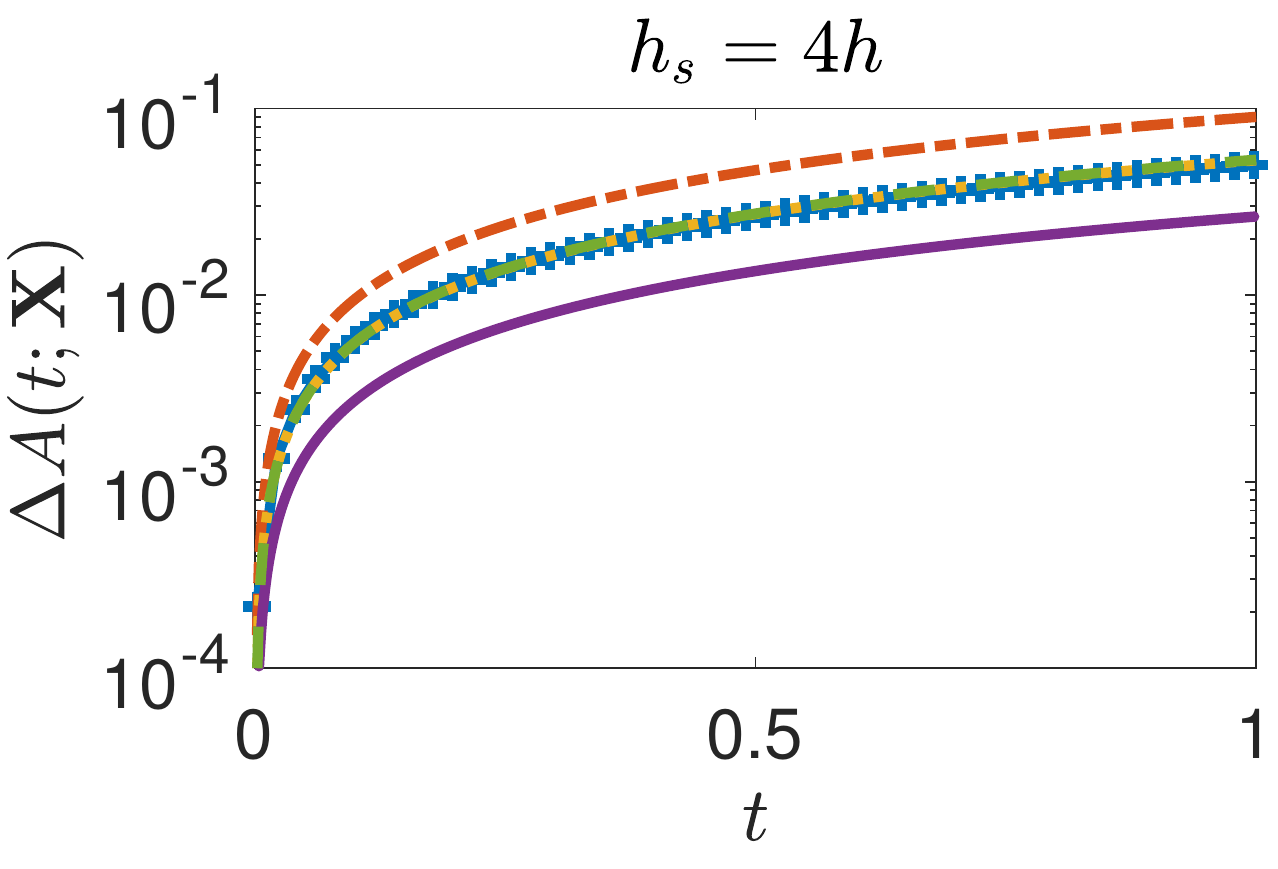} \label{fig:IBDFcircle_a}}
\subfloat[][]{\includegraphics[width =.33\linewidth]{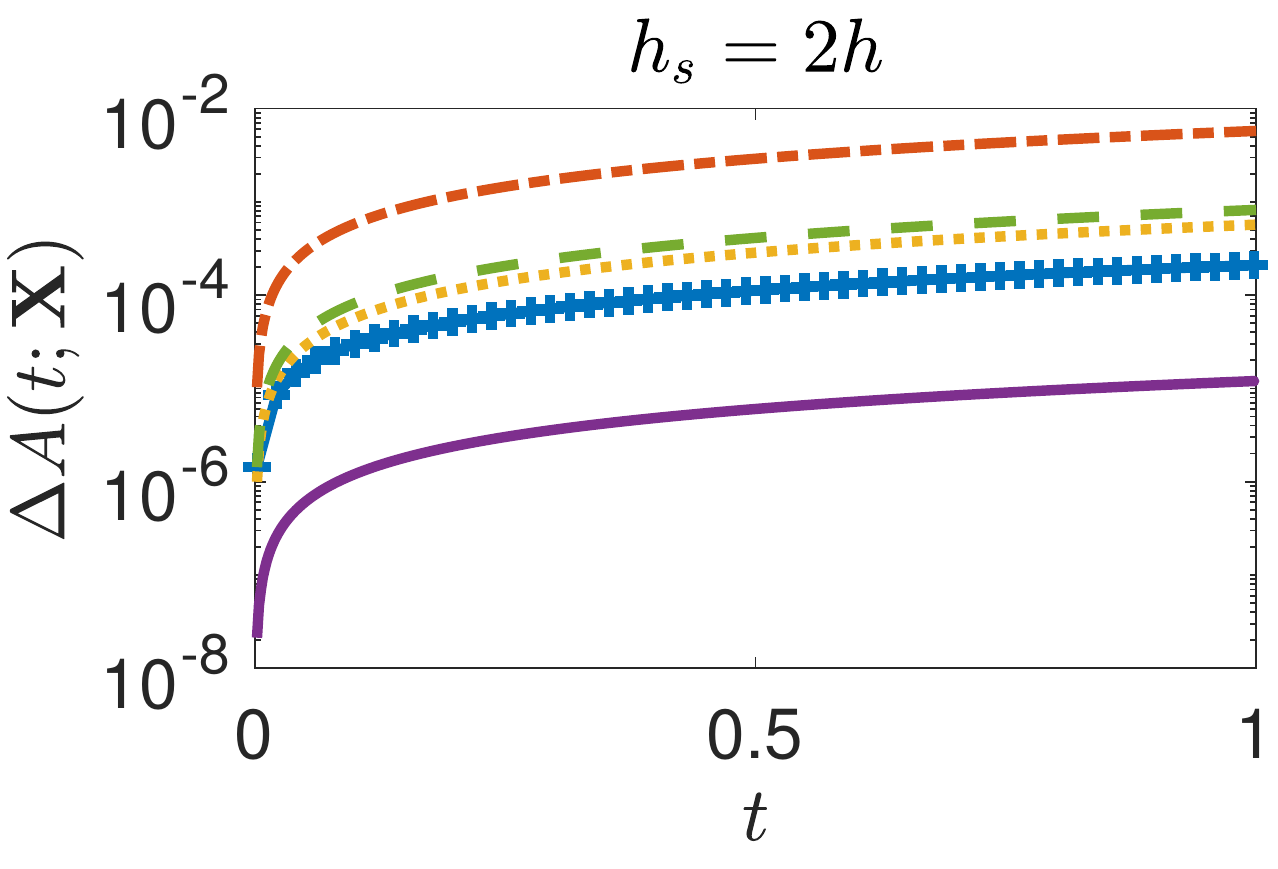} \label{fig:IBDFcircle_b}}
\subfloat[][]{\includegraphics[width =.33\linewidth]{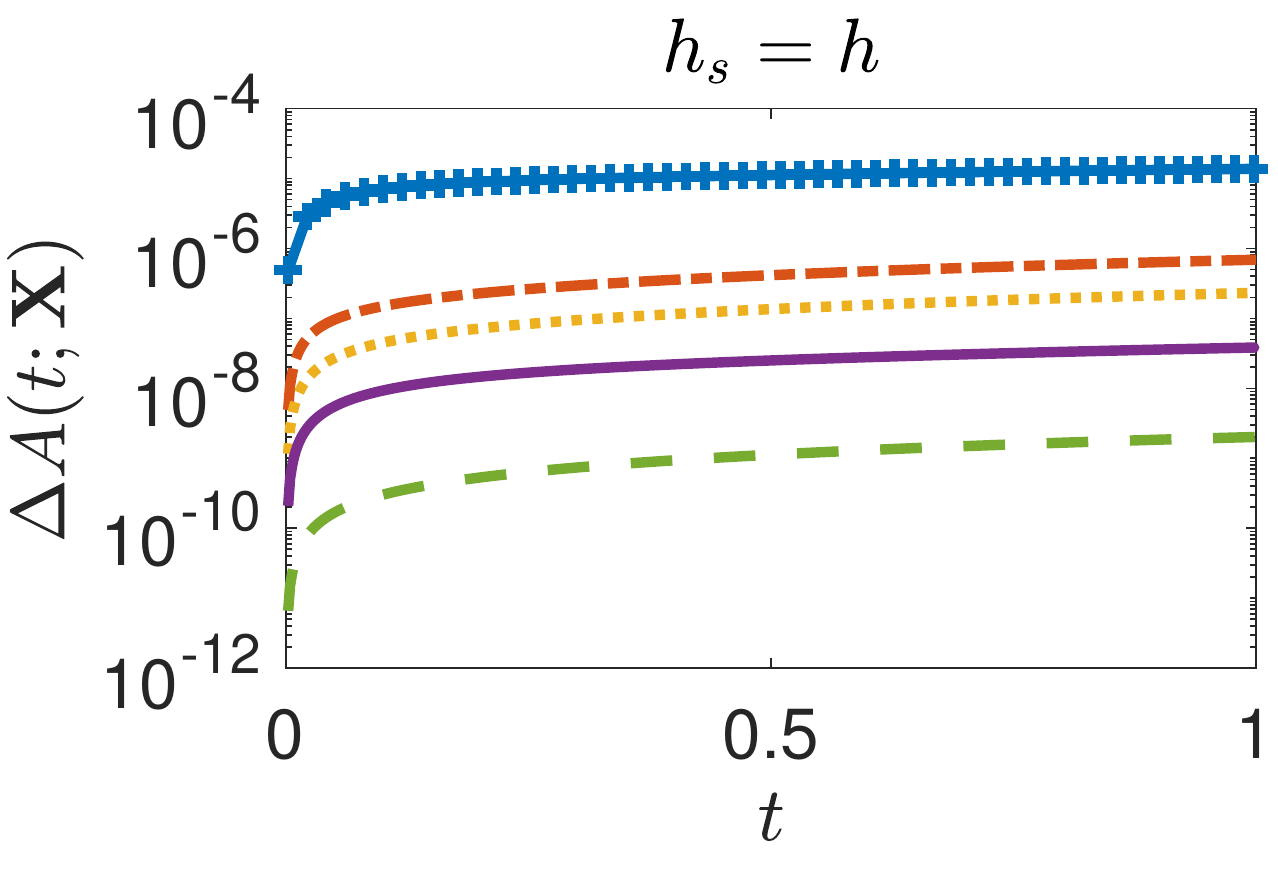} \label{fig:IBDFcircle_c}}
\caption{Normalized area errors of the pressurized circular membrane {(relative to the initial area, see \eqref{areaerr})} simulated by IBMAC (top panel) and DFIB (bottom panel) with the IB kernels:  $\stndfour \in \mathscr{C}^1$, $\bspline \in \mathscr{C}^2$, $\newfive \in \mathscr{C}^3$, $\newsix \in \mathscr{C}^3$ and $\bsplinesix \in \mathscr{C}^4$, and with Lagrangian marker spacings  $h_s \in \left\{4h, 2h, h\right\}$ indicated above each figure panel.}
\label{fig:circleCoarse}
\end{figure}

\begin{figure}
\subfloat[][]{\includegraphics[width =.33\linewidth]{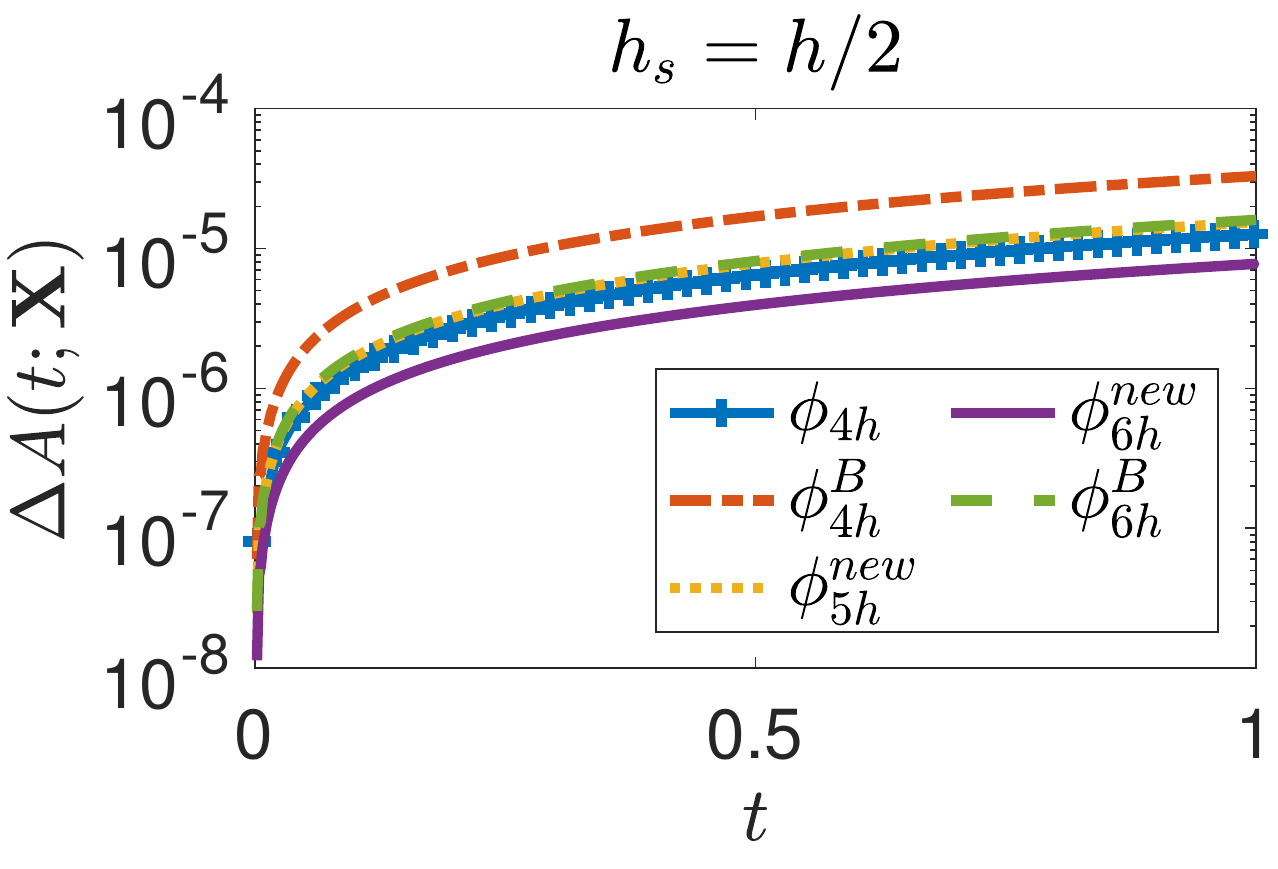} \label{fig:IMBACcircle_d} }
\subfloat[][]{\includegraphics[width =.33\linewidth]{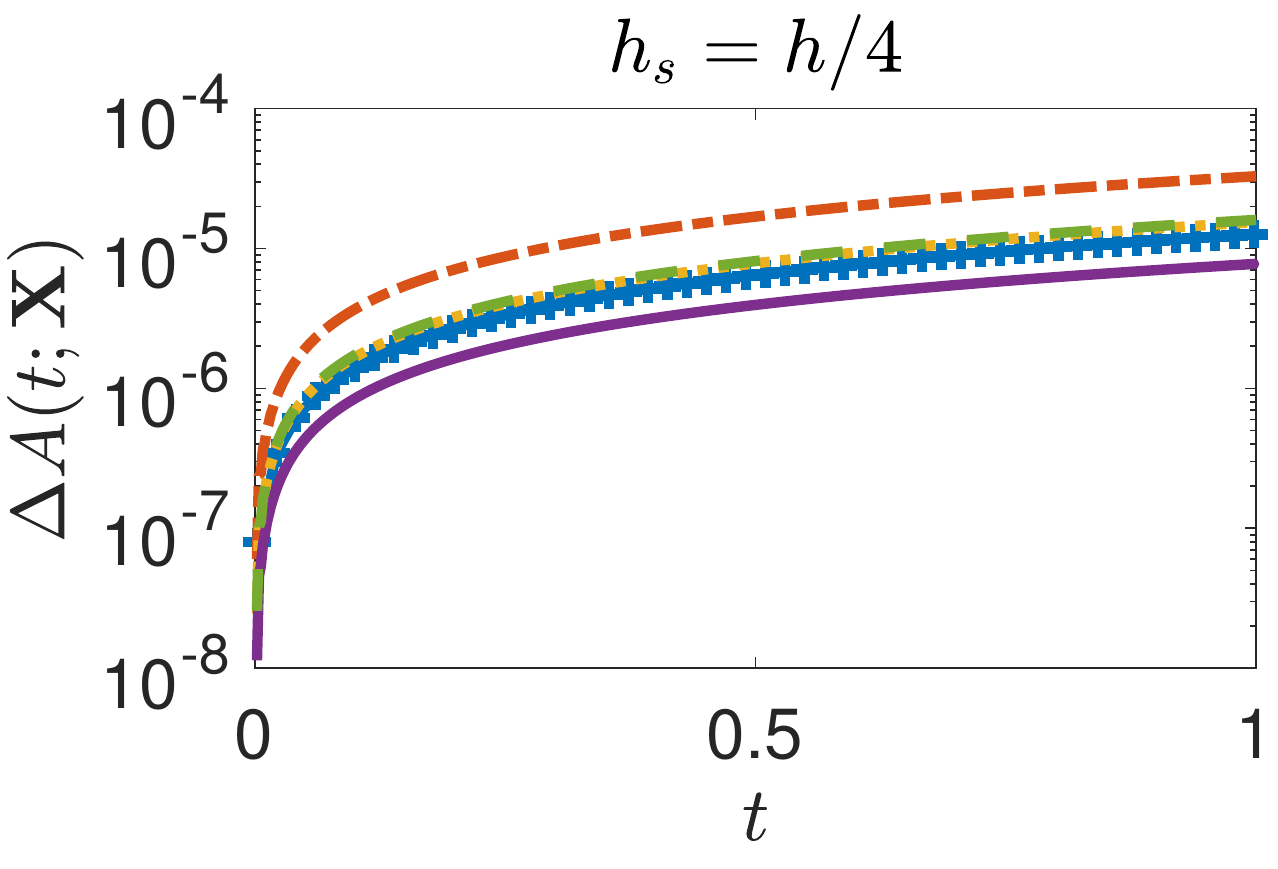} \label{fig:IBMACcircle_e} }
\subfloat[][]{\includegraphics[width =.33\linewidth]{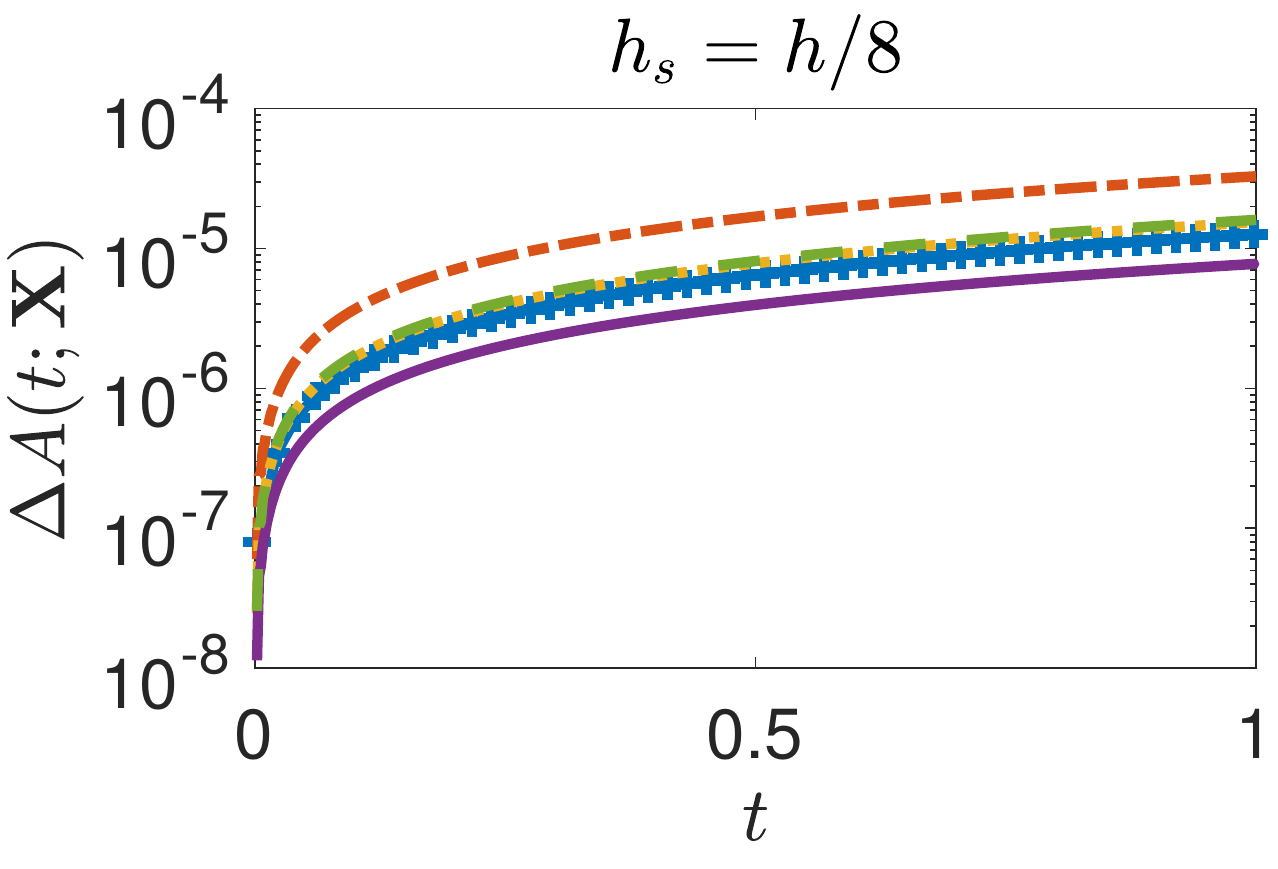} \label{fig:IBMACcircle_f}} \\
\subfloat[][]{\includegraphics[width =.33\linewidth]{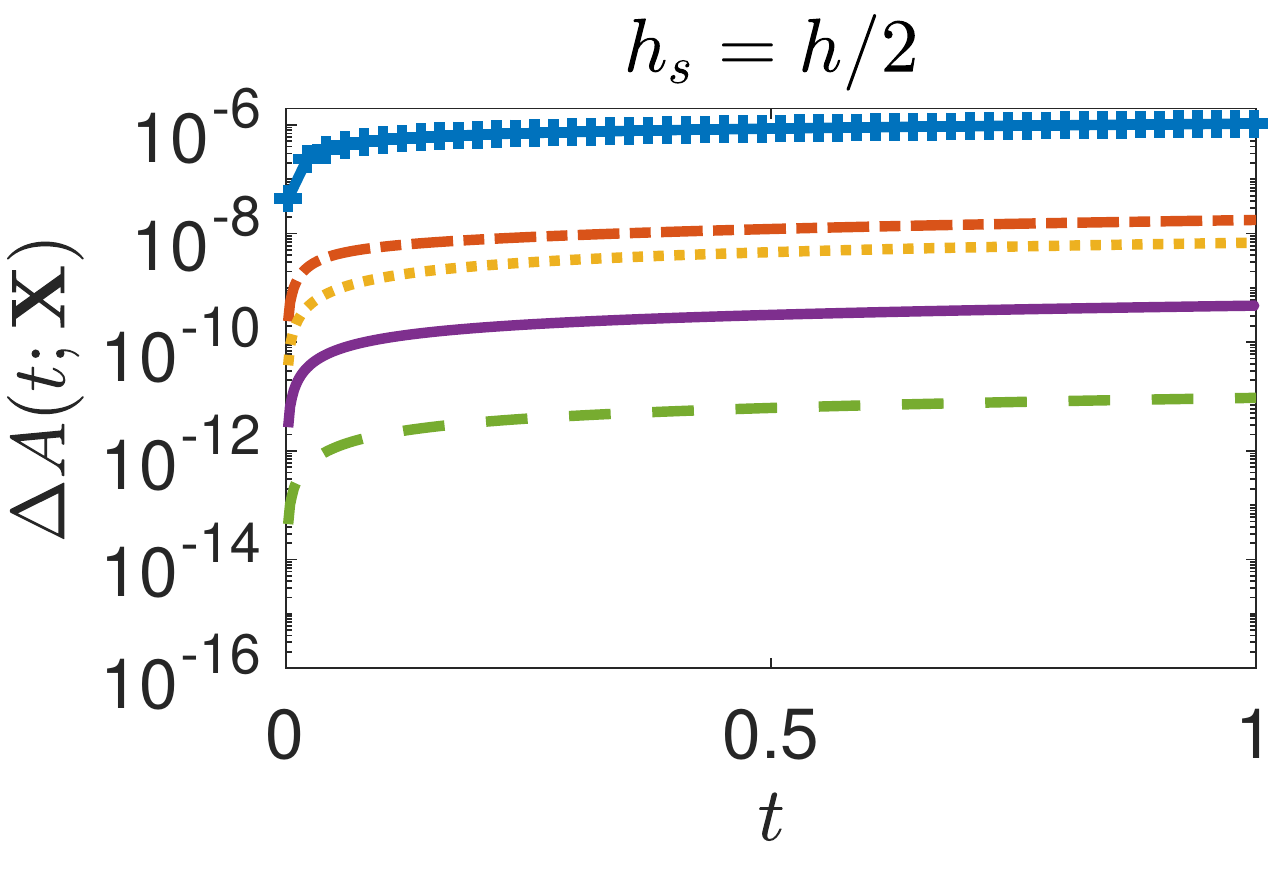} \label{fig:IBDFcircle_d}}
\subfloat[][]{\includegraphics[width =.33\linewidth]{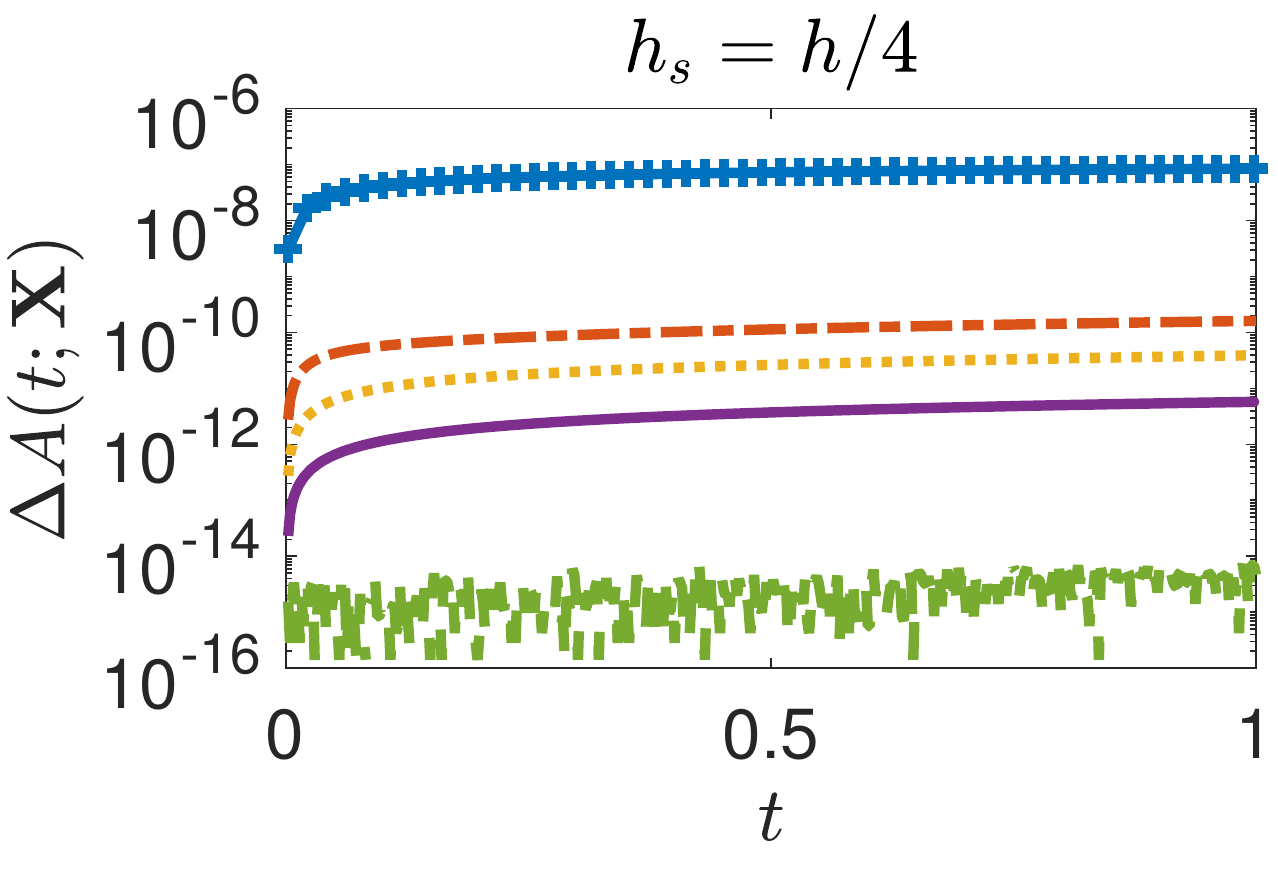} \label{fig:IBDFcircle_e}}
\subfloat[][]{\includegraphics[width =.33\linewidth]{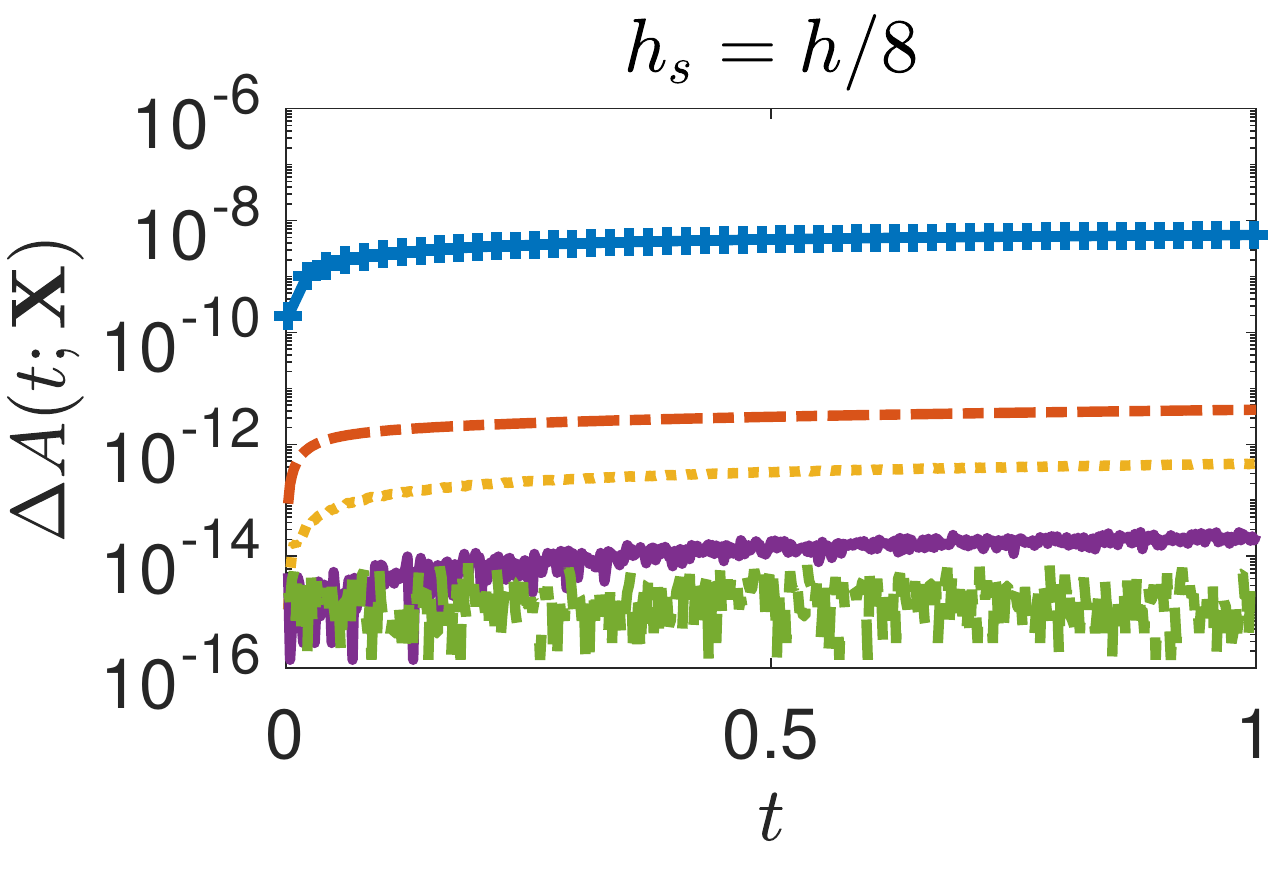} \label{fig:IBDFcircle_f}}
\caption{{Normalized area errors of the pressurized circular membrane (relative to the initial area, see \eqref{areaerr})} simulated by IBMAC (top panel) and DFIB (bottom panel) with the IB kernels:  $\stndfour \in \mathscr{C}^1$, $\bspline \in \mathscr{C}^2$, $\newfive \in \mathscr{C}^3$, $\newsix \in \mathscr{C}^3$ and $\bsplinesix \in \mathscr{C}^4$, and with Lagrangian marker spacings,  $h_s \in \left\{\frac{h}{2}, \frac{h}{4}, \frac{h}{8}\right\}$ indicated above each figure panel. As the Lagrangian mesh is refined, area errors for DFIB keep improving, even approaching the machine precision for $\bsplinesix$ at $h_s = \frac{h}{4}, \frac{h}{8}$ and for $\newsix$ at $h_s = \frac{h}{8}$. }
\label{fig:circleFine}
\end{figure}

We define the normalized area error with respect to the initial configuration
\begin{equation}
 \Delta A(t ; \mbf{X}) := \frac{|A(t;\mbf{X}) - A(0;\mbf{X})|}{A(0;\mbf{X})},
 \label{areaerr}
\end{equation}
where {the area enclosed by the Lagrangian markers $A(t;\mbf{X})$ is approximated by the area of the polygon formed by the Lagrangian markers $\mbf{X}=\{\mbf{X}_1,\dots,\mbf{X}_{M}\}$ at time $t$.}
\autoref{fig:circleCoarse} and \autoref{fig:circleFine} show the normalized area errors defined by \Cref{areaerr} for DFIB and IBMAC with  
different choices of the IB kernels: $\stndfour \in \mathscr{C}^1$, $\bspline \in \mathscr{C}^2$, $\newfive \in \mathscr{C}^3$, $\newsix \in \mathscr{C}^3$ and $\bsplinesix \in \mathscr{C}^4$. For the coarse Lagrangian marker spacings, for example, when $h_s = 2h,\, 4h$, the area errors for IBMAC and DFIB have similar orders of magnitude (compare \autoref{fig:IBMACcircle_a}, \ref{fig:IBMACcircle_b} to \autoref{fig:IBDFcircle_a}, \ref{fig:IBDFcircle_b}). As the Lagrangian marker spacing is reduced from $2h$ to $h$, we see a decrease in $\Delta A(t; \mbf{X})$ for IBMAC by approximately a factor of 10 (see \autoref{fig:IBMACcircle_b}, \ref{fig:IBMACcircle_c}) for all the IB kernels we consider in this set of tests. In contrast, the area errors for DFIB improve by at least a factor of $10^3$ for the IB kernels that are at least $\mathscr{C}^2$ (see \autoref{fig:IBDFcircle_b}, \ref{fig:IBDFcircle_c}), and in the best scenario,  $\Delta A(t ; \mbf{X})$ for  $\bsplinesix$ decreases from $10^{-4}$ to $10^{-9}$. Moreover, as the Lagrangian mesh is refined from $h$ to $\frac{h}{8}$, area errors for DFIB keep improving, even approaching the {machine epsilon in double precision} for $\bsplinesix$ at $h_s = \frac{h}{4}, \frac{h}{8}$ and for $\newsix$ at $h_s = \frac{h}{8}$ (see \autoref{fig:IBDFcircle_e}, \ref{fig:IBDFcircle_f}). For a moderate Lagrangian marker spacing, such as $h_s = h$ and $\frac{h}{2}$, area errors for DFIB are several orders of magnitude smaller than those of IBMAC. On the other hand, area errors for IBMAC stop improving around $10^{-5}$ for $h_s \leq h$, no matter how densely the Lagrangian mesh is refined (see \autoref{fig:IBMACcircle_c}, \ref{fig:IBMACcircle_f}). We remark that the smoothness of the IB kernel appears to play an important role in volume conservation of DFIB. In this study DFIB achieves the best volume conservation result for $h_s \leq h$ with $\bsplinesix$, and this kernel also has the highest regularity of the kernel functions considered in this work. 

\begin{figure}
\includegraphics[width =\linewidth]{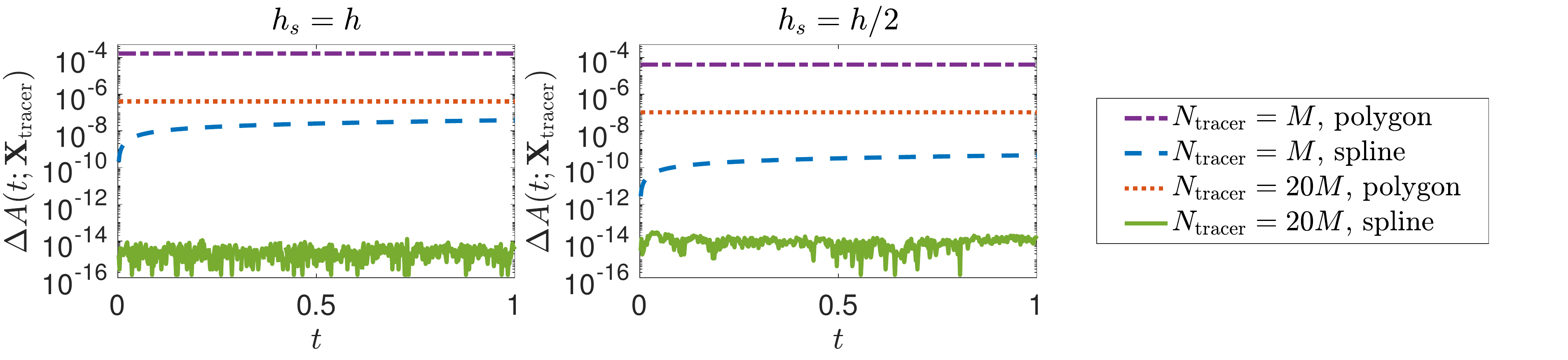} 
\caption{{Normalized area errors of the interface enclosed by the tracers that move passively with the interpolated velocity of the DFIB method (relative to the true area of the circle, see \Cref{areaerr_tracer}). The initial configuration of the interface is given by \Cref{membrane_circle}, and $\newsix$ is used for this computation. From left to right, the Lagrangian marker spacing is $h_s \in \left\{h, \frac{h}{2} \right\}$, as shown in each figure panel. In each case, the area error enclosed by the tracers is computed for tracer resolution $N_{\text{tracer}}=M \text{ and } 20M$ in two ways: (1) by the area of the polygon formed by the tracers, and (2) by the exact integration of the cubic spline interpolant of the tracer interface.} }
\label{fig:IBDFcircletracer}
\end{figure}

{The area errors of DFIB shown in \autoref{fig:circleCoarse} and \autoref{fig:circleFine} can be attributed to two sources of error. The first source of error is the time-stepping error from the temporal integrator, which is relatively small in the quasi-static circle test. The second source of area loss comes from discretizing the continuous curve (circle) as a polygon whose vertices are the IB markers. This kind of error can be substantially reduced by using a high-order representation of the interface, such as a periodic cubic spline. 
We define a normalized area error with respect to the true initial area of the interface $A_{\text{true}}$ by using the tracers,
\begin{equation}
 \Delta A(t ; \mbf{X}_{\text{tracer}}) := \frac{|A(t; \mbf{X}_{\text{tracer}} ) - A_{\text{true}}|}{A_{\text{true}}},
 \label{areaerr_tracer}
\end{equation}
where we compute $A(t; \mbf{X}_{\text{tracer}} )$ via exact integration of the cubic spline interpolant. 
\autoref{fig:IBDFcircletracer} shows that the area enclosed by the passive tracers using the cubic spline approximation is far more accurately preserved than the polygonal approximation. Furthermore, the area error approaches zero as the discretization of the tracer interface is refined, even if the marker and grid spacings are held fixed, as shown in \autoref{fig:IBDFcircletracer}. Indeed, the area of the spline interpolant through the tracers of spacing $h_s/20=h/20$ is conserved to the machine epsilon in double precision.}

\begin{figure}
\subfloat[][DFIB]{\includegraphics[width =.48\linewidth]{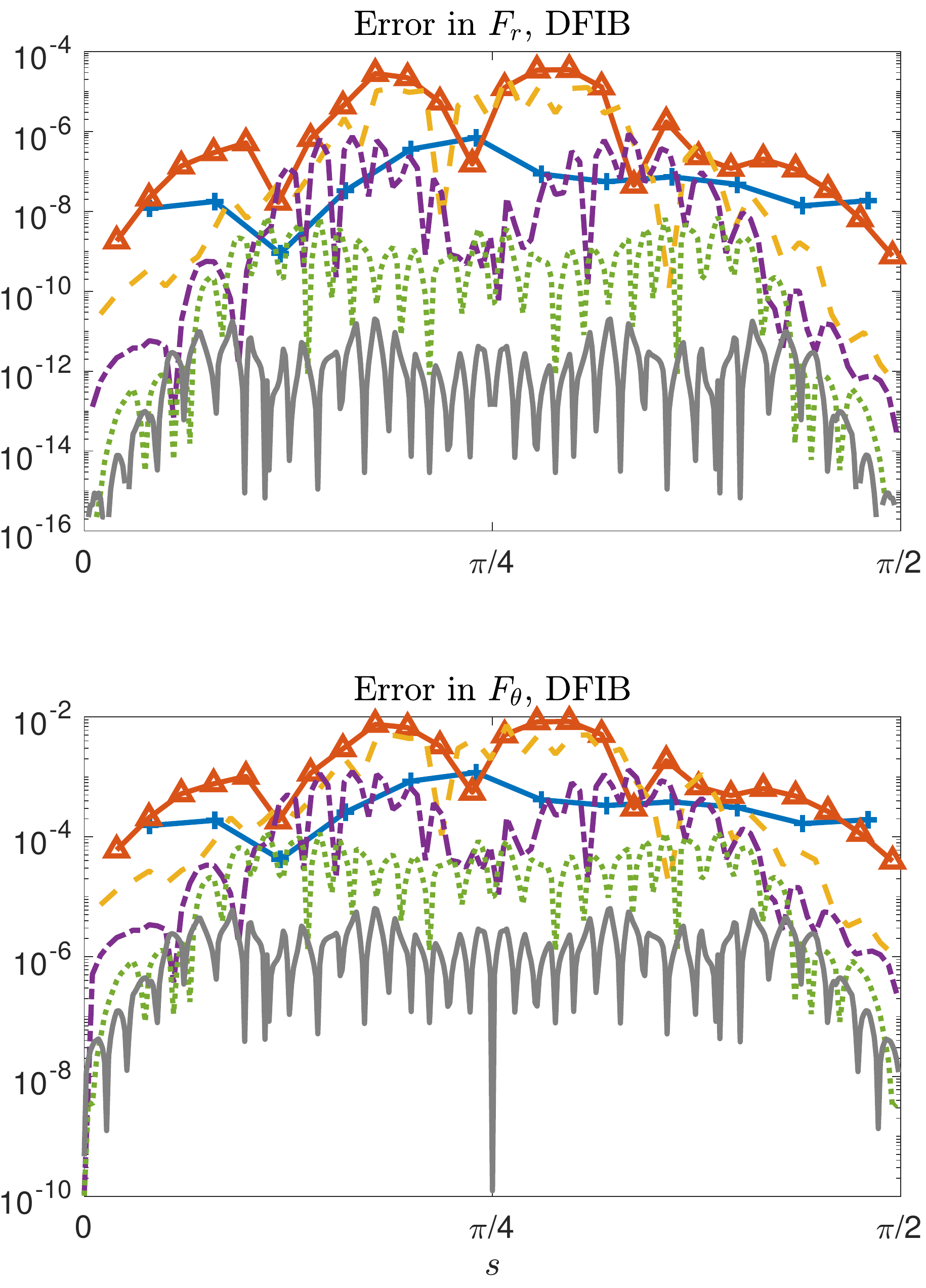} \label{fig:force_a} } \hspace{1em}
\subfloat[][IBMAC]{\includegraphics[width =.48\linewidth]{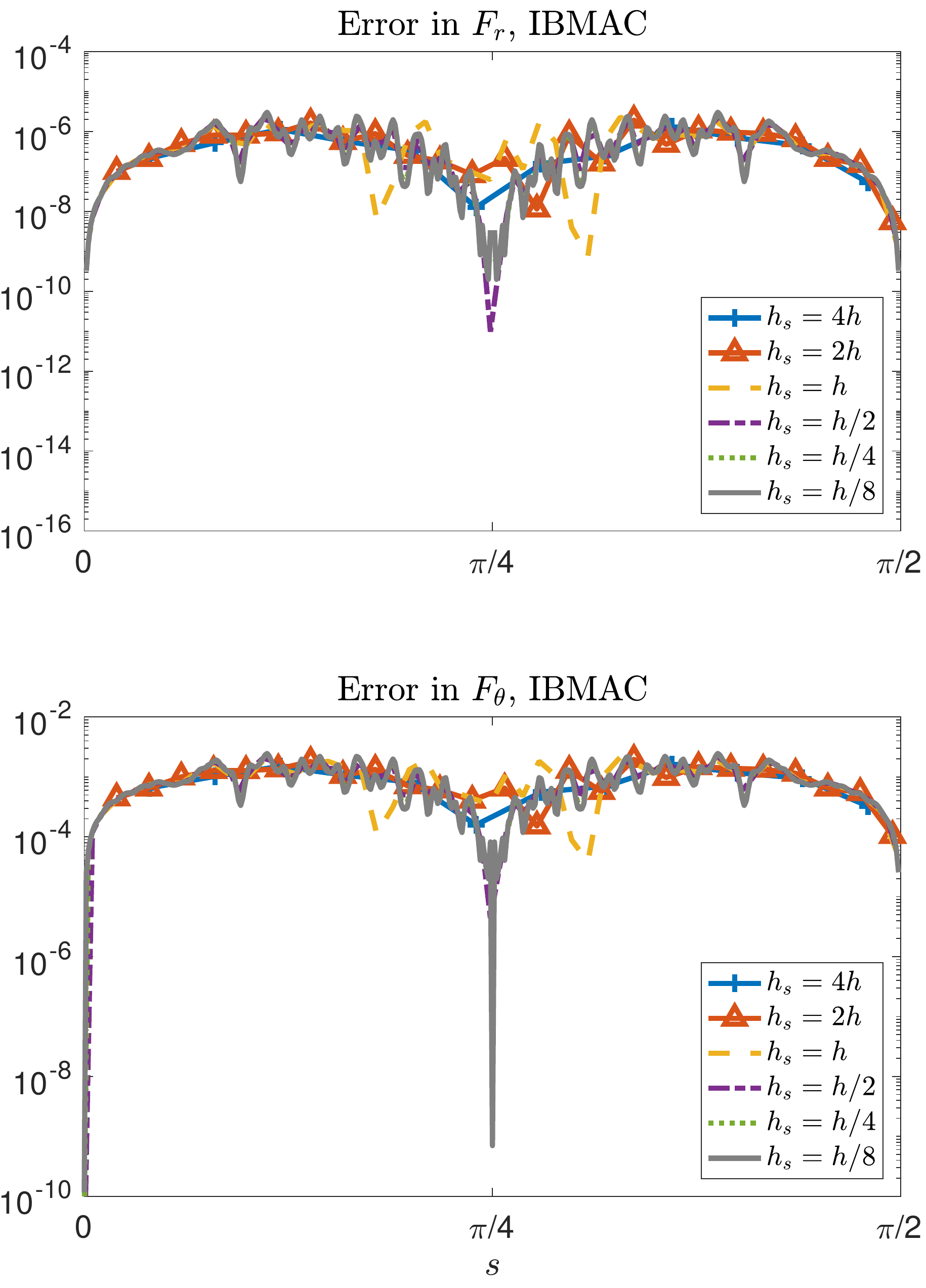} \label{fig:force_b}} 
\caption{{Normalized errors of the normal $F_r$ (top panels) and tangential $F_{\theta}$ (bottom panels) components of the Lagrangian force $\mbf{F}(s, t= 1)$ of the circular membrane for $s \in [0,\frac{\pi}{2}]$. The computations are performed using DFIB (left panels) and IBMAC (right panels) with $\newsix$, and with Lagrangian marker spacings $h_s \in \left\{ 4h, 2h, h, \frac{h}{2}, \frac{h}{4}, \frac{h}{8} \right\}$. }}
\label{fig:DFIBvsIBMAC_force}
\end{figure}

{It is well-known that the traditional IB method produces non-smooth surface tractions, and a number of improvements have been proposed \cite{IBStresses_Colonius,Tractions_Fauci,SmoothingDelta_IBM,Strychalski2016,LiskaColonius2017}. Somewhat unexpectedly, the divergence-free force spreading used in our DFIB method proposed here offers smoother and more accurate tractions without any post-processing such as filtering \cite{IBStresses_Colonius}. This is inherently linked to the reduced spurious flows compared to traditional methods \cite{Strychalski2016}. In \autoref{fig:DFIBvsIBMAC_force}, we compare the errors of the tangential and normal components of $\mbf{F}(s, t = 1)$ for DFIB and IBMAC for the quasi-static circle problem. We observe that the DFIB method dramatically improves the accuracy of Lagrangian forces by only refining the Lagrangian mesh, keeping the Eulerian grid fixed. By contrast, in the IBMAC method, the tractions do not improve as the Lagrangian grid is refined.}


\subsection{A thin elastic membrane with parametric resonance in 2D} 

\begin{figure}
\centering
\subfloat[][]{\includegraphics[width = .75\linewidth]{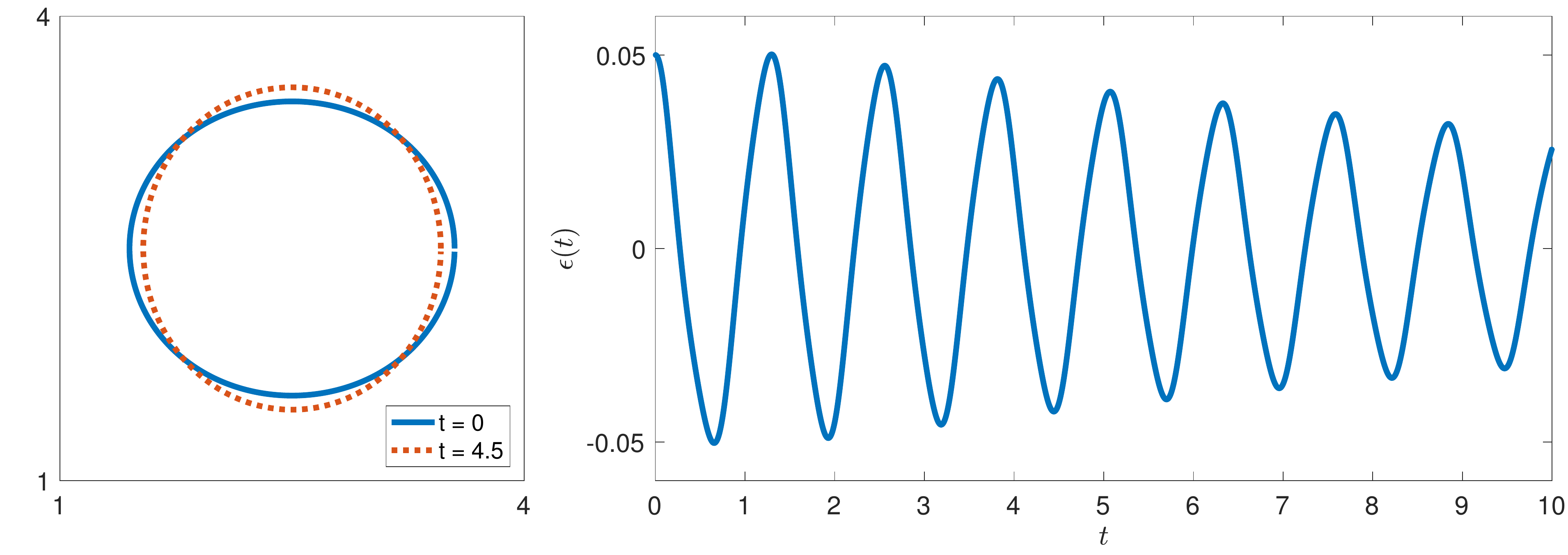} \label{fig:IBDFpara_stable}} \\
\subfloat[][]{\includegraphics[width = .75\linewidth]{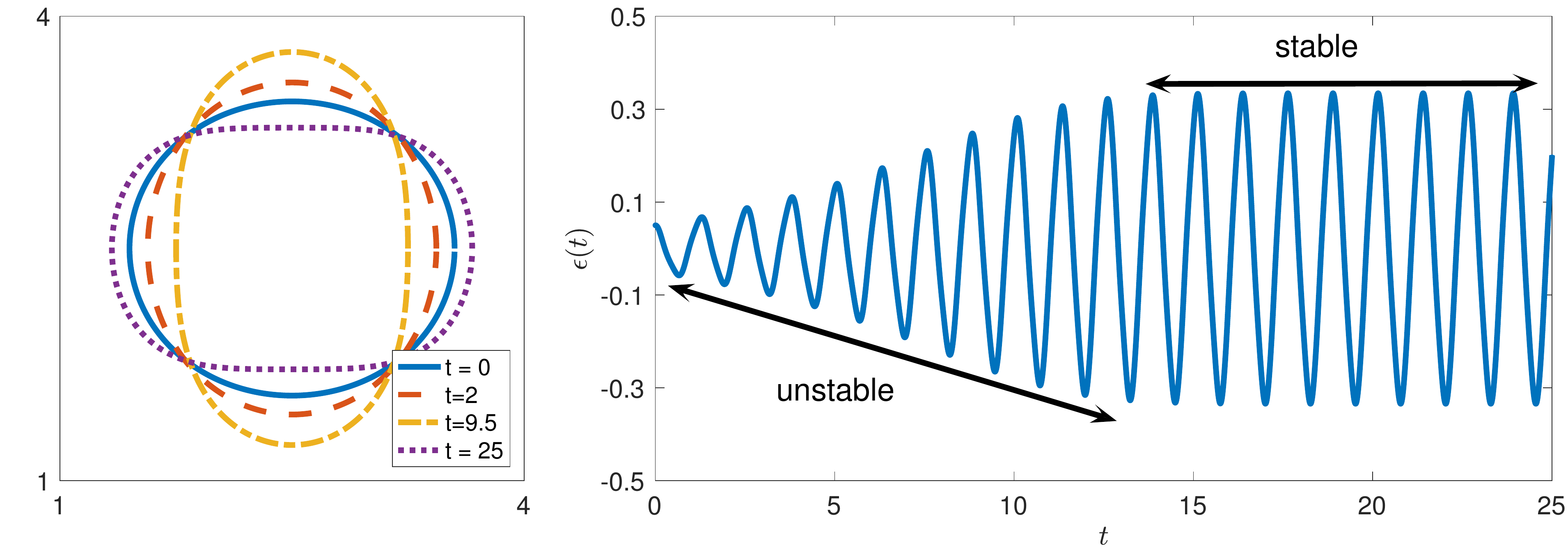}\label{fig:IBDFpara_unstable}} 
\caption{Left panel: snapshots of the 2D membrane with parametric resonance. Right panel: the time-dependent amplitude $\epsilon(t)$ of the perturbed mode in \Cref{soln_ansatz}. (a) Damped oscillation (b) Growing oscillation. }
\label{fig:IBDFpara_figs}
\end{figure}

In many biological applications, the immersed structure is an active material, interacting dynamically with the surrounding fluid and generating time-dependent motion. It has been reported that the simulation of active fluid-structure interactions using the conventional IB method may suffer from significant loss in the volume enclosed by {the} structure \cite{Peskin1993_IBmodified}. A simple prototype problem for active fluid-structure interaction is a thin elastic membrane that dynamically evolves in a fluid in response to elastic forcing with periodic variation in the stiffness parameter \cite{Cortez2004_2Dparametric, Ko2012}, that is, 
\begin{equation}
 \mbfit{F}(s,t) =  \kappa(t) \frac{\partial^2 {\mbfit{\mathcal{X}}}}{\partial s^2}, \label{force_parametric}
\end{equation}
where $\kappa(t)$ is a periodic time-dependent stiffness coefficient of the form
\begin{equation}
 \kappa(t) = K_c (1+2\tau \sin(\omega_0 t)). \label{para_stiffness}
\end{equation}
It is quite remarkable that such a purely temporal parameter variation can result in the emergence of spatial patterns, but that is indeed the case. We assume that the immersed structure is initially in a configuration that has a small-amplitude perturbation from a circle of radius $R$,
\begin{equation}
 {\mbfit{\mathcal{X}}}(s,0)  = R(1+\epsilon_0 \cos(ps)) \ \mbfit{\hat{r}}(s),
\end{equation}
where $\mbfit{\hat{r}}(s)$ denotes the position vector pointing radially from the origin. For certain choices of parameters, the perturbed mode in the initial configuration may resonate with the driving frequency $\omega_0$ in the periodic forcing, leading to large-amplitude oscillatory motion in the membrane. The stability of the parametric resonance has been studied in the IB framework using Floquet linear stability analysis for a thin elastic membrane in 2D \cite{Cortez2004_2Dparametric, Ko2012}, and recently for an elastic shell in 3D \cite{Ko2016_3Dparametric}. Motivated by the linear stability analysis of \cite{Cortez2004_2Dparametric, Ko2012}, we consider two sets of parameters listed in \autoref{table:IBparametric_parameters} for our simulations. The first set of parameters with $\tau = 0.4$ leads to a stable configuration in which the membrane undergoes damped oscillations (\autoref{fig:IBDFpara_stable}), and the second set with $\tau = 0.5$ leads to an unstable configuration in which the membrane oscillates with growing amplitude (\autoref{fig:IBDFpara_unstable}). 

\renewcommand{\arraystretch}{1.3}
\begin{table}
\centering
 \begin{tabular}{c c c c c c c c l}
\hline\hline
$\rho$ & $\mu$ & $L$ & $R$ & $K_c$ & $\omega_0$ & $p$ & $\epsilon_0$ & $\tau$  \\ 
\hline
\multirow{2}{*}{1} & \multirow{2}{*}{0.15}  & \multirow{2}{*}{5} & \multirow{2}{*}{1} & \multirow{2}{*}{10} & \multirow{2}{*}{10} & \multirow{2}{*}{2} &  \multirow{2}{*}{0.05}  &  0.4 (damped oscillation)  \\
& & & & & & & & 0.5 (growing oscillation) \\
\hline \hline
\end{tabular}
 \caption{Parameters used to simulate the motion of the 2D membrane with parametric resonance.}
\label{table:IBparametric_parameters}
\end{table}

The computational domain $\Omega = [0,L]^2$ is discretized by a $128 \times 128$ uniform Cartesian grid with meshwidth $h = \frac{L}{128}$. The number of Lagrangian markers is determined so that the distance between the Lagrangian markers is $h_s \approx \frac{h}{2}$ in the initial configuration. The discretization of the Lagrangian force density \Cref{force_parametric} is constructed in the same way as \Cref{springforce_2nd}. The time step size is $\Delta t = \frac{h}{10}$ to ensure the stability of computation. On the left panel of \autoref{fig:IBDFpara_figs} we show snapshots of the membrane configuration for each case, and on the right panel we plot the time-dependent amplitude $\epsilon(t)$ of the ansatz
\begin{equation}
{\mbfit{\mathcal{X}}}(s,t) = R(1+\epsilon(t)\cos(ps)) \ \mbfit{\hat{r}}(s) \label{soln_ansatz}
\end{equation}
by applying the FFT to the  Lagrangian marker positions $\mbf{X}$.
In the case of growing oscillation (\autoref{fig:IBDFpara_unstable}), the amplitude  of the perturbed mode increases from $0.05$ to $0.3$ until nonlinearities eventually stabilize the growing mode and the membrane starts to oscillate at a fixed amplitude.

\begin{figure}
\centering
\includegraphics[width =\linewidth]{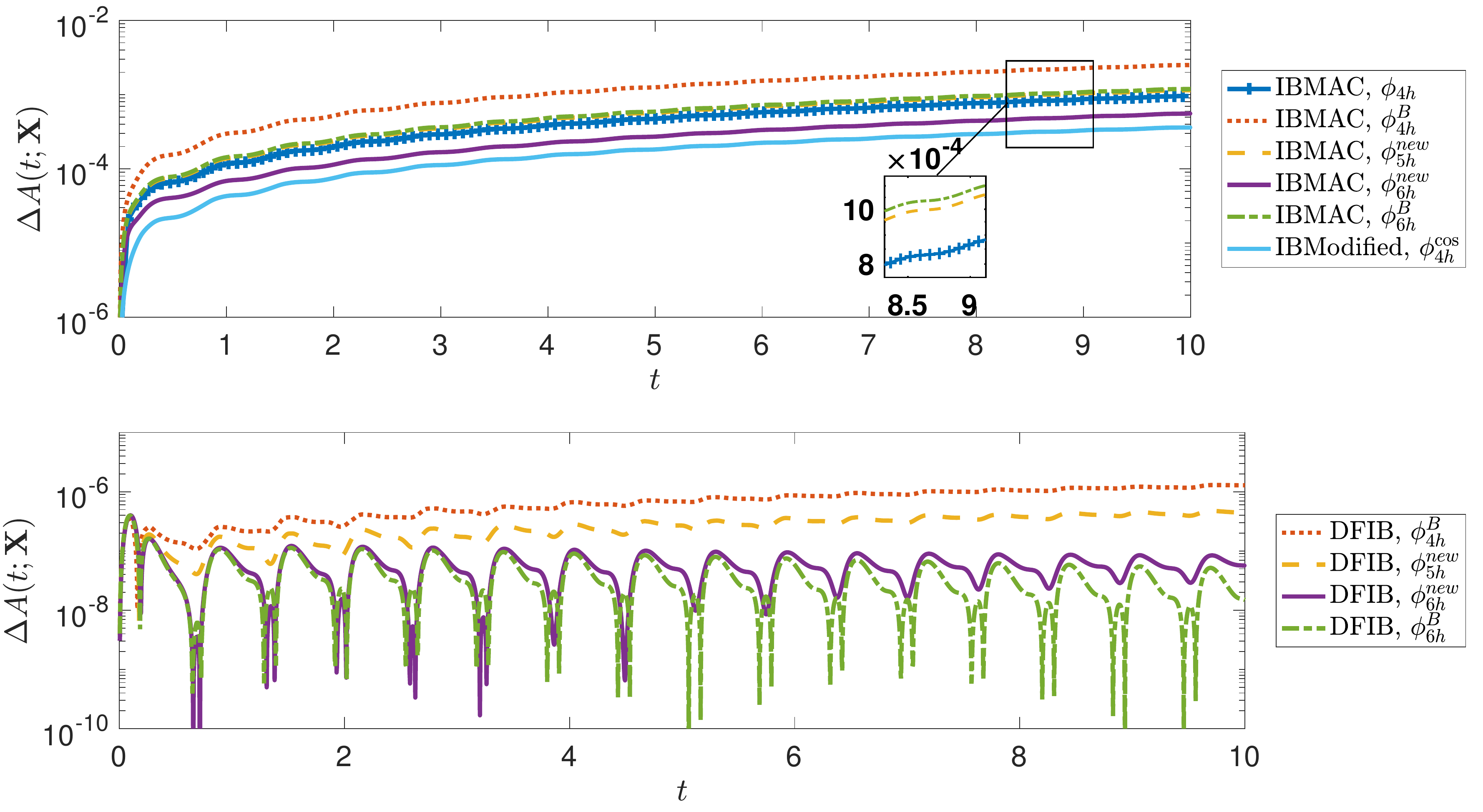} \\
\includegraphics[width =\linewidth]{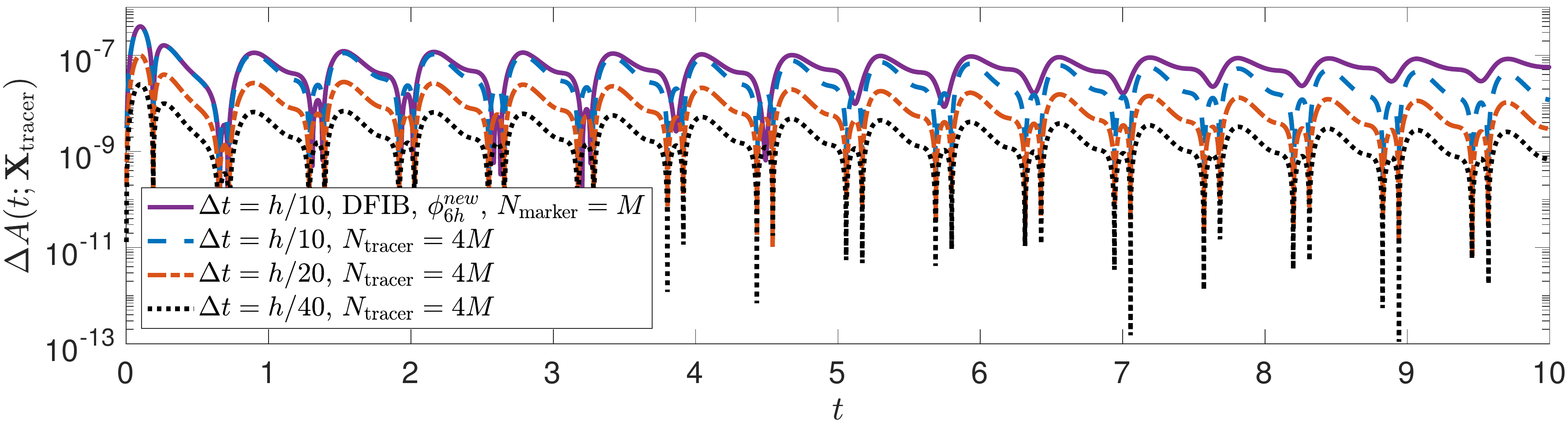}
\caption{{Top and middle panels:} normalized area errors $\Delta A(t; \mbf{X})$ {(relative to the initial area, see \eqref{areaerr})} of the 2D parametric membrane undergoing damped oscillatory motion (corresponding to the motion shown in \autoref{fig:IBDFpara_stable}) are plotted on the semi-log scale. {The area enclosed by the markers is computed by the exact integration of a cubic spline approximation.} The computations are performed using IBMAC and DFIB with the IB kernels: $\stndfour \in \mathscr{C}^1$, $\bspline \in \mathscr{C}^2$, $\newfive \in \mathscr{C}^3$, $\newsix \in \mathscr{C}^3$ and $\bsplinesix \in \mathscr{C}^4$, and IBModified with $\cosfour \in \mathscr{C}^1$. The top panel shows area errors for IBMAC and IBModified, and the {middle} panel shows area errors for DFIB. {The bottom panel shows improvement in area errors (relative to the true area of the ellipse, see \eqref{areaerr_tracer}) enclosed by the interface of $N_{\text{tracer}} = 4 M$ tracers by reducing the time step size: $\Delta t \in \left\{\frac{h}{10}, \frac{h}{20}, \frac{h}{40}\right\}$. } }
\label{fig:arealossIBpara_stable}
\end{figure}
\begin{figure}
\centering
\includegraphics[width =\linewidth]{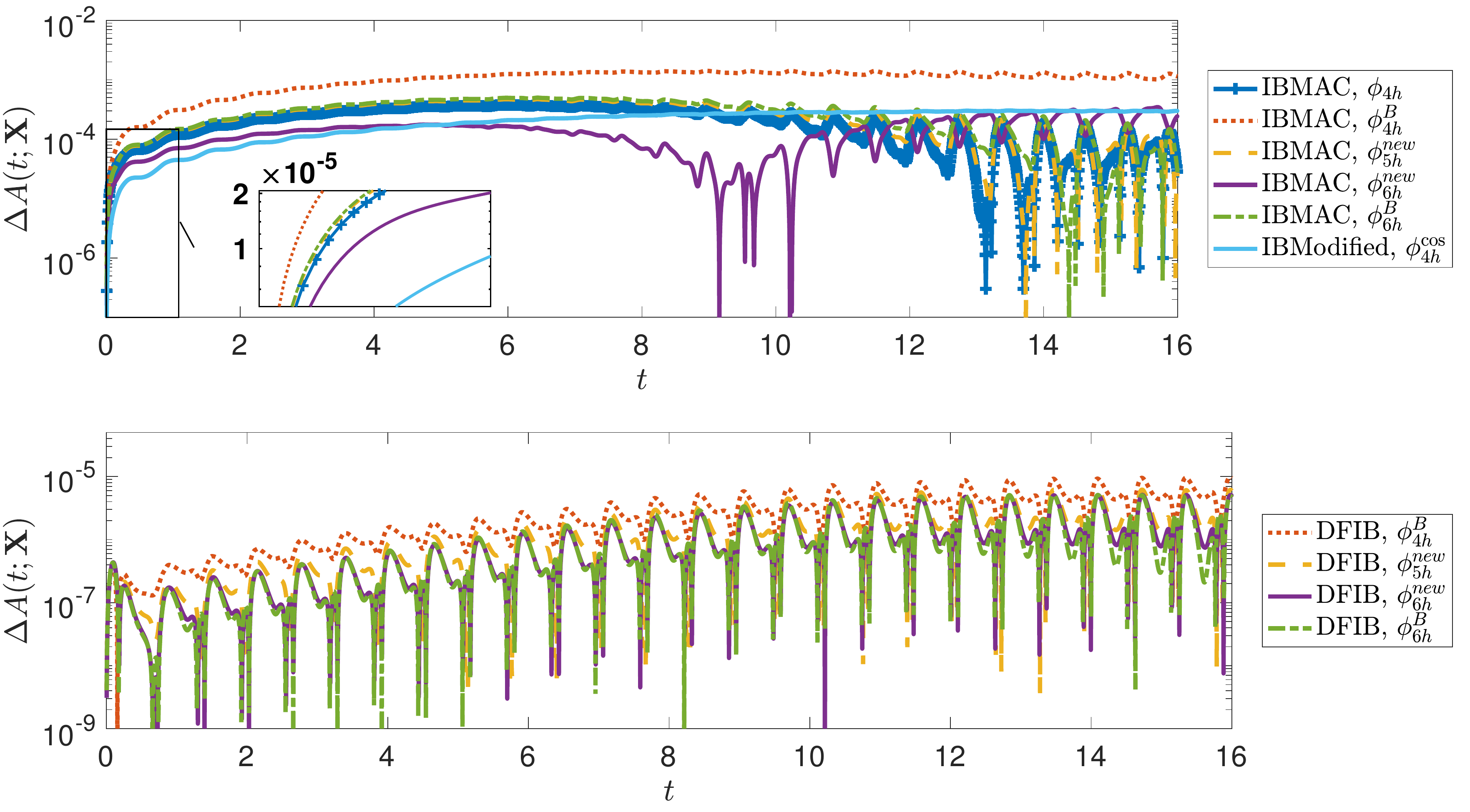} \\
\includegraphics[width =\linewidth]{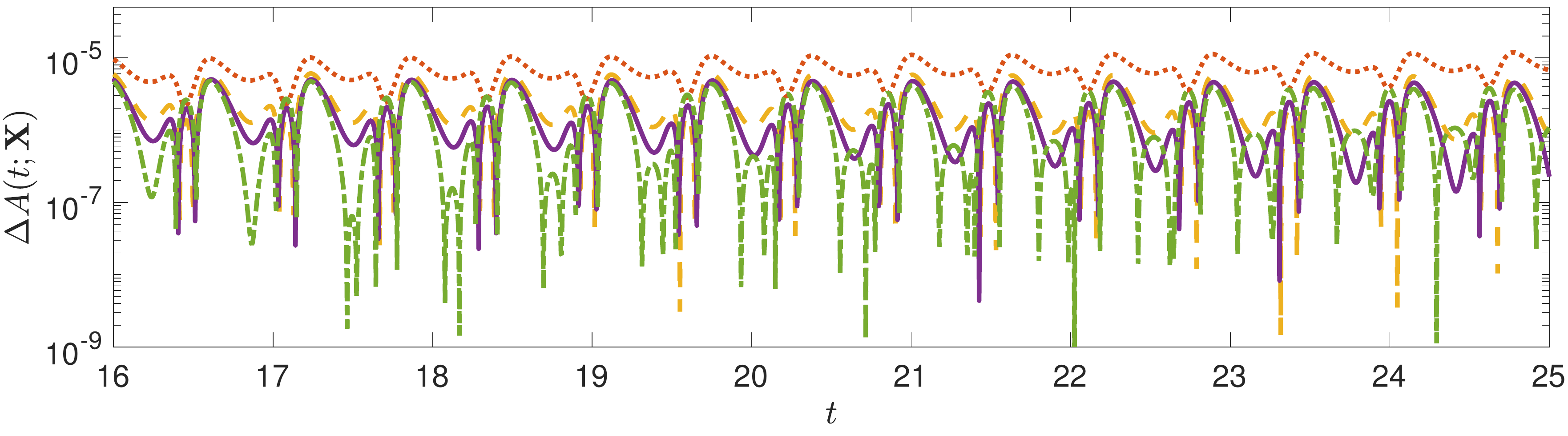} 
\caption{Normalized area errors $\Delta A(t; \mbf{X})$ of the 2D parametric membrane undergoing growing-amplitude oscillatory motion (corresponding to the motion shown in \autoref{fig:IBDFpara_unstable}) are plotted on the semi-log scale, as done in \autoref{fig:arealossIBpara_stable} for damped motion.  The computations are performed using IBMAC and DFIB with the IB kernels: $\stndfour \in \mathscr{C}^1$, $\bspline \in \mathscr{C}^2$, $\newfive \in \mathscr{C}^3$, $\newsix \in \mathscr{C}^3$ and $\bsplinesix \in \mathscr{C}^4$, and IBModified with $\cosfour \in \mathscr{C}^1$. The top panel shows area errors for IBMAC and IBModified, {the middle} panel shows the area errors for DFIB for $t = 0$ to 16, and the bottom panel extends the middle panel for $t = 16$ to 25.} 
\label{fig:arealossIBpara_unstable}
\end{figure}

We next give a direct comparison of area conservation of IBModified (with $\phi_{4h}^{\cos}$ \cite{Peskin1993_IBmodified}), IBMAC, and DFIB with  
the IB kernels $\stndfour \in \mathscr{C}^1$, $\bspline \in \mathscr{C}^2$, $\newfive \in \mathscr{C}^3$, $\newsix \in \mathscr{C}^3$ and $\bsplinesix \in \mathscr{C}^4$.
{In this test, the area enclosed the Lagrangian markers is computed by the cubic spline approximation discussed in \autoref{sec_spacing}.}
In \autoref{fig:arealossIBpara_stable} and \autoref{fig:arealossIBpara_unstable} we show time-dependent area errors enclosed by the parametric membrane for the damped-oscillation and the growing-amplitude cases  respectively. For the damped-oscillation case (\autoref{fig:arealossIBpara_stable}), we see that  area errors for DFIB are at least two orders of magnitude smaller than those of IBMAC and IBModified for IB kernels that are at least $\mathscr{C}^2$. The volume conservation of IBModified and IBMAC was not directly compared in the previous work \cite{Griffith2012_IBMACvolume}, but it was anticipated that they are similar. In our comparison, we find that IBModified is only slightly better than IBMAC in volume conservation, yet IBMAC is much simpler to use in practice. 

In this set of tests, the choice of IB kernel also plays a role in affecting area conservation. In particular, the area errors for DFIB using $\newsix$ and $\bsplinesix$ are  smaller than those of $\bspline$ and $\newfive$ by approximately one order of magnitude. Additionally, the error curves of DFIB with $\newsix$ and $\bsplinesix$  remain oscillating below $10^{-7}$ while apparent growth of error in time is observed  with $\bspline$ and $\newfive$, and in the other IB methods. {Unlike the quasi-static circle test in which the time-stepping error is negligible compared to the area loss due to moving a finite collection of Lagrangian markers,  the time-stepping error in this example can be observed by considering a dense collection of tracers with different time step sizes. With $N_{\text{tracer}} = 4M$, we first confirm that the area error of the tracer interface cannot be reduced by further including more tracers, but as we reduce the time step size $\Delta t \in \left\{\frac{h}{10}, \frac{h}{20}, \frac{h}{40}\right\}$, we observe an improvement in the area error, as shown in the bottom panel of \autoref{fig:arealossIBpara_stable}.}   
The improvements in area conservation of DFIB is consistently more than $10^4$ times over IBCollocated and about $10^3$ times over IBMAC. Similar results are obtained for the growing-amplitude case (see \autoref{fig:arealossIBpara_unstable}) except that the parametrically-unstable membrane has experienced some area loss due to the growing-amplitude oscillation before its motion is stabilized by the nonlinearities.

\subsection{A 3D thin elastic membrane with surface tension}
\begin{figure}
\centering
\subfloat[][]{\includegraphics[width = .32\linewidth]{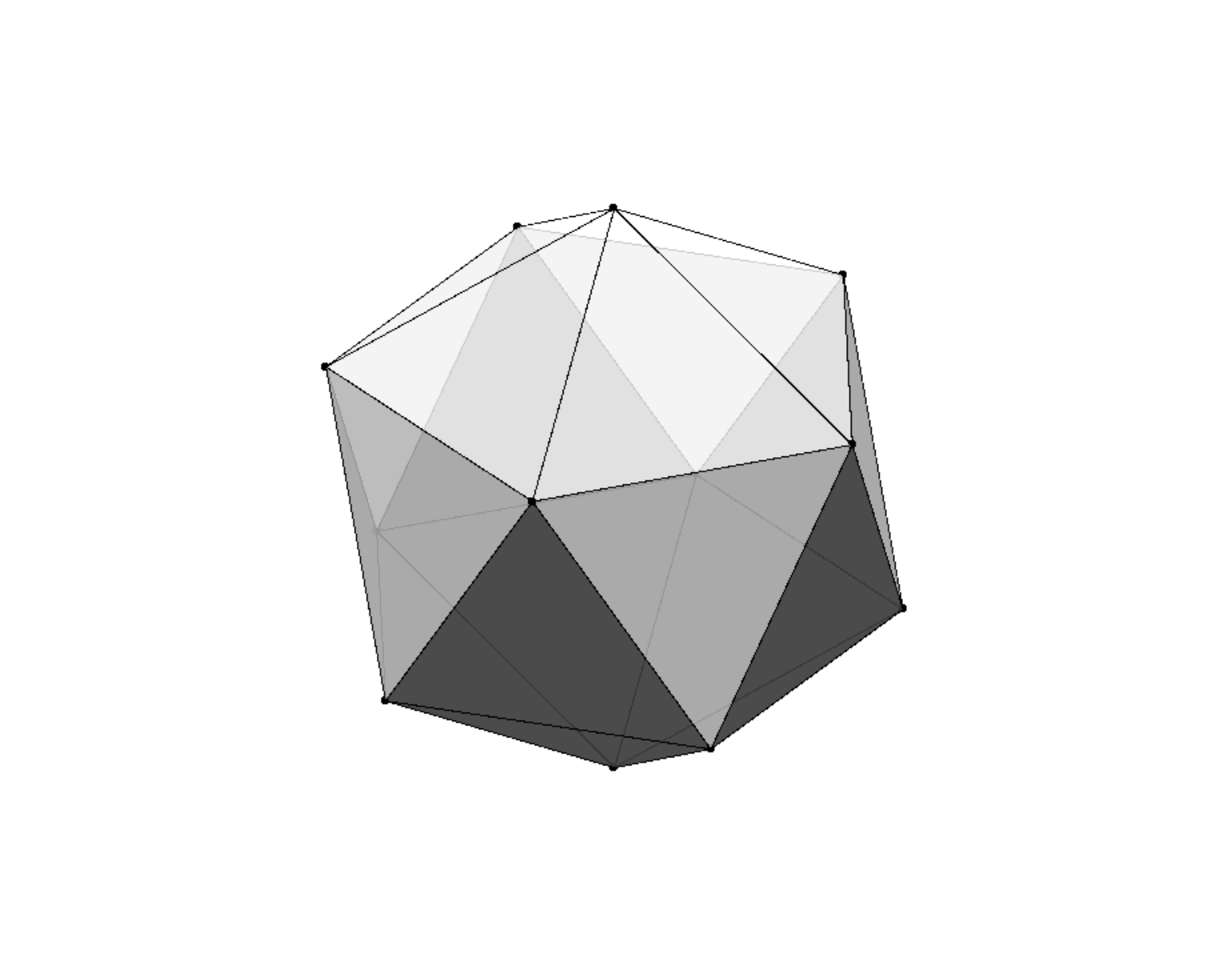}\label{fig:sphere0} } 
\subfloat[][]{\includegraphics[width = .32\linewidth]{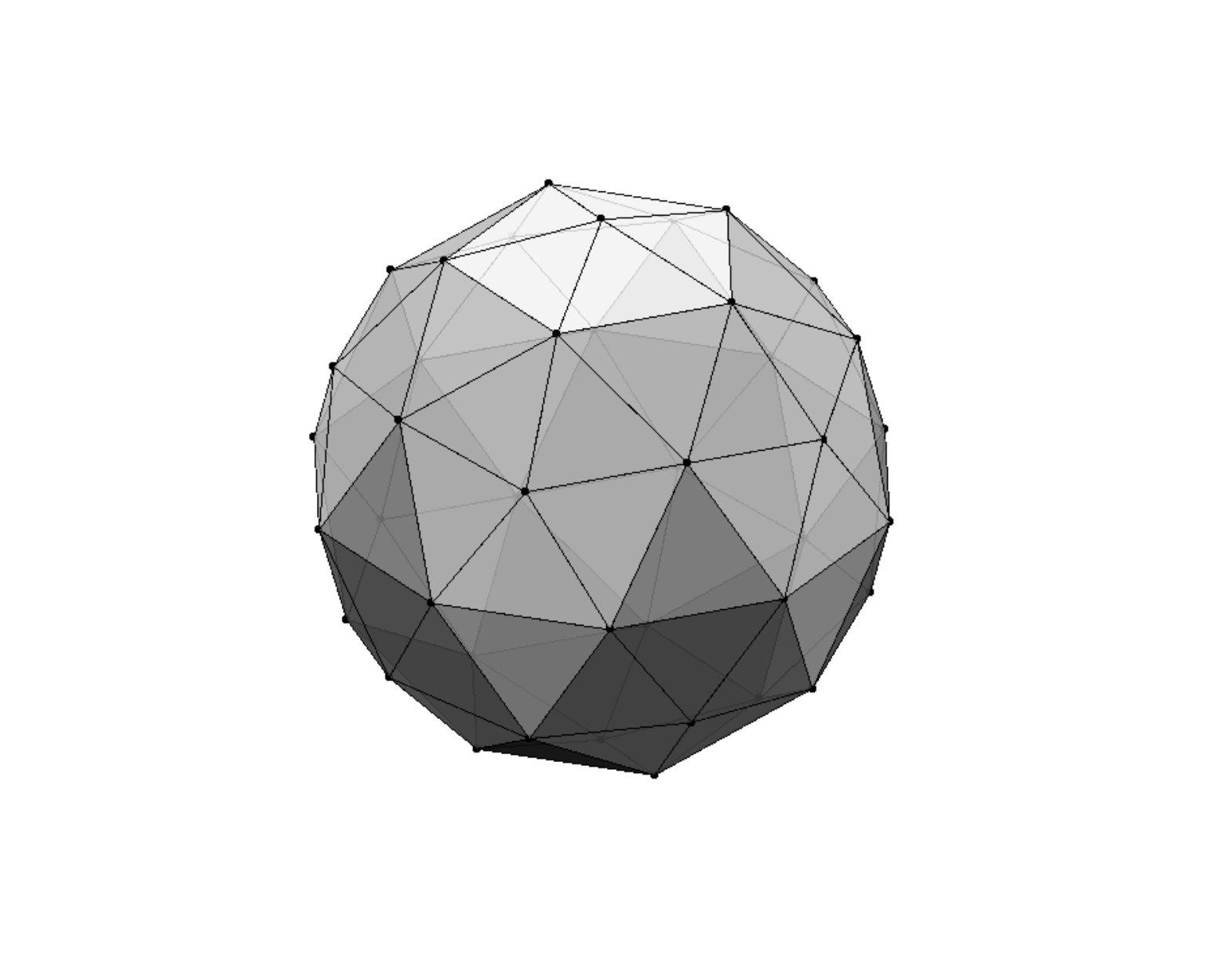} } 
\subfloat[][]{\includegraphics[width = .32\linewidth]{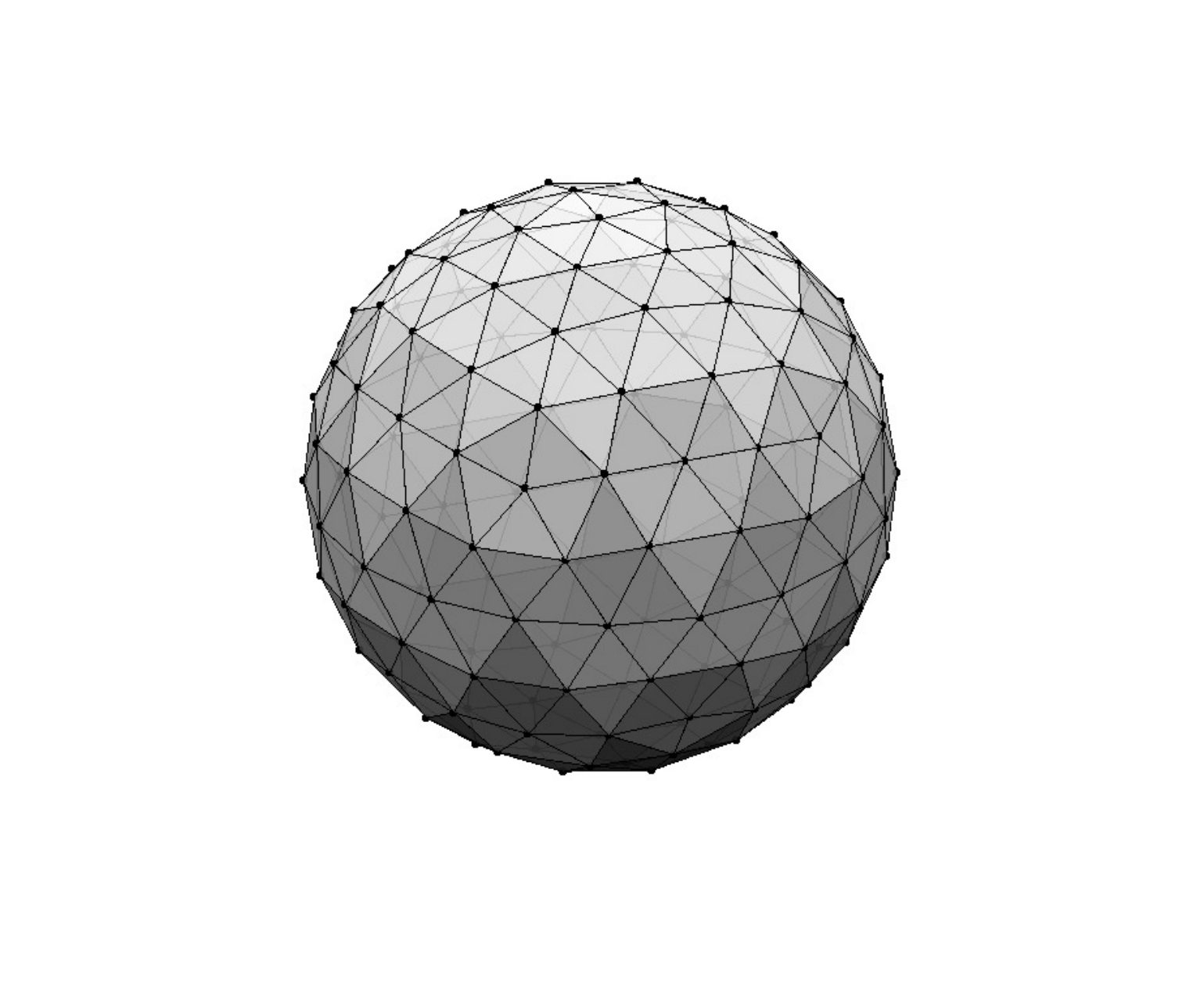} }
\caption{Triangulation of a spherical surface mesh via refinement of a regular icosahedron. (a) Regular icosahedron (b) Refined mesh after one level of refinement (c) Refined mesh after two levels of refinement. }
\label{fig:sphere_refine}
\end{figure}

In our final test problem, we examine volume conservation of the DFIB method by extending the surface tension problem to 3D. 
We consider in 3D a thin elastic membrane that is initially in its spherical equilibrium configuration. The spherical surface of the membrane is discretized by a triangulation consisting of approximately equilateral triangles with edge length approximately equal to $h_s$, constructed from successive refinement of a regular icosahedron by splitting each facet into four smaller equilateral triangles and projecting the vertices onto the sphere to form the refined mesh (see \autoref{fig:sphere_refine} for the first two levels of refinement). 
We use $\{\mbf{X}_1 , \mbf{X}_2, \dots, \mbf{X}_M \}$ and $\{\mbf{T}_1, \mbf{T}_2, \dots, \mbf{T}_P \}$ to denote the vertices (Lagrangian markers) and the triangular facets of the mesh respectively. The generalization of  discrete elastic energy functional of surface tension in 3D is the product of surface tension constant $\gamma$ (energy per unit area) and the total surface area of the triangular mesh \cite{Kim2014_foam3D}, that is, 
\begin{equation}
 E[\mbf{X}_1,\dots \mbf{X}_M] =  \gamma \sum_{p=1}^P |\mbf{T}_p| \, ,
\end{equation}
where $|\mbf{T}_p|$ is the area of the $p^{\text{th}}$ triangle. 
The Lagrangian force $\mbf{F}_k \Delta \mbf{s}$ at the $k^{\text{th}}$ vertex is minus the partial derivative of $E[\mbf{X}_1, \dots, \mbf{X}_M]$ with respect to $\mbf{X}_k$,
\begin{equation}
 \mbf{F}_k \Delta \mbf{s} = - \frac{\partial E}{\partial \mbf{X}_k} = - \gamma \sum_{l \, \in \nbor(k)} \frac{\partial |\mbf{T}_l| }{\partial \mbf{X}_k} \, , \label{lagrangianforce}
\end{equation}
where $\nbor(k)$ denotes the set of indices of triangles that share $\mbf{X}_k$ as a vertex\footnote{Here $\Delta \mbf{s}$ is the Lagrangian area associated with each node and $\mbf{F}_k$ is the Lagrangian force density with respect to Lagrangian area, but note that we do not need $\mbf{F}_k$ and $\Delta \mbf{s}$ separately; only their product is used in the numerical scheme.}. Each component of the term $\partial |\mbf{T}_l| / \partial \mbf{X}_k$ in \Cref{lagrangianforce} can be computed analytically \cite{Kim2014_foam3D}, 
\begin{align}
  \left(\frac{\partial |\mbf{T}_l|}{\partial \mbf{X}_k} \right)_{\alpha} &=  \frac{\partial}{ \partial \mbf{X}_{k,\alpha}} \left( \frac{1}{2} \left|(\mbf{X}_k - \mbf{X'}_k) \times (\mbf{X'}_k - \mbf{X}^{''}_k) \right| \right) \nonumber \\
    &= \frac{1}{2} \left(  (\mbf{X}^{'}_k - \mbf{X}^{''}_k) \times \mbf{\hat{n}}_l \right)_{\alpha} \ , \quad \alpha = 1,2,3,
\end{align}
where $\mbf{X}_k,\ \mbf{X}^{'}_k,\ \mbf{X}^{''}_k$ denote the three vertices of the triangle $\mbf{T}_{l}$ ordered in the counterclockwise direction and $\mbf{\hat{n}}$ is the unit outward normal vector of $\mbf{T}_{l}$.

\begin{figure}
\centering
\subfloat[][$t=0$]{\includegraphics[width = .5\linewidth]{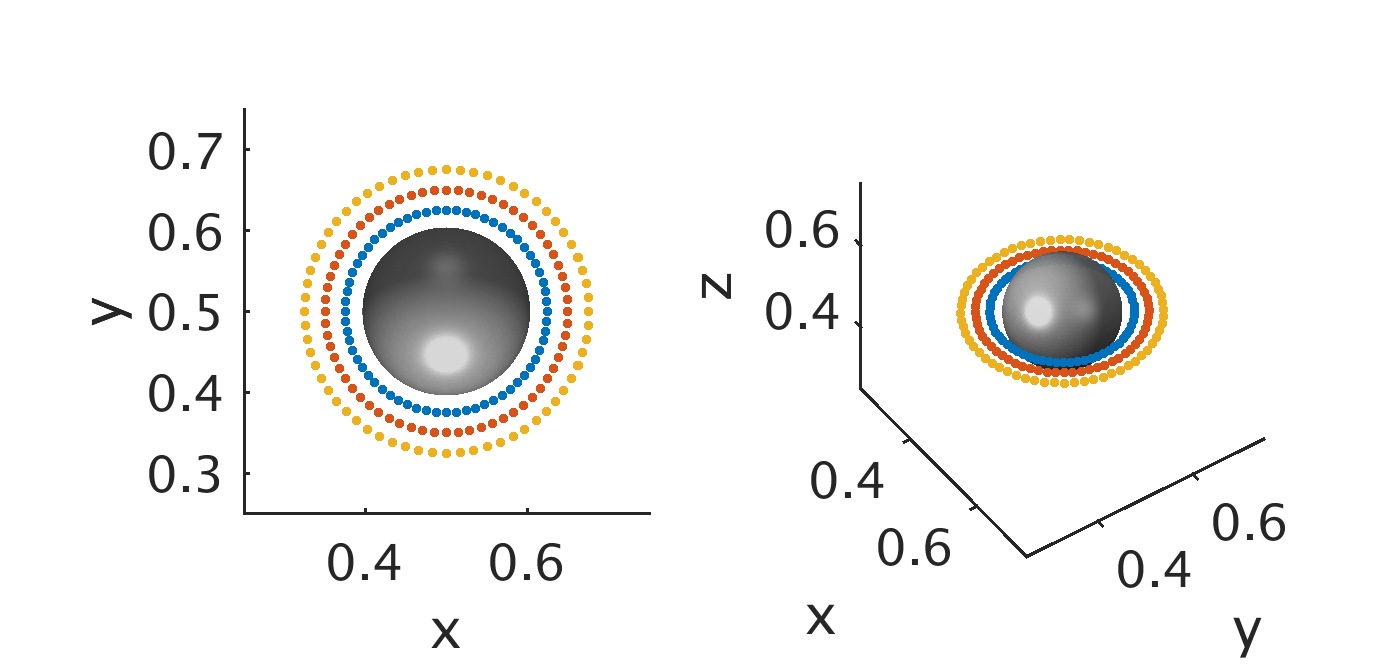}} \hspace{-1em}
\subfloat[][$t= \frac{1}{32}$]{\includegraphics[width = .5\linewidth]{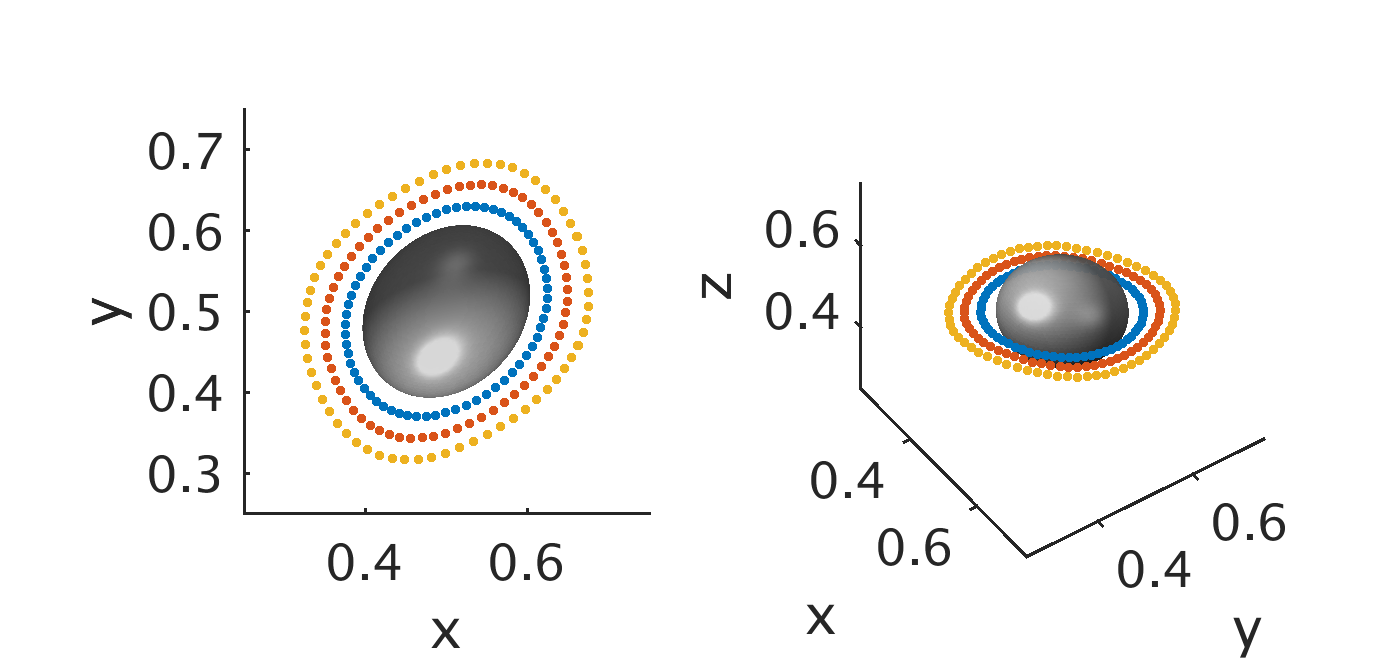}} \\
\subfloat[][$t= \frac{1}{4}$]{\includegraphics[width = .5\linewidth]{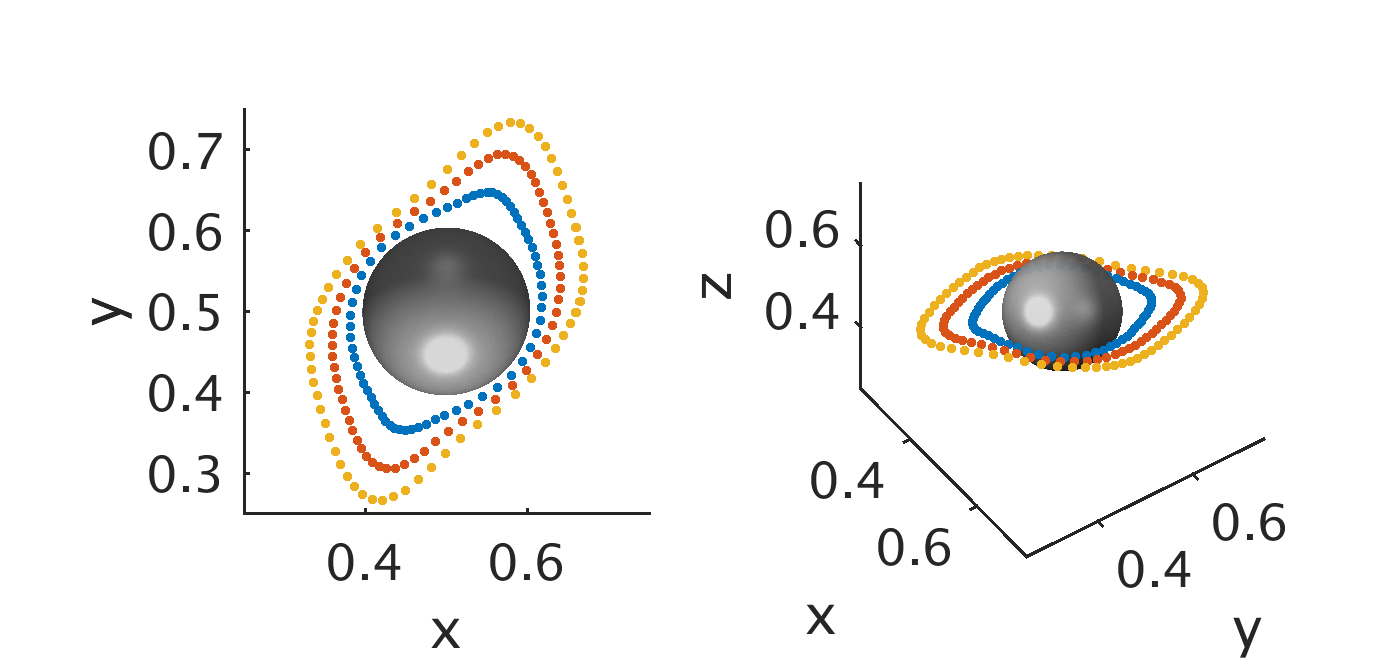} } \hspace{-1em}
\subfloat[][$t= \frac{1}{2}$]{\includegraphics[width = .5\linewidth]{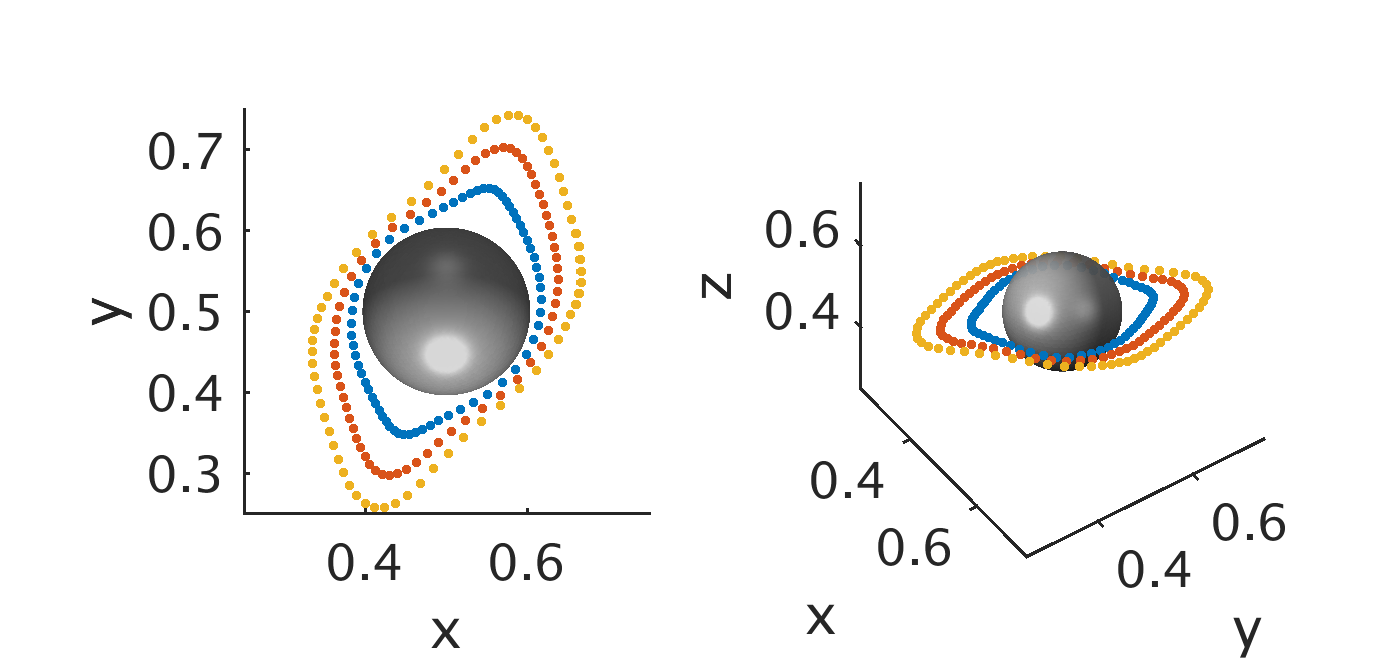} } 
\caption{Deformation of a 3D elastic membrane immersed in a viscous fluid with initial velocity $\mbfit{u}(\mbfit{x},t) = (0, \ \sin(4\pi x), \ 0)$ at $t = 0, \frac{1}{32},  \frac{1}{4}$ and $\frac{1}{2}$. The computation is performed using DFIB with $\newsix$  in the periodic box $\Omega = [0,1]^3$ with Eulerian meshwidth $h = \frac{1}{128}$. The elastic membrane, initially in spherical configuration with radius $R \approx 0.1$, is discretized by a triangular surface mesh with $M=10242$ vertices and $P=20480$ facets so that $h_s = \frac{h}{2}$ in the initial configuration.  Colored markers that move passively with the divergence-free interpolated fluid velocity  are added for visualizing the fluid flow in the vicinity of the membrane interface.}
\label{fig:surfsphere3D}
\end{figure}

The computation is performed in the periodic box $\Omega = [0,1]^3$ with Eulerian meshwidth $h = \frac{1}{128}$ using DFIB with $\newsix$. For the quasi-static test, the initial fluid velocity is set to be zero, and for the dynamic test, we set $\mbfit{u}(\mbfit{x},0) = (0, \ \sin(4\pi x), \ 0 )$. In the computational results shown in \autoref{fig:surfsphere3D}, the spherical membrane is discretized by triangulation (as shown in \autoref{fig:sphere_refine}) with 5 successive levels of refinement from the regular icosahedron (\autoref{fig:sphere0}), which results in a triangular mesh with $M=10242$ vertices and $P=20480$ facets. The radius of the spherical membrane is set to be $R \approx 0.1$ which corresponds to $h_s \approx \frac{h}{2}$. 
The remaining parameters in the computation are $\rho = 1, \ \mu = 0.05, \ \gamma = 1$ and the time step size $\Delta t = \frac{h}{4}$. In
\autoref{fig:surfsphere3D} we show snapshots of the 3D elastic membrane  
at $t=0, \ \frac{1}{32}, \ \frac{1}{4}$ and $\frac{1}{2}$ for the dynamic case. The elastic interface is instantaneously deformed by the fluid flow in the $y$-direction, and due to surface tension, the membrane eventually relaxes back to the spherical equilibrium configuration. Colored markers that move passively with the divergence-free interpolated fluid velocity are added for visualizing the fluid flow in the vicinity of the interface. 

\begin{figure}
\centering
\subfloat[][Quasi-static test]{\includegraphics[width = 0.48\linewidth]{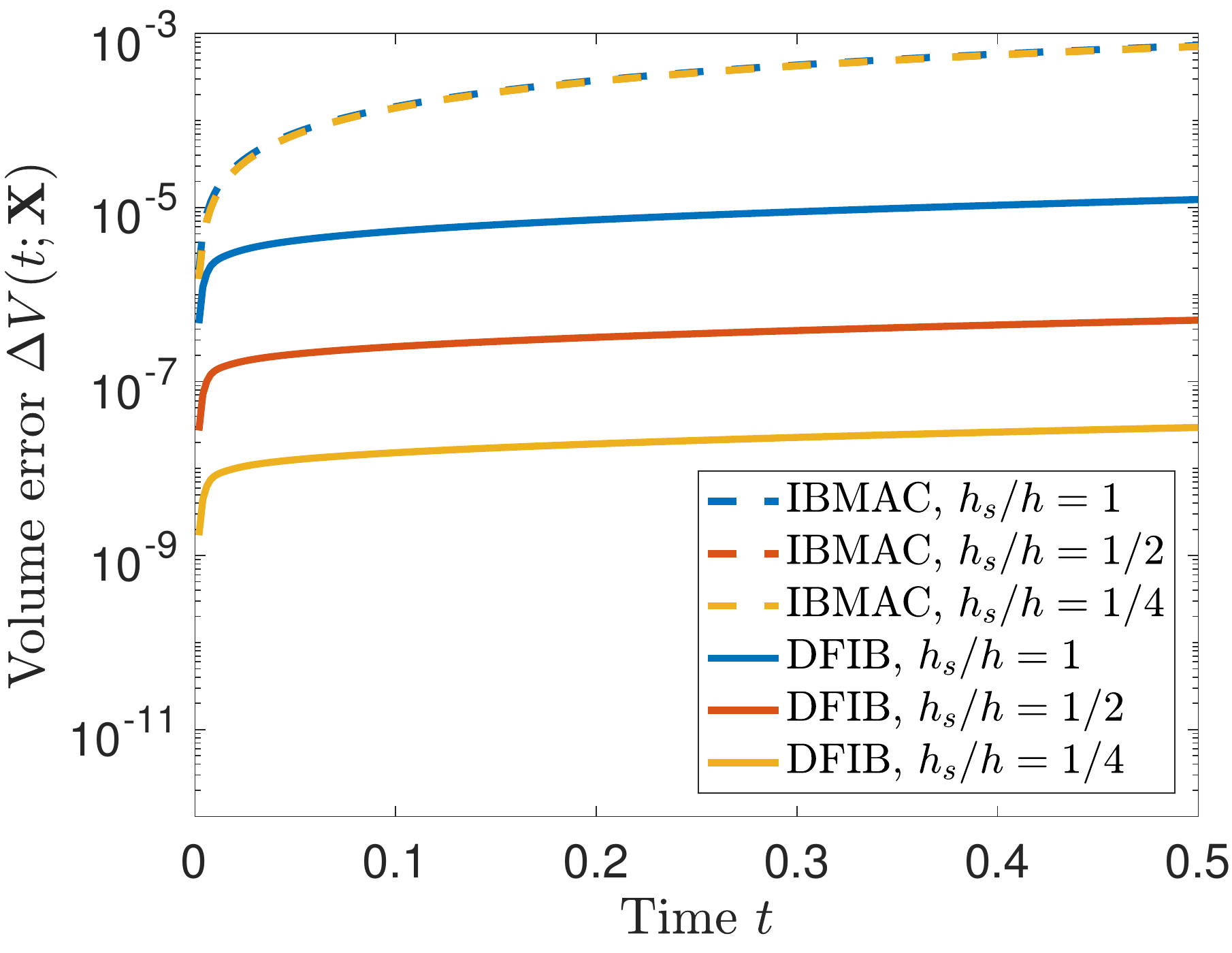} \hspace{1em}  \label{fig:volumeloss3D_zero} }
\subfloat[][Dynamic test]{\includegraphics[width = 0.48\linewidth]{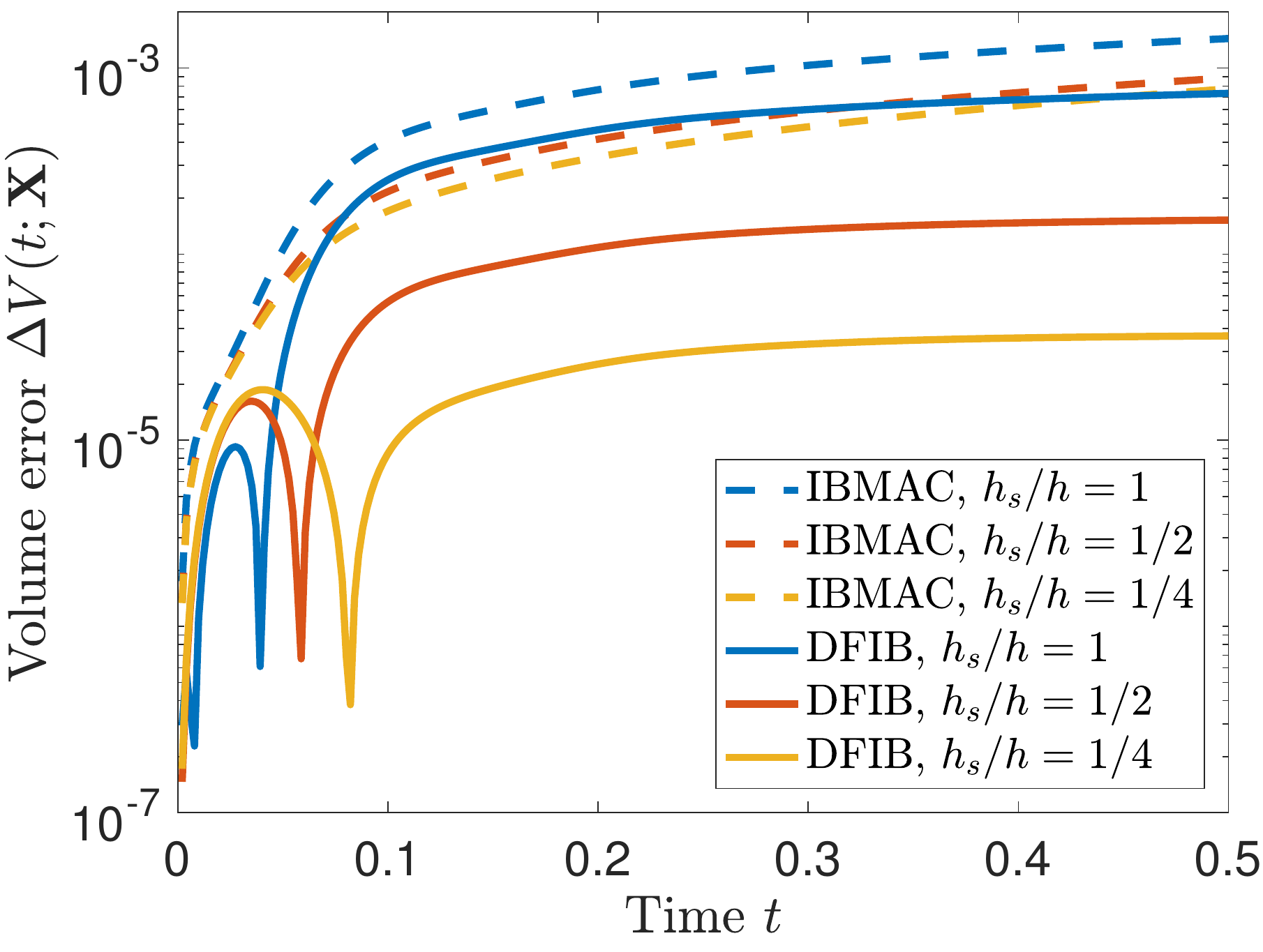} \label{fig:volumeloss3D_shear} }
\caption{Normalized volume error $\Delta V(t ; \mbf{X})$ of a 3D elastic membrane using IBMAC and DFIB with $h_s = h, \frac{h}{2}, \frac{h}{4}$, where $h=\frac{1}{128}$. For (a) the quasi-static test, $\Delta V(t,\mbf{X})$ of DFIB decreases with mesh refinement, while there is no improvement in volume error for IBMAC. For (b) the dynamic test, the volume error in DFIB remains (almost) steady in time for $h_s  = \frac{h}{2}, \frac{h}{4}$ as the membrane rests,  whereas we see no substantial improvement in volume conservation with mesh refinement for IBMAC, and the volume loss keeps increasing in time. For this set of computations, the $\mathscr{C}^3$ 6-point kernel $\newsix$ is used.}
\label{fig:volumelossDyn3D}
\end{figure}

The volume enclosed by the triangular surface mesh is approximated by the total volume of tetrahedra formed by each facet and one common reference point (e.g. the origin) using the scalar triple product. 
To study volume conservation of the DFIB method in 3D, we compare the normalized volume error defined by
\begin{equation}
 \Delta V(t ; \mbf{X}) := \frac{|\vol(t;\mbf{X}) - \vol(0;\mbf{X})|}{\vol(0;\mbf{X})} \label{volume_err}
\end{equation}
using IBMAC and DFIB with $h_s = h, \frac{h}{2}, \frac{h}{4}$, which correspond to triangular meshes with 4,5,6 levels of refinement from the regular icosahedron respectively. 
For the quasi-static case (\autoref{fig:volumeloss3D_zero}), volume errors for DFIB are at least 2 orders of magnitude smaller than those of IBMAC. Further, volume errors for DFIB keep decreasing as the Lagrangian mesh is refined from $h_s = h$ to $\frac{h}{4}$. For the dynamic case (\autoref{fig:volumeloss3D_shear}), both methods suffer a significant amount of volume loss arising from the rapid deformation at the beginning of simulation. The volume error of DFIB with $h_s = h$ is similar to those of IBMAC in magnitude, {but the volume error of DFIB  decreases as the Lagrangian mesh is refined for $h_s = \frac{h}{2}, \frac{h}{4}$.} 
It appears that the behavior of volume error changes in nature from $h_s = h$ to $\frac{h}{2}$, which coincides with the conventional recommendation that the best choice of Lagrangian mesh spacing in the IB method is $h_s = \frac{h}{2}$ in practice. 
{Similar to the two sources of error that contribute to the area loss in 2D, the volume error observed in \autoref{fig:volumelossDyn3D} can also be explained by contribution from the time-stepping error, and the volume loss due to only moving the vertices (Lagrangian markers) that constitute the triangular mesh. This kind of error in volume conservation decreases as the discretization of the surface is refined, even on a fixed Eulerian grid (as shown in \autoref{fig:volumelossDyn3D}).}
Finally, we remark that the improvement in volume conservation does not seem to be as substantial as the improvement in area conservation in 2D. We suspect that this may be attributed to the larger approximation error in computing the volume using the tetrahedral approximation (after the triangular mesh is deformed), whereas in two dimensions {we use a higher-order representation of the interface (cubic splines).} 
Nevertheless, the reduction in volume error from the Lagrangian mesh-refinement experiments indeed confirms that the DFIB method can generally achieve better volume conservation if the immersed boundary is sufficiently resolved ($h_s \leq \frac{h}{2}$).


\section{Conclusions}
In this paper, we introduce an IB method with divergence-free velocity interpolation and force spreading. Our IB method makes use of staggered-grid discretization to define an edge-centered discrete vector potential.
By interpolating the discrete vector potential in the conventional IB fashion, we obtain a continuum vector potential whose curl directly yields a continuum Lagrangian velocity field that is exactly divergence-free by default. The corresponding force-spreading operator is constructed to be the adjoint of velocity interpolation so that energy is preserved in the interaction between the fluid and the immersed boundary. Both the new interpolation and spreading schemes require solutions of discrete vector Poisson equations which can be efficiently solved  by a variety of algorithms. The transfer of information from the Eulerian grid to the Lagrangian mesh (and vice versa) is performed using $\grad \delta_h$ on the edge-centered staggered grid $\edgegrid$. We have found that volume conservation of DFIB improves with the smoothness of the IB kernel used to construct $\delta_h$, and  we have numerically tested that IB kernels that are at least $\mathscr{C}^2$ are good candidate kernels that can be used to construct the regularized delta function in the DFIB method.

We have incorporated the divergence-free interpolation and spreading operators in a  second-order time-stepping scheme, and applied it to several benchmark problems in two and three spatial dimensions. First, we have tested that our method achieves second-order convergence in both the fluid velocity and the Lagrangian deformation map for the 2D surface tension problem, which is admittedly a special case, since its continuum solution has a continuous normal derivative of the tangential velocity across the immersed boundary. The highlight of the DFIB is its capability of substantially reducing volume error in the immersed structure as it moves and deforms in the process of fluid-structure interaction. Through numerical simulations of quasi-static and dynamic membranes, we have confirmed that the DFIB method improves volume conservation by several orders of magnitude compared to IBMAC and IBModified. Furthermore, owing to the divergence-free nature of its velocity interpolation, the DFIB method reduces volume error with Lagrangian mesh refinement while keeping the Eulerian grid fixed. A similar refinement study would not yield improved volume conservation when using the conventional IB method. {Although the numerical examples considered in this paper only involve thin elastic structures, we note that the DFIB method can also be directly applied to model thick elastic structures without any modification to the method, other than representing the thick elastic structure by a curvilinear mesh of Lagrangian points \cite{Griffith2005_IBaccuracy,Image_IB-FEM_2014,Hybrid_FDFEM_2017,HyperelasticIB}.  We also remark that the current version of the DFIB method is accompanied with a single-fluid Navier-Stokes fluid solver. This is not a fundamental limitation in our approach, and as a direction of future research, the DFIB method may be extended to work with variable-viscosity and variable-density fluid solvers \cite{Fai2013}. }

Unlike other improved IB methods that either use non-standard finite-difference operators (IBModified \cite{Peskin1993_IBmodified}) that complicate the implementation of the fluid solver, or rely on analytically-computed correction terms (IIM \cite{Li2001, Lee2003_IIM} or Blob-Projection method \cite{Cortez2000_blobprojection}) that may not be readily accessible in many applications, the DFIB method is generally applicable, and it is straightforward to implement in both 2D and 3D from an existing IB code that is based on the staggered-grid discretization. Moreover, the additional costs of performing the new interpolation and spreading do not increase the overall complexity of computation and are modest compared to the existing IB methods.

We point out two limitations of our present work.
A first limitation of the current version of DFIB method is based on the assumption of periodic boundary conditions. Extending the method to include physical boundary conditions at the boundaries of the computational domain is one possible direction of future work, {but there are several challenges to overcome. First,  a special treatment of spreading and interpolation is required near the boundaries since the support of the IB kernel can extend outside of the physical domain \cite{IBMDelta_Boundary,RigidIBM}.
Second, with non-periodic BCs, instead of FFTs, the resulting linear system needs to be solved by geometric or algebraic multigrid method to achieve high performance. For unbounded domains, a lattice Green's function technique was recently proposed as an alternative approach \cite{LiskaColonius2017}.
Third, for domains with physical boundaries, the use of projection-based fluid solvers to eliminate pressure introduces splitting errors near the physical boundaries, and instead one ought to solve a coupled velocity-pressure system at every time step \cite{NonProjection_Griffith}.    In the DFIB method, it is also nontrivial to specify boundary conditions for the vector potential $\mbf{a}(\mbf{x})$ at physical boundaries, which would lead to a Poisson equation with non-periodic BCs. It may be that volume conservation as the structure passes near a boundary requires solving a coupled velocity-potential system. A second limitation of our DFIB method is that the pressure gradient generated by the Lagrangian forces is part of the resulting Eulerian force density because force spreading is also constructed to be discretely divergence-free. However, this may also be an important advantage of our method from the standpoint of accuracy, since it means that jumps in pressure across the interface do not require any explicit representation. We do not yet see an obvious way to extract the pressure from the Eulerian force density in case it is needed for output purposes, or for the purposes of imposing physical boundary conditions involving tractions or avoiding splitting errors near boundaries \cite{NonProjection_Griffith}.}

\section*{Acknowledgement}
Y. Bao and A. Donev were supported in part by the National Science Foundation under award DMS-1418706, and by the U.S. Department of Energy Office of Science, Office of Advanced Scientific Computing Research, Applied Mathematics program under Award Number DE-SC0008271. B.E. Griffith acknowledges research support from the National Science Foundation (NSF awards ACI 1450327, DMS 1410873, and CBET 1511427) and the National Institutes of Health (NIH award HL117063).

\appendix


\section{Vector identities of discrete differential operators}
\label{appendix_vector_identity}
Suppose $\vfun{\varphi}{\mbf{x}}$ is a scalar grid function defined on $\cellgrid$, and $\vfun{u}{x}$ and $\vfun{a}{x}$ are  vector grid functions defined on $\facegrid$ and $\edgegrid$ respectively. The following discrete vector identities are valid on the periodic staggered grid just as in the continuum case,
\begin{align}
 & \dcurl{h} \dgrad{h} \varphi  = 0, \label{curlgrad} \\
 & \ddiv{h} ( \dcurl{h} \mbf{u} ) = 0, \label{divcurl} \\
 & \dcurl{h} ( \dcurl{h} \mbf{u} ) = \dgrad{h}(\ddiv{h} \mbf{u}) - \vlapl{h} \mbf{u} \label{curlcurl}, \\
 & \sum_{\mbf{x} \in \facegrid} \vfun{u}{x} \cdot (\dgrad{h}\varphi) (\mbf{x}) h^3 = -\sum_{\mbf{x} \in \cellgrid} 
(\ddiv{h} \mbf{u})(\mbf{x}) \, \varphi(\mbf{x}) h^3, \label{sumbyparts1} \\
 & \sum_{\mbf{x} \in \edgegrid} \vfun{a}{x} \cdot (\dcurl{h} \mbf{u})(\mbf{x}) h^3 =  \sum_{\mbf{x} \in \facegrid} (\dcurl{h} \mbf{a})(\mbf{x}) \cdot \vfun{u}{x}  h^3. \label{sumbyparts2}
\end{align}
\Cref{curlgrad,curlcurl} are merely discrete analogues of well-known vector identities involving gradient, divergence and curl. These identities can be proved in the same manner as their continuous counterparts. \Cref{sumbyparts1,sumbyparts2} can be verified via ``summation by parts''. Note that \Cref{divcurl,curlcurl} also hold if we replace $\mbf{u}$ (which lives on $\facegrid$) by $\mbf{a}$ (which lives on $\edgegrid$).

\section{Existence of discrete vector potential}
\label{appendix_existence}

\begin{lem}
 Suppose $\ddiv{h} \mbf{u} = 0$ and $\dcurl{h} \mbf{u} = 0$ for $\mbf{x} \in \facegrid$, then $\mbf{u}(\mbf{x})$ is a constant function on $\facegrid$. \label{lemma1}
\end{lem}

\begin{pf}
 To prove this statement, we use  \Cref{curlcurl,sumbyparts1,sumbyparts2},
 \begin{align*}
  \sum_{\mbf{x} \in \edgegrid} (\dcurl{h} \mbf{u})(\mbf{x}) \cdot (\dcurl{h} \mbf{u}) (\mbf{x}) h^3 &= \sum_{\mbf{x} \in \facegrid} \mbf{u}(\mbf{x}) \cdot (\dcurl{h} (\dcurl{h} \mbf{u}) ) h^3 \\
  &= \sum_{\mbf{x} \in \facegrid} \mbf{u}(\mbf{x}) \cdot \dgrad{h} (\ddiv{h} \mbf{u}) h^3 - \sum_{\mbf{x} \in \facegrid} \mbf{u}(\mbf{x}) \cdot (\vlapl{h} \mbf{u}) h^3 \\
  &= -\sum_{\mbf{x} \in \cellgrid} ( \ddiv{h} \mbf{u} )^2 h^3 + \sum_{\substack{ \mbf{x} \in \edgegrid, i\neq j \\ \mbf{x} \in \cellgrid, i = j}} \left( D^h_j u_i \right)^2 h^3.
 \end{align*}
Thus,
\begin{equation}
 \sum_{ \substack{ \mbf{x} \in \edgegrid, i\neq j \\ \mbf{x} \in \cellgrid, i = j}} \left(D^h_j u_i\right)^2 h^3 = \sum_{\mbf{x} \in \edgegrid } \left|(\dcurl{h} \mbf{u})\right|^2 h^3 + \sum_{\mbf{x} \in \edgegrid} ( \ddiv{h} \mbf{u} )^2 h^3. \label{lemma1_eq4}
\end{equation}
Since $\ddiv{h} \mbf{u} = 0$ and $\dcurl{h} \mbf{u} = 0$ by hypothesis, the left-hand side of \Cref{lemma1_eq4} is also zero. But this implies $u_i$ is constant for $i=1,2,3$. 
\end{pf}
\begin{lem}
 If $\psi$ is a scalar grid function that lives on one of the staggered grids, such that
 \begin{equation}
  \sum_{\mbf{x}} \psi({\mbf{x}}) h^3 = 0, \label{psi_zero}
 \end{equation}
 then there exists a grid function $\varphi$ such that
 \begin{equation}
  \dlapl{h} \varphi = \psi. \label{appendix_poisson}
 \end{equation}
 \label{lemma2}
\end{lem}
\begin{pf}
This lemma states  the solvability of the discrete Poisson problem \Cref{appendix_poisson}. Since $\dlapl{h} = \ddiv{h} \dgrad{h}$ is symmetric with respect to the inner product on the periodic grid
\begin{equation}
	(\varphi, \psi) = \sum_{\mbf{x}} \varphi(\mbf{x}) \psi(\mbf{x}) h^3, 
\end{equation}
what we have to show is that any $\psi$ satisfying \Cref{psi_zero} is orthogonal to any $\varphi_0$ in the null space of $\vlapl{h}$. But the null space of $\vlapl{h}$ with periodic boundary conditions contains only the constant function, and hence $(\psi , \varphi_0) = 0$ because of \Cref{psi_zero} as required. 
\end{pf}

Now we are ready to state the theorem that guarantees the existence of a discrete vector potential $\vfun{a}{x}$ for $\mbf{x} \in \edgegrid$ given a discretely divergence-free velocity field $\vfun{u}{x}$ for $\mbf{x} \in \facegrid$.
\begin{thm} \label{appendix_vpot_existence}
Suppose $\vfun{u}{x}$ is a periodic grid function for $\mbf{x} \in \facegrid$, and $\vfun{u}{x}$ satisfies 
\begin{equation}
\sum_{\mbf{x} \in \facegrid} \vfun{u}{x} h^3 = 0 \quad \text{and} \quad \ddiv{h} \mbf{u} = 0,
\end{equation}
then there exists a grid function $\vfun{a}{x}$ for $\mbf{x} \in \edgegrid$ such that 
\begin{equation}
\mbf{u} = \dcurl{h} \mbf{a}.
\end{equation}
\end{thm}

\begin{pf}
We choose $\vfun{a}{x}$ to be any solution of 
\begin{equation}
-\vlapl{h} \, \mbf{a} = \dcurl{h} \mbf{u}. \label{thm1_poisson}
\end{equation}
Such an $\vfun{a}{x}$ exists by Lemma \ref{lemma2}, because 
\begin{align*}
\sum_{\mbf{x} \in \edgegrid} \left( \dcurl{h} \mbf{u} \right)_i(\mbf{x}) &= \epsilon_{ijk} \sum_{\mbf{x} \in \edgegrid} 1 \cdot D_j u_k  \, h^3  \\
&= -\epsilon_{ijk} \sum_{\mbf{x} \in \facegrid} (D_j 1) u_k \, h^3 \\
&= 0.
\end{align*}
By applying $\ddiv{h}$ to \Cref{thm1_poisson} and using the property that $\vlapl{h}$ and  $\ddiv{h}$ commute, we also have 
\begin{equation}
-  \vlapl{h} (\ddiv{h} \mbf{a}) = \ddiv{h} (\dcurl{h} \mbf{u}) = 0. \label{thm1_nullspace}
\end{equation}
Because the null space of $\vlapl{h}$ contains only the constant function, it follows that 
\begin{equation}
\dgrad{h} (\ddiv{h} \mbf{a}) = 0. \label{thm1_dgrad_ddiv_a}
\end{equation}
If we use \Cref{thm1_dgrad_ddiv_a} and \Cref{curlcurl} for $\mbf{a}$ , we can rewrite \Cref{thm1_poisson} as 
\begin{equation}
	\dcurl{h} (\dcurl{h} \mbf{a}) = \dcurl{h} \mbf{u}, 
\end{equation}
or 
\begin{equation}
	\dcurl{h} (\dcurl{h} \mbf{a} - \mbf{u}) = 0. \label{thm1_dcurl}
\end{equation}
But we also know from \Cref{divcurl} and the requirement that $\ddiv{h} \mbf{u} = 0$ that 
\begin{equation}
	\ddiv{h} (\dcurl{h} \mbf{a} - \mbf{u}) = 0.  \label{thm1_ddiv}
\end{equation} 
From \Cref{thm1_dcurl,thm1_ddiv} and Lemma \ref{lemma1}, it follows that
\begin{equation}
 \dcurl{h} \mbf{a} - \mbf{u} = \text{constant}.
\end{equation}
The constant must be zero, however, since $\dcurl{h} \mbf{a}$ has zero sum by ``summation by parts'', and $\mbf{u}$ has zero sum by assumption. This completes the proof of the existence of a vector potential satisfying $\mbf{u} = \dcurl{h} \mbf{a}$.
\end{pf}

\bibliographystyle{elsart-num-sort}
\bibliography{References,DFIB}

\end{document}